\documentclass[10pt]{amsart}
\usepackage{amsthm}
\usepackage{amsmath,amsfonts,amssymb,mathrsfs}
\usepackage[dvipsnames]{xcolor}
\usepackage{geometry}
\usepackage{fancyhdr}
\usepackage[all]{xy}
\usepackage{graphicx}

\def\sqtimes{\setbox0=\hbox{\kern-.13em$\times$\kern-.13em}
     \dimen0=\ht0 \advance\dimen0 -.09em \ht0=\dimen0
     \dimen0=\dp0 \advance\dimen0 -.09em \dp0=\dimen0
     \mathbin{\vcenter{\hrule\kern-.4pt
       \hbox{\vrule\kern-.4pt$\box0$\kern-.4pt\vrule}\kern-.4pt\hrule}}}

\def\e{\varepsilon}

\def\ML{\mathscr{M\! L}}
\def\Map{\operatorname{Mod}}

\def\supp{\operatorname{supp}}

\def\stab{\mathrm{stab}}
\def\fill{^{\mathrm{bind}}}
\def\fh{^{\mathrm{fh}}}
\def\Curr{\mathscr{C}}
\def\CurrK{\Curr^{\leq 1}}

\def\Cufill{\Curr\fill}
\def\Cufh{\Curr\fh}
\def\Cutriple{\Curr\fh}

\def\Thu{m_{\mathrm{Th}}}
\def\bThu{\ol{m}_{\mathrm{Th}}}

\def\sys{\mathrm{sys}}

\def\MLfh{\ML\fh}
\def\RR{\mathbb{R}}
\def\Geod{\mathcal{G}}

\def\GeodK{\Geod^{\leq 1}}

\def\IG{\mathcal{IG}}
\def\S{S}
\def\R{R}
\def\A{A}
\def\a{\alpha}
\def\ALPH{\eta}
\def\Arcs{\mathfrak{A}}

\def\K{\mathcal{K}}

\def\dis{\displaystyle}

\def\DD{\mathbb{D}}

\def\cpt{M}

\def\inte#1{\mathring{#1}}

\def\hull{\mathrm{hull}}
\def\PP{\mathbb{P}}

\def\pa{\partial}

\def\ol#1{\overline{#1}}

\def\wti#1{\widetilde{#1}}
\def\ti#1{\tilde{#1}}

\def\rar{\rightarrow}
\def\lra{\longrightarrow}

\def\Zcal{\mathcal{Z}}
\def\arcsinh{\mathrm{arcsinh}}

\def\Diff{\mathrm{Diff}}
\def\ZZ{\mathbb{Z}}

\def\Teich{\mathcal{T}}
\def\Liou{\mathcal{L}}

\def\NN{\mathbb{N}}
\def\Coll{\mathrm{Col}}

\theoremstyle{definition}

\newtheorem{thm}{Theorem}[section]
\newtheorem{lemma}[thm]{Lemma}
\newtheorem{prop}[thm]{Proposition}
\newtheorem{cor}[thm]{Corollary}

\newtheorem{example}[thm]{Example}
\newtheorem{remark}[thm]{Remark}

\newtheorem{defi}[thm]{Definition}

\newtheorem{maintheorem}{Theorem}

\newtheorem{mainprop}[maintheorem]{Proposition}
\newtheorem{maincor}[maintheorem]{Corollary}
\newtheorem{mainlemma}[maintheorem]{Lemma}
\newtheorem*{notation*}{Notation}

\newtheorem*{namedtheorem}{\theoremname}
\newcommand{\theoremname}{testing}
\newenvironment{named}[1]{\renewcommand{\theoremname}{#1}\begin{namedtheorem}}{\end{namedtheorem}}
\numberwithin{equation}{section}

\begin{document}

\title[]{Ergodic invariant measures \\ on the space of geodesic currents}

\author{Viveka Erlandsson}
\address{Viveka Erlandsson - School of Mathematics, University of Bristol and Department of Mathematics and Statistics UiT Arctic University of Norway}
\email{v.erlandsson@bristol.ac.uk}

\author{Gabriele Mondello}
\address{Gabriele Mondello - Department of Mathematics, ``Sapienza'' Universit{\`a} di Roma}
\email{mondello@mat.uniroma1.it}

\begin{abstract}
Let $\S$ be a compact, connected, oriented surface, possibly with boundary,  of negative Euler characteristic.
In this article we extend Lindenstrauss-Mirzakhani's and Hamenst{\"a}dt's classification of 
locally finite mapping class group invariant ergodic measures
on the space of measured laminations $\ML(\S)$ to the space of geodesic currents $\Curr(\S)$,
and we discuss the homogeneous case.
Moreover, we extend
Lindenstrauss-Mirzakhani's classification of orbit closures to $\Curr(\S)$. 
Our argument relies on their results and
on the
decomposition of a current into a sum of three currents with
isotopically disjoint supports: a measured lamination without closed leaves, a simple multi-curve and a current that binds its hull. 
\end{abstract}



\maketitle


%
    


\section{Introduction}\label{sec:intro}

%

\subsection{Setting}
Let $\S $ be a smooth, compact, connected, oriented surface of negative Euler characteristic, possibly
with boundary, and let $\Map(\S )$ be its mapping class group, i.e.~the group of isotopy classes of orientation-preserving diffeomorphisms $\S \to\S $ that send each boundary curve of $\S$ to itself.

Consider an auxiliary hyperbolic metric on $\S$ such that $\pa\S$ is geodesic.
A geodesic current on $\S$ is a $\pi_1(\S )$-invariant Radon measure on the space $\Geod(\tilde{\S })$ of bi-infinite geodesics in the universal cover $\tilde{\S}$ of $\S$. The space $\Curr(\S )$ of all geodesic currents, 
naturally endowed with the weak$^\star$-topology, can also be viewed as the completion of the set of weighted closed curves on $\S$ in the same way as the space  $\ML(\S )$ of measured laminations is the completion of the set of weighted simple closed curves. Recall that a measured lamination is a closed subset of $\S$ foliated by complete geodesics and endowed
with a transverse measure of full support. Hence a measured lamination can be viewed as a current and $\ML(\S )$ can be viewed as a subspace of $\Curr(\S )$. The geometric intersection number of closed curves has a unique continuous extension to a symmetric, bi-homogenous intersection form $\iota(\cdot, \cdot) : \Curr(\S )\times\Curr(\S )\to\RR_{\geq 0}$ (see \cite{bonahon:currents}). The subspace of measured laminations consists exactly of those currents $c$ for which $\iota(c,c)=0$.

The aim of this paper is to provide a classification
of locally finite ergodic measures on $\Curr(\S)$ that are invariant under the natural
action of $\Map(\S)$ and of closures of $\Map(\S)$-orbits on $\Curr(\S)$.

\subsection{Motivation}
The impetus for the present paper -- in addition to the classification theorem for ergodic invariant measures on $\ML(\S)$ proven in
Lindenstrauss-Mirzakhani \cite{LM08} (and almost completely in Hamenst{\"a}dt \cite{hamenstadt:measures}) -- was a series of articles on counting problems 
of closed curves and of currents on surfaces, originating with Mirzakhani \cite{LM08, mirzakhani-orbit} in the hyperbolic case
and generalized to other settings by Erlandsson-Souto \cite{ES}, Erlandsson-Parlier-Souto \cite{EPS} and Rafi-Souto \cite{rafi-souto}.
Part of Mirzakhani's argument
was to study a sequence of measures on $\ML(\S)$ converging to
a multiple of the Thurston measure $\Thu$. The main ingredient
in its generalizations was to analyze the corresponding sequence of measures on $\Curr(\S)$ and show that they in fact limit to a homogeneous measure supported on $\ML(\S)$.
In fact, we recover the result in \cite{ES} as a consequence of our classification of
invariant ergodic homogeneous measures on $\Curr(\S)$.

\subsection{Notational conventions}\label{sec:conventions}

All surfaces we consider are smooth, compact, oriented and possibly with boundary; we also require that each component of a surface has negative Euler characteristic.
Moreover, a subsurface of a surface is always meant to be closed and with smooth boundary.

As a rule, the surface $\S$ is connected, while a subsurface $\R$ of $\S$ may be disconnected, unless differently specified. For such an $\R$ with connected components $\{\R_i\}$
we define its space of geodesic currents as the product $\Curr(\R):=\prod_i \Curr(\R_i)$ and
its mapping class group as $\Map(\R):=\prod_i\Map(\R_i)$.
Many results we are going to state for connected surfaces can be easily extended to disconnected ones in an obvious way. 
We will occasionally stress in the hypothesis that $\S$ is connected when it is particularly relevant.

Throughout the paper, subsurfaces, simple closed curves, laminations and supports of currents will often be
considered {\it{up to isotopy}}. So, for instance, we
will say that the subsets $\{X_k\}$ of $\S$ are {\it{isotopically disjoint}} if there exist isotopies $f_k:\S\rar\S$
such that the subsets $\{f_k(X_k)\}$ are pairwise disjoint. 
As another example, if $c$ is a current and $h$ is a hyperbolic metric on $\S$, then
$\supp_h(c)$ is the union of all $h$-geodesics in $\S$ in the support of $c$; but
the {\it{support}} $\supp(c)$ is the isotopy class of $\supp_h(c)$, which is independent of the choice of $h$.

Simple-closed-curve-free (scc-free) currents and a-laminational currents defined below will be
particularly important in the formulation of our results.

\begin{defi}[Scc-free and a-laminational currents]\label{def:scc-free}
A geodesic current $c$ on $\S$ is {\it{scc-free}} if it cannot be written as a sum $c=\Gamma+c'$
of two currents with isotopically disjoint supports, where $\Gamma\neq 0$ is a weighted simple multi-curve.
Such $c$ is {\it{a-laminational}} if it is scc-free and no connected component of $\supp(c)$ is a lamination.
\end{defi}


Finally, we also introduce the hull of a current.

\begin{defi}[Hull of a current]\label{def:hull}
The \emph{surface hull} of an scc-free current $\check{c}\in\Curr(\S)$ 
is the isotopy class $\hull(\check{c})$ of the smallest closed subsurface of $\S$ that contains the support of $\check{c}$.
\end{defi}

In Section \ref{sec:hull} we will see that the hull is well-defined
and we will discuss some of its properties.

\subsection{Invariant measures and orbit closures in $\ML$}
Lindenstrauss-Mirzakhani \cite{LM08} and (almost completely) Hamenst{\"a}dt \cite{hamenstadt:measures} independently classified all $\Map(\S )$-invariant, locally finite, ergodic measures on $\ML(\S)$. Here we describe such classification
and we adopt the terminology used in \cite{LM08}, since this is more in alignment with our result. 

Their main theorem states that any such measure $m$ is a positive multiple of a measure associated to a so-called complete pair $(\R,\Gamma)$ (see Theorem \ref{thm:LM} and \cite[Theorem 1.1]{LM08}). 
Here, a complete pair $(\R, \Gamma)$ consists of a simple multi-curve $\Gamma$ and a subsurface $\R\subset \S$ such that 
$\R$ and $\supp(\Gamma)$ are isotopically disjoint, and each boundary curve of $\R$ is 
homotopic either to a curve in the support of $\Gamma$ or to a boundary curve on $\S$. 
Viewing the space $\ML_0(\R)$ of measured laminations supported in the interior of $\R$
(i.e. without simple closed leaves homotopic to boundary circles of $\R$)
as a subspace of $\ML(\S)$, the measure determined by the pair $(\R,\Gamma)$ is just
the sum of the Thurston measure on $\ML_0(\R)+\Gamma$ and of all its $\Map(\S)$-translates
%
(see  Section \ref{sec:ML} of this paper and \cite[Section 3]{LM08} for more details). 
In particular, the case $\R=\emptyset$ 
corresponds to an atomic measure on $\ML(\S)$ supported on the translates of $\Gamma$.

We recall that each measured lamination admits a unique decomposition $\lambda+\Gamma$, which we call ``standard'',
as a sum of two measured laminations with isotopically disjoint supports,
where $\Gamma$ is a simple multi-curve and $\lambda$ is scc-free.
A way to detect the nature of an $\Map(\S)$-invariant ergodic measure $m$ on $\ML(\S)$ is to consider the standard decomposition 
$\lambda+\Gamma$ of a general element in $\supp(m)$ and let $\R$ be the hull of $\lambda$.
Such standard decomposition of a measured lamination is also the key to understand the closure of its 
$\Map(\S)$-orbit (see Theorem \ref{thm:orbit-lamination} in this paper and \cite[Theorem 8.9]{LM08}).

\subsection{Main results}\label{sec:main-results}
Viewing $\ML(\S)$ as a subspace of $\Curr(\S)$ it is natural to ask what the possible $\Map(\S )$-invariant, locally finite, ergodic measures on $\Curr(\S )$ are. We will show that a classification
of such measures very much analogous to the above one holds.
%

As an example, consider a current $c$ on $\S$ such that all connected component of $\supp(c)$ which are measured laminations
are in fact weighted simple closed curves.
Consider the counting measure supported on the $\Map(\S )$-orbit of $c$, which is clearly $\Map(\S )$-invariant and ergodic. We will see in Lemma \ref{lemma:Gamma+a-bound} that the orbit of such a current $c$ 
cannot accumulate anywhere, and hence the above measure is locally finite. 

We will prove that any $\Map(\S )$-invariant, locally finite, ergodic measure supported on $\Curr(\S )$ is essentially a combination of the Thurston measure and of a counting measure supported on a current $c$ as in the above example. 

We first need a canonical way to decompose a current into more elementary pieces.

\begin{defi}[Standard decomposition]\label{def:standard}
A {\it{standard decomposition}} of a current $c\in\Curr(\S)$
is a decomposition of $c$ as a sum $c=\lambda+\Gamma+\a$ of three currents with isotopically disjoint supports
such that $\Gamma$ is a weighted simple multi-curve with support $C$, $\lambda$ is an scc-free measured lamination
and $\a$ is an a-laminational current with hull $\A$.
\end{defi}

The following result will be very useful.

\begin{mainprop}[Standard decomposition of a current]\label{mainprop:standard} 
Every current on $\S$ admits a unique standard decomposition.
\end{mainprop}

\begin{remark}
In Definition \ref{def:standard} we choose the name ``standard''
to suggest that such decomposition is well-behaved, meaning that
it is canonical and it is compatible with the action of the mapping class group.
\end{remark}


The previous statement allows us to formulate our first main result.

\begin{maintheorem}[Orbit closure of a geodesic current]\label{thm:orbit-closure}
Let $c\in\Curr(\S)$ be a non-zero geodesic current 
with standard decomposition $c=\lambda+\Gamma+\a$.
Then
\[
\ol{\Map(\S)\cdot c}=
\Map(\S)
\cdot \left(\ML_{\R}(\S)+\Gamma+\a\right)
\]
where $\R$ is the union of the components of $\S\setminus(C\cup\A)$ that intersect the support of $\lambda$. Moreover, 
the subgroup $\stab(\R,C\cup\A)\subset\Map(\S)$ of mapping classes that pointwise fix $C\cup\A$ and send $\R$ to itself
is contained inside the stabilizer
$\stab(\ML_{\R}(\S)+\Gamma+\a)$ of the locus $\ML_{\R}(\S)+\Gamma+\a$
as a finite-index subgroup.
%
\end{maintheorem}

To state the classification of ergodic invariant measures on $\Curr(\S)$,
following Lindenstrauss-Mirzakhani, we extend the notion of complete pair and of the measure it defines to our setting.

\begin{defi}[Pairs and complete pairs]\label{def:pair}
Let $\R\subset \S$ be a subsurface and let $c\in\Curr(\S)$
be a current that standardly decomposes as a sum $c=\Gamma+\a$ of a simple multi-curve $\Gamma$
and an a-laminational $\a$.
The couple $(\R,c)$ is a {\it{pair}} if $\supp(c)$ and $\R$ are isotopically disjoint.
Moreover, $(\R,c)$ is a {\it{complete pair}}
if it is a pair and each boundary curve of $\R$ is 
homotopic either to a boundary curve of $\S$, or
to a curve in the support of $\Gamma$, or
to a boundary curve of $\hull(\a)$.
%
\end{defi}

Note that Definition \ref{def:pair} reduces to Lindenstrauss-Mirzakhani's definition of a complete pair for $\a=0$, and that the case $c=0$ is not excluded. We emphasize that the couple $(R,\Gamma+\alpha)$
that appears in Theorem \ref{thm:orbit-closure} is indeed a complete pair.


Given a pair $(R,c)$, we define the corresponding measure on $\Curr(\S)$ as follows.
If $\R=\emptyset$, denote by $m^{(\emptyset,c)}=\delta_c$ the  Dirac measure
supported on the current $c$.
If $\R\neq\emptyset$,
denote by $m^{(\R,c)}$ the push-forward of the Thurston measure through
the map $\ML_0(\R)\rightarrow \Curr(\S)$ that sends $\lambda\mapsto \lambda+c$, where $\ML_0(\R)$ denotes the set of laminations supported on the interior of $\R$ (see Section \ref{subsec:ML}).

\begin{defi}[Subsurface measures]
Given a pair $(\R,c)$, the {\it{subsurface measure of type $[\R,c]$}} on $\Curr(\S)$ is
\[
m^{[\R,c]}:=\sum_\varphi m^{(\varphi(\R),\varphi(c))}
\]
where $\varphi$ ranges over 
$\Map(\S)/\stab(m^{(\R,c)})$.
\end{defi}

Again, these are the measures on $\ML(\S)$ considered by Lindenstrauss-Mirzakhani
in the case $\a=0$.


The second main result of the paper is the following.

\begin{maintheorem}[Classification of ergodic invariant measures on $\Curr$]\label{thm:main} 
The measure $m^{[\R,c]}$ on $\Curr(\S)$ is ergodic, $\Map(\S)$-invariant and locally finite
for every complete pair $(\R,c)$. Moreover, if $m$ is a locally finite, $\Map(\S )$-invariant, ergodic measure on $\Curr(\S )$,
then $m$ is a positive multiple of $m^{[\R,c]}$ for some complete pair $(\R,c)$.
\end{maintheorem}

\begin{remark}
The space $\Curr(\S)$ is $\sigma$-locally compact and metrizable, and so
completely metrizable and separable (see Theorem \ref{thm:topological} and \cite{bonahon:currents}).
We will deal with spaces obtained from 
Borel subsets of spaces of geodesic currents
by taking images via continuous maps with finite fibers, products and disjoint unions.
On such spaces every locally finite non-negative measure is a Radon measure.
We will only consider locally finite non-negative measures without further mention.
\end{remark}

We comment briefly on the ingredients in the proofs of the main results. 
The proof of Theorem \ref{thm:orbit-closure} basically relies on the following facts:
\begin{itemize}
\item
the standard decomposition of a current exists
and is unique (Proposition \ref{mainprop:standard});
\item
$\Map(\S)$ acts properly discontinuously on the locus $\Cufill(\S)$ of binding currents (Proposition \ref{properly});
\item
the $\Map(\S)$-orbit of a measured lamination with full hull is dense in $\ML(\S)$ (Theorem \ref{thm:orbit-lamination} and
\cite[Theorem 8.9]{LM08}).
\end{itemize}

The proof of Theorem \ref{thm:main} relies on
\begin{itemize}
\item
a $\Map(\S)$-invariant partition of $\Curr(\S)$ provided by Corollary \ref{cor:invariant-partition} (explained in Section \ref{sec:intro-inv});
\item
the discontinuity of the action of $\Map(\S)$ on $\Cufill(\S)$ (Proposition \ref{properly});
\item
the classification of locally finite,
ergodic,  $\Map(\S)$-invariant measures on $\ML(\S)$ obtained in \cite{LM08} and \cite{hamenstadt:measures}.
\end{itemize}

As $\RR_+$ acts on $\Curr(\S)$ by multiplication, it makes sense to speak of
$d$-homogeneous measures, namely of measures $m$ such that
$m(t\cdot U)=t^d\cdot m(U)$ for all Borel subsets $U\subset\Curr(\S)$.
Notice that the Thurston measure $\Thu$ is $N(\S)$-homogeneous,
with $N(\S):=-3\chi(\S)-n$.

In \cite[Proposition 8.2]{LM08} it is shown that a locally finite $d$-homogeneous
$\Map(\S)$-invariant measure supported on $\ML(\S)$ must satisfy
$d\neq N(\S)$.

Because of the relevance of $\Map(\S)$-invariant
homogeneous measures to curve counting problems,
we also provide a sharpening of \cite[Proposition 8.2]{LM08}
and an almost complete classification of such measures.

For every $d\in\RR$ 
consider the measure
\[
m^{(\R,c)}_d:=
\begin{cases}
m^{(\S,0)} & \text{if $(\R,c)=(\S,0)$ and $d=N(\S)$}\\
\int_0^{+\infty}t^{d-N(\R)-1}m^{(\R,tc)}dt & \text{if $c\neq 0$}
\end{cases}
\]
on $\Curr(\S)$, where $(\R,c)$ is understood to be a complete pair.
Moreover, set
\[
m^{[\R,c]}_d:=\sum_{\varphi}m^{(\varphi(\R),\varphi(c))}
\]
as $\varphi$ ranges over $\Map(\S)/\stab(\R,c)$.

\begin{maintheorem}[Classification of ergodic invariant homogeneous measures on $\Curr$]\label{thm:locally-finite-homogeneous}
Every locally finite $\Map(\S)$-invariant $d$-homogeneous ergodic measure
on $\Curr(\S)$ is a positive multiple of one of the following:
\begin{itemize}
\item[(i)]
the Thurston measure $m^{[\S,0]}_{N(\S)}=\Thu$
\item[(ii)]
the measure $m^{[\R,c]}_d$ with $c\neq 0$ and $d>N(\S)$ large enough.
\end{itemize}
In part (ii) every $d>N(\S)+N(\R)$ works.
\end{maintheorem}

In particular, the Thurston measure $\Thu$ is the $d$-homogeneous
measure with smallest $d$, and actually the only one (up to multiples)
with $d=N(\S)$. This explains its frequent occurrence in 
problems analogous to Mirzakhani's simple closed curve counting
theorem \cite[Theorem 1.1]{Mir08}. In fact, using Theorem \ref{thm:locally-finite-homogeneous} we recover one of the main
results of \cite{ES} (Proposition 4.1).

However, we emphasize that there are quite natural
counting problems that give rise to homogeneous measures
of degree higher than $N(\S)$, and which thus cannot be governed
by Thurston measure: see, for instance, Example \ref{ex:higher}.

The proof of Theorem \ref{thm:locally-finite-homogeneous} relies on
Theorem \ref{thm:main}, on an estimate of the volume of unit balls
in $\ML_\R(\S)$ for a wandering subsurface $\R$ of $\S$
(Lemma \ref{lemma:estimate-ball})
and on the following result that will be proven in Appendix \ref{app}.

\begin{mainlemma}[Asymptotic growth of $b_h^c$]\label{lemma:bound-orbit-c}
Let $(\S,h)$ be a hyperbolic surface
and let $c$ be a current of type $c=\Gamma+\alpha$.
Denote by $b_h^c([L_1,L_2])$ the number
of points in the $\Map(\S)$-orbit of $c$ with $h$-length in $[L_1,L_2]$.
Then there exists $q>1$ such that
\[
{\textstyle\frac{1}{v}} \cdot L^{N(\S)}<b_h^c([0,L]),\,b_h^c([L, qL])<v \cdot L^{N(\S)}\qquad\text{for all $L$,}
\]
for a suitable constant $v>1$ (that depends only on $\S$, $h$ and $c$).
\end{mainlemma}

The above statement 
is much weaker than Mirzakhani's Theorem 1.1 in \cite{mirzakhani-orbit}, which gives the exact asymptotics of $b_h^c$.
We mention that the upper bound contained in 
Lemma \ref{lemma:bound-orbit-c} was also proven in
\cite[Lemma 2.4]{sapir} (for closed surfaces) and
\cite[Lemma 5.6]{mirzakhani-orbit} (for binding currents).

\subsection{An invariant partition of $\Curr(\S)$}\label{sec:intro-inv}

The existence and uniqueness of the standard decomposition
(Proposition \ref{mainprop:standard}) and
the classification of locally finite, ergodic,
$\Map(\S)$-invariant measures (Theorem \ref{thm:main}) rely on a partition of $\Curr(\S)$ 
into $\Map(\S)$-invariant Borel subsets.
The key step in the construction of such partition
is the analysis of the locus $\Cufh(\S)$ of currents
of {\it{full hull}} (namely, of hull equal to $\S$) which contains
two special subsets:
the locus $\MLfh(\S)$ of laminations of full hull
and the locus $\Cufill(\S)$ of {\it{binding currents}}, i.e. of currents
$c$ that intersect every geodesic which is not asymptotic to $\pa\S$
(see Definition \ref{def:binding}).

The proof of the following result is also contained in Burger-Iozzi-Parreau-Pozzetti \cite{burger:currents}-\cite{burger:currents2}.

\begin{maintheorem}[Partition of $\Cufh$]\label{MLfill-intro}
A current of full hull on the connected surface $\S$ is either a measured lamination
or a binding current. In other words,
\[
\Cufh(\S)=\MLfh(\S)\dot{\bigcup}\Cufill(\S)
\]
in the set-theoretical sense. Moreover, both $\MLfh(\S)$ and $\Cufill(\S)$ are
$\Map(\S)$-invariant Borel subsets.
\end{maintheorem}

\begin{notation*}
Given topological subspaces $\{X_k\}$ of $X$, we will denote by
the dotted symbol $\dot{\bigcup}_k \{X_k\}$ the subspace of $X$ obtained as the union of all $X_k$'s
if such $X_k$'s are pairwise disjoint.
\end{notation*}

As a consequence of Theorem \ref{MLfill-intro}, we get the following partition of the full space of geodesic currents (see Corollary \ref{partition2}). For a subsurface $\R$ of $\S$, let $\MLfh_{\R}(\S)$ denote the subset of measured laminations 
supported in the interior of $\R$ and with hull $\R$
and similarly define $\Cufill_{\R}(\S)$ to be the subset of currents 
supported in the interior of $\R$ and that bind $\R$. Then 
\[
\Curr(\S)=\dot{\bigcup_{(\R,C,\A)}}\Cutriple_{(\R,C,\A)}(\S)
\]
with
\[
\quad \Cutriple_{(\R,C,\A)}(\S):=\MLfh_{\R}(\S) \oplus\Cufh_{C}(\S)\oplus \Cufill_{\A}(\S)
\]
where $\R,\A\subseteq\S$ are disjoint subsurfaces and
$\Cufh_{C}$ is the subspace of simple multi-curves whose support is the unweighted simple multi-curve $C\subset \S$ disjoint from $\R\cup\A$.  

The existence of the above decomposition quickly leads to the proof of
Proposition \ref{mainprop:standard}. In fact, a current $c$ must belong to a unique
$\Cutriple_{(\R,C,\A)}(\S)$, and so $c=\lambda+\Gamma+\a$ 
with $\lambda$ being a lamination of full hull in $\R$, $\Gamma$ a simple multi-curve with support $C$ and
$\a$ an a-laminational current with hull $\A$.

Finally, denoting  by $[\R,C,\A]$ a {\it{type}},
that is an equivalence class of triples $(R,C,A)$ under the action of $\Map(\S)$,
and by
\[
\Cutriple_{[\R,C,\A]}(\S):=\bigcup_{\varphi\in\Map(\S)} \Cutriple_{(\varphi(\R),\varphi(C),\varphi(\A))}(\S),
\]
we also obtain the following invariant partition of the space of currents.

\begin{maincor}[$\Map(\S)$-invariant partition of $\Curr$]\label{cor:invariant-partition}
The space $\Curr(\S)$ can be decomposed into a union over all types $[\R,C,\A]$ in $\S$ 
\[
\Curr(\S)=\dot{\bigcup_{[\R,C,\A]}}\Cutriple_{[\R,C,\A]}(\S)
\]
of the $\Map(\S)$-invariant, pairwise-disjoint, Borel subsets
$\Cutriple_{[\R,C,\A]}(\S)$. 
\end{maincor}

\subsection{Outline of the paper}

In Section \ref{sec:currents} we give the necessary background on geodesic currents. In section \ref{sec:partitionmain} we construct the partition described above
and prove Theorem \ref{MLfill-intro} as well as Corollaries \ref{partition2} and \ref{cor:invariant-partition}.
In Section \ref{sec:action} we
study the action of the mapping class group on subsets of $\Curr(\S)$ and prove Proposition \ref{properly} and Theorem \ref{thm:orbit-closure}. 
In Section \ref{sec:family} we recall the classification of invariant measures on $\ML(\S)$ 
by Lindenstrauss-Mirzakhani and Hamenst{\"a}dt, 
we construct the ergodic $\Map(\S)$-invariant subsurface measures $m^{[\R,c]}$ on $\Curr(\S)$
and we show that $m^{[\R,c]}$ is locally finite if and only if the pair $(\R,c)$ is complete.
Finally, in Section \ref{sec:proof} we prove 
Theorem \ref{thm:main} and Theorem \ref{thm:locally-finite-homogeneous} is proven in Section \ref{sec:homogeneous}.
Appendix \ref{app} contains some estimates that are used
in Theorem \ref{thm:locally-finite-homogeneous}
and the proof of Lemma \ref{lemma:bound-orbit-c}.


\subsection{Acknowledgements} 
While writing up the present article, we learnt that
a proof of the dichotomy for currents of full hull (Theorem \ref{MLfill-intro})
was independently obtained by Burger-Iozzi-Parreau-Pozzetti in 
\cite{burger:currents}.
%
%

The authors would like to thank Fran\c{c}ois Labourie and Chris Leininger for interesting conversation on the topics of this paper.
We are particularly indebted to Juan Souto for very valuable feedback, remarks and suggestions. 
We are also grateful to Marc Burger, Alessandra Iozzi, Anne Parreau and Beatrice Pozzetti for a very useful exchange of ideas
after both our preprints had appeared and for spotting some mistakes and incongruences in the first version of this work.
We also thank an anonymous referee for carefully reading our paper
and for useful comments.

V.E.~was partially supported by Academy of Finland project \#297258 and EPSRC 
grant EP/T015926/1. G.M.~was partially supported by
INdAM GNSAGA research group.

\section{The space $\Curr$ of geodesic currents}\label{sec:currents}



Let $\S$ be a smooth, compact, connected, oriented surface, possibly with boundary $\pa\S$ consisting of the closed curves $\beta_1,\dots,\beta_n$.
Assume that $\chi(\S)<0$
and let $\pi:=\pi_1(\S)$ be its fundamental group and
$\tilde{\S}\rar \S$ its universal cover.

Throughout the paper we will
call $h$ a {\it{hyperbolic metric}} on $\S$
if $h$ is a metric of curvature $-1$ on $\S$ 
and the boundary $\pa\S$ is $h$-geodesic.

\subsection{The space of bi-infinite geodesics in $\tilde{\S}$}

Fix an auxiliary hyperbolic metric $h$ on $\S$ 
and let $\tilde{h}$ be its lift to $\tilde{\S}$.
Then $\tilde{\S}$ can be identified with a subset of $\DD^2$ and we define the {\it{finite boundary}} $\pa_{f}\tilde{\S}$ to be 
the locus of points in $\tilde{\S}$ that project to $\pa\S$.
The {\it{ideal boundary}} $\pa_\infty\tilde{\S}$ is the locus of points of $\pa{\DD}^2$ in the closure of $\tilde{\S}$.
%
%
The {\it{boundary}} $\pa\tilde{\S}=\pa_{f}\tilde{\S}\cup\pa_{\infty}\tilde{\S}$ is homeomorphic to $S^1$
and it inherits an orientation from $\tilde{\S}$. Given
three distinct points $x,y,z\in S^1$ we write $x\prec z\prec y$
if a path travelling from $x$ to $y$ in the positive direction meets $z$.
If $\pa\S$ is non-empty, then $\pa_f\tilde{\S}$ is the union of countably many open intervals
and $\pa_\infty\tilde{\S}$ is a closed subset of $\pa\tilde{\S}$ with no internal part.
If $y_1,y_2\in\pa_\infty\tilde{\S}$, then we denote
by $[y_1,y_2]_\infty$ the subset of points $y\in\pa_\infty\tilde{\S}$
such that $y_1\preceq y\preceq y_2$,
and by $(y_1,y_2)_\infty$ the subset $[y_1,y_2]_\infty\setminus\{y_1,y_2\}$.

%

\begin{defi}
The {\it{space of bi-infinite geodesics on $\tilde{\S}$}} is
the space $\Geod(\tilde{\S})$ of unordered pairs of distinct points in $\pa_\infty\tilde{\S}$.
\end{defi}



\begin{remark}\label{rmk:geodesics-double}
If $\S$ has non-empty boundary, then we can view
$\S$ as a subsurface of its double $D\S$.
Thus, a bi-infinite geodesic in $\S$ gives rise to a bi-infinite geodesic in $D\S$
that does not hit $\pa\S$. Hence, we can view
$\Geod(\tilde{\S})$ as a closed subset of $\Geod(\widetilde{D\S})$.
\end{remark}

Given a compact subset $\cpt$ of $\S$ and a hyperbolic metric $h$ on $\S$,
we denote by $\Geod_{h,\cpt}(\tilde{\S})$ the subset of $\Geod(\tilde{\S})$
representing $\tilde{h}$-geodesics whose projection to $\S$ is contained in $\cpt$.
We also denote by $\GeodK(\tilde{\S})\subset\Geod(\tilde{\S})$ 
the subset of all geodesics of $\tilde{\S}$ such that
their (parametrized) projection
$\gamma:\RR\rar\S$ is {\it{reduced}}, meaning that $\gamma$ does not
contain a closed subcurve $\gamma|_{[t_1,t_2]}$
homotopic to $\beta_j\ast\beta_j$ for any $j=1,\dots,n$.

We omit the proof of the following simple observation.

\begin{lemma}\label{lemma:curve-selfint-support}
Fix a hyperbolic metric $h$ on $\S$. Then
\begin{itemize}
\item[(a)]
there exists a compact subset $\cpt$ of the interior of $\S$
such that $\GeodK(\tilde{\S})\subset\Geod_{h,\cpt}(\tilde{\S})$;
\item[(b)]
for every $s\geq 0$ there exists a compact $\cpt$ in the interior of $\S$
such that every lift of an $h$-geodesic
with at most $s$ self-intersections
is contained inside $\Geod_{h,\cpt}(\S)$.
\end{itemize}
\end{lemma}

The following observation will be useful in Section \ref{sec:binding}.

\begin{lemma}\label{lemma:reduced-curve}
Let $h$ be a hyperbolic metric on $\S$ and let $\gamma\subset\S$ be a bi-infinite geodesic.
Then there exists a reduced bi-infinite geodesic $\gamma^{red}$ whose support is
isotopic to a subset of the support of $\gamma$.
Moreover, if no end of $\gamma$ spirals about a boundary component of $\S$, the same is true of $\gamma^{red}$.
\end{lemma}
\begin{proof}
Construct a curve $\hat{\gamma}^{red}$ starting from $\gamma$
by replacing every closed subcurve homotopic to $\beta_j^{\ast l}$ with $l\geq 2$ by $\beta_j$
(resp. replacing $\beta_j^{\ast (-l)}$ with $l\geq 2$ by $\beta_j^{-1}$).
The geodesic representative $\gamma^{red}$ of $\hat{\gamma}^{red}$ is easily seen
to satisfy all the requirements.
\end{proof}

For a disconnected surface $\R=\coprod_i\R_i$ we define $\tilde{\R}:=\coprod_i\tilde{\R}_i$
and $\Geod(\tilde{\R}):=\coprod_i\Geod(\tilde{\R}_i)$. We also define
$\Geod_{h,M}(\tilde{\R})$ and $\GeodK(\tilde{\R})$ analogously.

\subsection{Geodesic currents}\label{subsec:currents}

%


Note that the group $\pi$ naturally acts on $\Geod(\tilde{\S})$ via the diagonal action
on $(\pa_\infty\tilde{\S})^2$.

%
%



\begin{defi}[Geodesic current on a closed surface]
A {\it{geodesic current}} on a surface $\S$ without boundary is 
a $\pi$-invariant locally finite measure $c$ on $\Geod(\tilde{\S})$.
\end{defi}

Given a surface $\S$ with non-empty boundary, we view
$\Geod(\tilde{\S})$ as a closed subset of $\Geod(\widetilde{D\S})$
by Remark \ref{rmk:geodesics-double}. As a consequence,
$\pi_1(D\S)$-invariant locally finite measures on $\Geod(\widetilde{D\S})$ 
can be restricted to $\pi_1(\S)$-invariant locally finite measures on $\Geod(\tilde{\S})$.

\begin{defi}[Geodesic current on a surface with non-empty boundary]
A geodesic current on a surface $\S$ with non-empty boundary
is a $\pi$-invariant, locally finite measure on $\Geod(\tilde{\S})$
obtained as the restriction of a current in $\Curr(D\S)$, which is
invariant under the natural involution of $D\S$ and whose support
does not transversely intersect $\pa\S$.
\end{defi}

We denote by $\widetilde{\supp}(c)\subset\Geod(\tilde{\S})$
the support a the geodesic current $c\in\Curr(\S)$.
Given a hyperbolic metric $h$ on $\S$,
we denote by $\wti{\supp}_h(c)\subseteq \tilde{\S}$
the union of all $\tilde{h}$-geodesics in $\wti{\supp}(c)$
and by $\supp_h(c)\subseteq\S$ the projection
of $\wti{\supp}_h(c)$ to $\S$.

\begin{remark}\label{rmk:Xi}
Let $\Xi\subset T^1\S$ be the unit tangent vectors to geodesics
that do not transversally hit $\pa\S$ (or, equivalently, whose lifts
to $\tilde{\S}$ have endpoints in $\pa_\infty\tilde{\S}$).
The datum of a geodesic current is equivalent to
a locally finite measure on $\Xi$ which is invariant under the geodesic flow.
If $\pa\S=\emptyset$, then $\Xi=T^1\S$ and
so a geodesic current can be seen as a locally finite measure on $T^1\S$
which is invariant under the geodesic flow.
%
\end{remark}
%
%
%

The {\it{space of geodesic currents $\Curr(\S)$}} on the surface $\S$
is endowed with the weak$^\star$-topology, meaning that
\[
c_k\lra c \in\Curr(\S)\qquad\iff \qquad \int_{\Geod(\tilde{\S})}f\cdot c_k\lra \int_{\Geod(\tilde{\S})}f\cdot c
\]
for all continuous functions $f:\Geod(\tilde{\S})\rar\RR$ with compact support.

%
%
%



\begin{example}[Weighted sums of closed curves]\label{ex:weightedsums}
Let $\gamma$ be a homotopically nontrivial closed curve on $\S$.
Each of its lifts $\tilde{\gamma}_i$ to $\tilde{\S}$ determines
a point in $\Geod(\tilde{\S})$. Thus, $\sum_i \delta_{\tilde{\gamma}_i}$ is
a $\pi$-invariant measure on $\Geod(\tilde{\S})$, and so a geodesic current
which is denoted by $\gamma$ with little abuse.
We can thus view the set of homotopy classes of closed curves on $\S$ as a subset
of $\Curr(\S)$. 
Clearly, given homotopically nontrivial closed curves $\gamma_1,\dots,\gamma_k$
and real numbers $w_1,\dots,w_k>0$, the linear combination
$\sum_j^k w_j\gamma_j$ is again a geodesic current,
which we call a {\it{(weighted) multi-curve}} and its support is the {\it unweighted multi-curve} 
$\cup_{j=1}^k \gamma_j$. When the curves $\gamma_j$ are simple and pairwise disjoint, we call such a current
a {\it (weighted) simple multi-curve}, and similarly, its support an {\it{unweighted simple multi-curve}}. 
\end{example}

\begin{example}[Current attached to a measured foliation]\label{example:MF}
Let $\mathcal{F}$ be a foliation on $S$ (possibly with singularities of type
$\mathrm{Re}(z^k dz^2)=0$ with $k\geq -1$) 
such that no leaf of $\mathcal{F}$ is transverse
to $\pa\S$ or spirals about some component of $\pa \S$.
If $\mathcal{F}$ is endowed with a transverse measure,
then it determines a geodesic current on $\S$ (see 
\cite{thurstonnotes}, \cite{FLP}, \cite{levitt:laminations}, \cite{Bon86},
\cite{bonahon:currents}).
\end{example}

\begin{remark}[Non-spiralling behavior of atomic leaves]
Since geodesic currents are locally finite measures, leaves with an end
that spirals about a simple closed curve cannot carry an atomic measure.
\end{remark}

A geodesic current $c\in\Curr(\S)$ is {\it{supported on the boundary of $\S$}}
if it can be written as a linear
combination $c=\sum_{j=1}^n u_j\beta_j$ of the boundary
curves with all $u_j\geq 0$. The current $c$ is {\it{internal in $\S$}} if $c(\tilde{\beta}_j)=0$
for all lifts $\tilde{\beta}_j$ of $\beta_j$ and all $j$.

\begin{remark}[Support on internal current can reach the boundary]
The support of an internal current $c$ need not be disjoint from
$\pa\S$. Consider, for example $c=\sum_j w_j\gamma_j$
a weighted
sum of all closed non-peripheral curves $\{\gamma_j\}$ in $\S$
(with rapidly decaying weights $w_j>0$ so that the sum makes sense).
\end{remark}

We denote by $\Curr_{\pa\S}(\S)$ the subset of currents supported on the boundary of $\S$
and by $\Curr_0(\S)$ the subset of currents which are internal in $\S$.
Moreover, we call $\CurrK(\S)$ the subset of internal currents
that are supported on the closure of $\GeodK(\tilde{\S})$.

Given a hyperbolic metric $h$ on $\S$
and a compact subset $\cpt$ of the interior of $\S$,
we denote by $\Curr_{h,\cpt}(\S)$ the 
subset of currents $c$
such that $\supp_h(c)\subseteq \cpt$. 
Note that $\CurrK(\S)$ is contained in $\Curr_{h,\cpt}(\S)$,
for some compact subset $\cpt$ that depends on $h$.

All of the above definitions immediately extend to disconnected surfaces.\\

We will often use the following decomposition.

\begin{lemma}[Interior+boundary decomposition of a current]\label{lemma:interior+boundary}
For every surface $\S$ with boundary $\pa\S=\bigcup_{j=1}^n\beta_j$ we have the algebraic decomposition
\[
\Curr(\S)=\Curr_{\pa\S}(\S)\oplus\Curr_0(\S)
\]
meaning that each $c\in\Curr(\S)$ can be uniquely written as a sum of a current
supported on the boundary and an internal current.
Moreover, $\Curr_{\pa\S}(\S)$ is a closed subset and
$\Curr_0(\S)$ is a dense Borel subset of $\Curr(\S)$.
\end{lemma}
\begin{proof}
Let $c\in\Curr(\S)$. We want to find $u_1,\dots,u_n\geq 0$
such that $c=c_0+\left(\sum_j u_j\beta_j\right)$ with $c_0\in\Curr_0(\S)$.
If $\tilde{\beta}_j$ is a lift of $\beta_j$ with endpoints $x_j\prec y_j$,
then it is enough to set $u_j:=c(\{x_j,y_j\})$. The uniqueness
of such choice is immediate.
To show that $\Curr_0(\S)$ is dense, it is enough to show that
each $\beta_j$ belongs to the closure of $\Curr_0(\S)$.
Now, fix a closed curve 
$\ALPH$ based at a point of $\beta_j$,
which is not homotopic to a power of $\beta_j$,
and let $\gamma_k$ be the concatenation of 
$k\cdot\beta_j$ and $\ALPH$.
Such a $\gamma_k$ is non-simple for $k\geq 2$
and so non-peripheral: it follows that $\frac{1}{k}\gamma_k$
belongs to $\Curr_0(\S)$. Clearly, $\frac{1}{k}\gamma_k\rar\beta_j$.

Since $\Curr_{\pa\S}(\S)$ is clearly closed, we are left to show that
$\Curr_0(\S)$ is Borel.
Fix a lift $\tilde{\beta}_j\in\Geod(\tilde{\S})$ of $\beta_j$
and let $(U_k)$ and $(V_k)$ be countable fundamental systems of neighbourhoods of $\tilde{\beta}_j$ such that $\ol{U}_k\subset V_k$
and $\ol{V}_k$ is compact.
Moreover let $f_k:\Geod(\tilde{\S})\rar[0,1]$ be a continuous
function with support in $V_k$ and such that $f_k|_{\ol{U}_k}\equiv 1$.
The the subset of currents $c$ whose mass fades to zero
near $\tilde{\beta}_j$ is given by
\[
\bigcap_{l\geq 1}\bigcup_{k\geq 1}\left\{c\in\Curr(\S)\ \Big|\ \int_{\Geod(\tilde{\S})} f_k\cdot c<1/l\right\}
\]
which is then a Borel subset of $\Curr(\S)$.
We conclude by observing that $\Curr_0(\S)$ is obtained by intersecting
countably many
similar subsets for all lifts of $\beta_1,\dots,\beta_n$.
\end{proof}

\subsection{The mapping class group}

Let $\Diff_+(\S)$ be the topological group of orientation-preserving
diffeomorphisms of $\S$ that send every boundary component
to itself, and let $\Diff_0(\S)$ be the subgroup
of diffeomorphisms isotopic to the identity, which is a connected component of $\Diff_+(\S)$. 
For a disconnected surface $\coprod_i \R_i$ we 
moreover require the diffeomorphisms to send every component to itself, so that
$\Diff_+(\coprod_i\R_i)\cong\prod_i \Diff_+(\R_i)$.

The mapping class group is the discrete group $\Map(\S)=\Diff_+(\S)/\Diff_0(\S)$.
If $\R$ is a subsurface of $\S$ such that all components $\R_i$ of $\R$ have negative
Euler characteristic, we denote by $\Map(\S,\R)$ the subgroup
of elements in $\Map(\S)$ that can be represented by diffeomorphisms which are the identity on $\R$,
and we define similarly $\Map(\S,C)$ if $C$ is an unweighted simple multi-curve.

We also denote by $\stab(\R)\subset\Map(\S)$ the subgroup of mapping classes
that send $\R$ to itself up to isotopy,
and by $\stab(\R,\pa\R)$ the finite-index subgroup of $\stab(\R)$
consisting of elements that send each boundary component of $\R$ to itself.
%
%
If $C$ is an unweighted simple multi-curve in $\S$, then 
$\Map(\S,C)$ is a finite-index subgroup
of $\stab(C)$.

\begin{notation*}
Suppose that $\R,\A$ are disjoint subsurfaces of $\S$ and
that $C\subset\S$ is an unweighted simple multi-curve disjoint from $\R$ and $\A$.
By slight abuse, we will denote by $\stab(\R,C\cup\A)$
the subgroup of elements of $\Map(\S)$ which send $\R$ to itself
and which restrict to the identity on $C$ and on $\A$.
We incidentally remark that Dehn twists along simple closed curves supported on $C$ belong to $\stab(\R,C\cup\A)$.
\end{notation*}

Finally, we note that
the mapping class group $\Map(\S)$
acts on $\Geod(\tilde{\S})$ and hence on $\Curr(\S)$ by self-homeomorphisms.
We denote by $\stab(c)$ the stabilizer of a current $c\in\Curr(\S)$.
Similarly, $\Map(\S)$ also acts on the space of measures on $\Curr(\S)$ by push-forward
and $\stab(m)$ denotes the stabilizer of a measure $m$ on $\Curr(\S)$.

\subsection{Push-foward of currents}\label{sec:push-forward}

Let $\R$ be a subsurface of $\S$, possibly disconnected and with boundary,
such that every connected component of $\R$ has negative Euler characteristic.

Fix an auxiliary hyperbolic metric on $\S$. 
A {\it{geodesic realization}} of $\R$ inside $\S$ is
a map $I:\R\hookrightarrow \S$ that sends the interior of $\R$
homemorphically onto its image and each boundary curve of $\pa\R$
homeomorphically onto a closed geodesic of $\S$.

Note that two boundary curves of $\R$ can be mapped to the same geodesic of $\S$.

\begin{lemma}[Geodesics in a subsurface]\label{lemma:injective}
The map $I$ induces a closed continuous map
$\tilde{I}:\Geod(\tilde{\R})\rightarrow\Geod(\tilde{\S})$.
If $\R$ is connected, then
\begin{itemize}
\item
$\tilde{I}$ is injective;
\item
given lifts $\tilde{\gamma}_1,\tilde{\gamma}_2$
of two distinct geodesics $\gamma_1,\gamma_2\subset\R$,
the image $\tilde{I}(\tilde{\gamma}_1)$ is $\pi_1(\S)$-conjugate
to $\tilde{I}(\tilde{\gamma}_2)$ if and only if 
$\gamma_1,\gamma_2\subset\pa\R$ and $I(\gamma_1)=I(\gamma_2)$.  
\end{itemize}
\end{lemma}
\begin{proof}
Clearly, it is enough to prove the statement for $\R$ connected.
As before, let $h$ be a hyperbolic metric on $\S$ and
let $I$ map every boundary component of $\R$ to a geodesic on $\S$. Endow $\R$
with the pull-back metric.

%
%
The induced map $\tilde{\R}\rightarrow \tilde{\S}$ is a local isometry onto its image. Thus we obtain a proper continuous map
$\tilde{\R}\cup \pa\tilde{\R} \rar \tilde{\S}\cup\pa\tilde{\S}$, which restricts then to
a closed map $\pa_\infty \tilde{\R}\rar \pa_\infty\tilde{\S}$.

The injectivity of $\tilde{I}$ follows from the injectivity
of $I_*:\pi_1(\R)\rar\pi_1(\S)$ and the last claim from
the identification of conjugacy classes in $\pi_1(\S)$ with
free homotopy classes of loops in $\S$.
\end{proof}

The above lemma allows us to define
a push-forward map
\[
I: \Curr(\R)\longrightarrow \Curr(\S)
\]
which we denote still by $I$ with little abuse of notation.
If $\R$ is connected, we set $\displaystyle I(c):=\sum_{[g]} g\cdot \tilde{I}(c)$, where $[g]$ ranges over $\pi_1(\S)/I_*\pi_1(\R)$
and $\tilde{I}(c)$ is the push-forward of the measure $c$ via the map $\tilde{I}$. 
If $\R=\coprod_i \R_i$ and $c_i\in\Curr(\R_i)$, then we simply let
$I(\sum_i c_i):=\sum_i I(c_i)$.
The second claim of Lemma \ref{lemma:injective}
guarantees that the restriction of
the above push-forward map to $\Curr_0(\R)$ is injective.

\begin{defi}
The subset $\Curr_{\R}(\S)$ of currents on $\S$ {\it{internal in $\R$}}
is the image of $\Curr_0(\R)$ via the push-forward map $I$.
\end{defi}


\begin{cor}\label{cor:R-subset}
The locus $\Curr_{\R}(\S)$ is a Borel subset of $\Curr(\S)$.
\end{cor}
\begin{proof}
Let $h$ be an auxiliary hyperbolic metric on $\S$
and let $\ol{\R}$ be the $h$-realization of $\R$ inside $\S$
(which is not homeomorphic to $\R$ if two boundary circles of $\R$
are isotopic two each other inside $\S$).
Clearly, $\Curr_{\R}(\S)$ is contained in the closed locus
of currents $c\in\Curr(\S)$ that have support inside $\ol{\R}$ 
and that do not transversally intersect any boundary circle of $\R$.
Adapting the proof of Lemma \ref{lemma:interior+boundary},
one can easily show that $\Curr_{\R}(\S)$ is a Borel subset inside such closed locus.
\end{proof}


%


\subsection{Intersection pairing}\label{sec:intersection}

Two geodesics $\ALPH,\ALPH'\in\Geod(\tilde{\S})$ with endpoints
$x_1,x_2$ and $x'_1,x'_2$ in $\pa_\infty\tilde{\S}$ intersect transversely if $x_1\prec
x'_1\prec x_2\prec x'_2$ or $x_1\prec x'_2\prec x_2\prec x'_1$. We denote by $\IG(\tilde{\S})$ the open subset of $\Geod(\tilde{\S})\times\Geod(\tilde{\S})$
consisting of pairs of transversely intersecting geodesics. The diagonal action of $\pi$ on $\Geod(\tilde{\S})\times\Geod(\tilde{\S})$
preserves $\IG(\tilde{\S})$ and we denote by $\widetilde{\IG}(\S)\subset \IG(\tilde{\S})$ a fundamental domain. 


\begin{defi}[Geometric intersection of currents]
Given two geodesic currents $c_1,c_2\in\Curr(\S)$, their geometric intersection number is
\[
\iota(c_1,c_2):=\int_{\widetilde{\IG}(\S)}c_1\times c_2.
\]
\end{defi}

Given two distinct closed curves $\gamma_1,\gamma_2$, the intersection number
$\iota(\gamma_1,\gamma_2)$ counts the minimal number of 
intersection points
between homotopic representatives of $\gamma_1$ and $\gamma_2$ in general position.
If $\gamma_1$ and $\gamma_2$ are non-isotopic to each other,
such minimal number is actually attained by choosing geodesic representatives
with respect to an auxiliary hyperbolic metric on $\S$.

Note that if $\gamma$ is an open geodesic arc in the hypebolic surface $(\S,h)$, then it makes sense
to speak of the intersection of $\gamma$ with a current $c$,
namely $\iota(\gamma,c):=c(\widetilde{\IG}_\gamma)$
where $\widetilde{\IG}_\gamma$ is the subset of geodesics in $\Geod(\tilde{\S})$
that transversely intersect a fixed lift of $\gamma$.


We recall the following result by Bonahon \cite{bonahon:currents}.

\begin{thm}[Continuity of geometric intersection]\label{thm:intersection}
The intersection pairing
\[
\iota:\Curr(\S)\times\Curr(\S)\longrightarrow\RR_{\geq 0}
\]
is continuous.
In particular, the function $\ell_c=\iota(c,\cdot):\Curr(\S)\rightarrow\RR_{\geq 0}$ associated to
any $c\in\Curr(\S)$ is continuous.
\end{thm}

Though the restriction of $\iota$ to $\Curr_0(\S)$ is non-degenerate,
$\iota$ itself is degenerate if $\Curr(\S)$ has boundary. In fact, for every
boundary curve $\beta_j$ of $\S$ we have $\iota(\beta_j,c)=0$ for all $c\in\Curr(\S)$.
Such $\iota$ can be modified in order to make it non-degenerate by
considering arcs that meet the finite boundary of $\tilde{\S}$. 
We will not need
such a construction here and so we
refer to \cite{duchin-leininger-rafi:flat} for further details. 

%
Let $c$ be a non-simple closed curve on $S$ that intersects every closed curve in $S$ 
(or more generally, let $c$ be a binding current as defined in Section \ref{sec:binding}).
We will see in Section \ref{sec:binding} that the function 
$\ell_c=\iota(c,\cdot)$ is strictly positive on $\Curr_{h,\cpt}(\S)\setminus\{0\}$
for every hyperbolic metric $h$ and every compact subset
$\cpt$ contained in the interior of $\S$. In particular, $\ell_c$ will be strictly
positive on $\CurrK(\S)\setminus\{0\}$.


\subsection{Liouville current attached to a metric} 

Let $\S$ be a closed surface and let $h$ be a hyperbolic metric on $\S$. Let $\Omega_g$ be the natural volume form on the unit tangent bundle $T^1\S$ of $\S$
(that pushes down to $2\pi$ times the area form $dA_g$ on $\S$).
Since $\Omega_g$ is invariant under the geodesic flow,
it defines a geodesic current $\Liou_h\in\Curr(\S)$
by Remark \ref{rmk:Xi}.
%
%

\begin{defi}[Liouville current]
The current $\Liou_h$ on $\S$ is called the {\it{Liouville current}} associated to the hyperbolic metric $h$.
\end{defi}


Here we recall an important property of Liouville currents.

\begin{prop}[Liouville current and length of closed geodesics]\label{prop:liouville}
Let $h$ be a hyperbolic metric on 
the closed surface $\S$. 
Then
\[
\iota(\Liou_h,\gamma)=\ell_h(\gamma)
\]
for every closed geodesic $\gamma$ in $\S$.
\end{prop}

The above construction is due to many authors, building on the work
of Bonahon \cite{bonahon:currents} in the hyperbolic case.
For example, Otal \cite{otal} treated the case of a smooth metric of negative curvature,
Duchin-Leininger-Rafi \cite{duchin-leininger-rafi:flat} and Bankovic-Leininger \cite{bankovic-leininger}
dealt with flat surfaces with conical points and Constantine \cite{constantine} with non-positively curved metric
with conical points.\\

Consider now a surface $\S$ with non-empty boundary.
It is possible to define the length function $\ell_h$
attached to a hyperbolic metric $h$ as follows.

\begin{remark}[Length function attached to a hyperbolic metric with geodesic boundary]\label{rmk:hyperbolic-length-function}
Let $\S$ be a surface with non-empty boundary
and consider $\S$ as embedded in its double $D\S$.
Given a hyperbolic metric $h$ on $\S$ such that $\pa\S$ is geodesic,
we can endow $D\S$
with the metric $Dh$ induced by $h$ which is invariant under the natural involution.
Identify $\Curr(\S)$ to the closed subset of $\Curr(D\S)$ supported inside $\S\subset D\S$
and let $\ell_h:\Curr(\S)\rar\RR$ be the restriction of 
continuous function $\ell_{Dh}:\Curr(D\S)\rar\RR$ to $\Curr(\S)$,
so that $\ell_h(\gamma)$ is in fact the $h$-length of $\gamma$
for every closed geodesic $\gamma$ in $\S$.
It will follow from Proposition \ref{prop:compactness} that
$\ell_{Dh}$ is proper. As a consequence, $\ell_h$ is proper too.
\end{remark}

By contrast with Proposition \ref{prop:liouville}, note that
the length function $\ell_h$ associated to
a hyperbolic metric $h$ on a surface $\S$ 
with boundary $\pa\S=\bigcup_j \beta_j$ as in the above remark
is not induced from a Liouville-type current $\Liou_h$ on $\S$
that fits our definitions, since $\Liou_h$ must satisfy
$\iota(\Liou_h,\beta_j)=0$ whereas $\ell_h(\beta_j)\neq 0$
for all $j$.



\begin{example}[Hyperbolic metrics with cusps]
Let $\S'$ be a punctured surface and let $h'$ be
a hyperbolic metric with cuspidal ends on $\S'$.
A Liouville current $\Liou_{h'}$ can be defined quite in the same way
as above, but it is not locally finite since $\S'$ has cusps.
Fix a homeomorphism $f:\S\rar\S'$ of the surface
with boundary $\S$ onto
its image, which is also a homotopy equivalence.
Such maps lifts to $\tilde{f}:\tilde{\S}\rar\tilde{\S}'$
which extends to the respective boundaries.
In particular, if $x'\in\pa\tilde{\S}'$ corresponds to a cusp,
then $\tilde{f}^{-1}(x')$ consists of two points in $\pa_\infty\tilde{\S}$
that bound an interval of $\pa_f\tilde{\S}$; otherwise
$\tilde{f}^{-1}(x')$ consists of one point.
Since the subset of geodesics in $\Geod(\tilde{\S}')$ with any fixed common endpoint
in $\pa\tilde{\S}'$ has $\Liou_{h'}$-measure $0$, 
one can easily define a pull-back measure $\tilde{f}^*\Liou_{h'}$.
Such a measure is not locally finite though.
\end{example}


\subsection{Measured laminations}\label{subsec:ML}

An important subspace of $\Curr(S)$ is given by measured laminations. Here we recall a few facts about this space. 

\begin{defi}[Measured geodesic laminations]
A {\it{geodesic lamination}} on the hyperbolic surface $(\S,h)$
is a closed subset $\Lambda\subset \S$
that is foliated by complete geodesics.
A {\it{measured geodesic lamination}} is a 
geodesic lamination $\Lambda$ endowed with a measure $\lambda$
on the space $\Arcs(\Lambda)$ of arcs
that are transverse to $\Lambda$ and with endpoints in $\S\setminus\Lambda$ such that
\begin{itemize}
\item[(i)]
$\lambda$ is non-negative and $\lambda(\ALPH)>0$ if and only if $\ALPH\cap\Lambda\neq\emptyset$;
\item[(ii)]
if $\ALPH,\ALPH'\in\Arcs(\Lambda)$ and the endpoint of $\ALPH$ agrees with the starting point of $\ALPH'$,
then $\lambda(\ALPH\ast\ALPH')=\lambda(\ALPH)+\lambda(\ALPH')$;
\item[(iii)]
if $(\ALPH_t)_{t\in[0,1]}$ is a continuous family of arcs in $\Arcs(\Lambda)$, then
$\lambda(\ALPH_0)=\lambda(\ALPH_1)$.
\end{itemize}
\end{defi}

\begin{example}[Simple multi-curves]
Let $\gamma_1,\dots,\gamma_l\subset\S$ be pairwise disjoint simple closed geodesics which are homotopically nontrivial and let $w_1,\dots,w_l>0$. Let $\Lambda=\bigcup_{i=1}^l \gamma_i$ and $\lambda$ be
the transverse measure defined by $\lambda(\ALPH)=\sum_{i=1}^l w_i\cdot |\ALPH\cap \gamma_i|$
for every $\ALPH\in\Arcs(\Lambda)$. Then $(\Lambda,\lambda)$ is a measured lamination
of special type, namely a (weighted) simple multi-curve. As in Example \ref{ex:weightedsums}, we often denote it by
$w_1\gamma_1+\dots+w_l\gamma_l$ and its support is the unweighted simple multi-curve $\bigcup_{i=1}^l\gamma_l$.
\end{example}

A geodesic lamination
determines a $\pi$-invariant closed subset $\tilde{\Lambda}$ of $\Geod(\ti{\S})$ and
a measured geodesic lamination $(\Lambda,\lambda)$ determines a geodesic current
in $\Curr(\S)$ supported on $\tilde{\Lambda}$. By abuse of notation, we will denote
such a current just by $\lambda$.

Given two hyperbolic metrics $h,h'$, 
there is a canonical correspondence between
$h$-geodesic laminations and $h'$-geodesic laminations,
and hence it makes sense just to speak of ``laminations'' on $\S$.
Similarly, the concept of ``measured laminations'' is independent of
the chosen hyperbolic metric.

The {\it{space of measured laminations}} $\ML(\S)$
is the locus of currents in $\Curr(\S)$ induced
by a measured geodesic lamination on $\S$.
We denote by $\ML_0(\S)$ the subset of measured laminations
internal in $\S$, so that $\ML(\S)=\ML_0(\S)\oplus \Curr_{\pa\S}(\S)$.



\begin{remark}
For every auxiliary hyperbolic metric $h$ on $\S$, 
there is a compact subset of the internal part of $\S$ (that depends on $h$ only)
which contains the support of every measured geodesic lamination in
$\ML_0(\S)$. Hence, $\ML_0(\S)$ is a closed subset of $\Curr(\S)$.
\end{remark}

Since measured geodesic
laminations are currents
supported on a set of pairwise non-intersecting simple geodesics, 
the following characterization holds (Bonahon \cite{bonahon:currents}).

\begin{lemma}[$\ML$ as a quadratic cone in $\Curr$]
The locus of geodesic measured laminations $\ML(\S)$ 
can be characterized inside $\Curr(\S)$ as
the closed quadratic $\RR_+$-cone $\ML(\S)=\{c\in\Curr(\S) \ |\ \iota(c,c)=0\}$.
\end{lemma}

Lastly, recall that the space $\ML_0(\S)$ can be described using charts given by train tracks. This allowed Thurston to prove that $\ML_0(\S)$ can be given the structure of a manifold, piecewise-linearly homeomorphic to a Euclidean space of dimension $-3\chi(\S)-n$ (where $n$ is the number of boundary components).

%


\subsection{Currents and spikes}

Here we recall a basic well-known property of geodesic currents, namely the fact that
no mass can be supported on a subset of geodesics which enter a spike.

\begin{defi}
Let $\tilde{\ALPH},\tilde{\ALPH}'\subset\DD^2$ be two distinct geodesics
which are asymptotic to the same point $x\in\pa\DD^2$, and let $\tilde{\sigma}\subset\DD^2$
be the region bounded by $\tilde{\ALPH},\tilde{\ALPH}'$.
A {\it{spike}} is a hyperbolic surface isometric to the end of $\tilde{\sigma}$ that is asymptotic to $x$.
\end{defi}

Endow the surface $\S$ with a hyperbolic metric $h$ and
consider two semi-infinite oriented geodesic rays $\ALPH,\ALPH':[0,\infty)\rar\S$ that are asymptotic to each other.
Two lifts $\tilde{\ALPH},\tilde{\ALPH}'\subset\tilde{\S}$ which are asymptotic to the same
point in $\pa\tilde{\S}$ determine a spike $\tilde{\sigma}$.


\begin{lemma}[Geodesics constrained inside a spike have measure zero]\label{lemma:spike}
Let $x$ be a point in $\pa_\infty\tilde{\S}$
and let $y_1,y_2\in\pa_\infty\tilde{\S}$ such that
$x\notin [y_1,y_2]_\infty$.
Suppose that
the geodesic $\ALPH_y$ in $\S$ determined by  
$\{x,y\}\in\Geod(\tilde{\S})$ is not closed for any $y\in [y_1,y_2]_{\infty}$.
Then every current $c\in\Curr(\S)$ satisfies $c(\{x\}\times (y_1,y_2)_\infty)=0$.
\end{lemma}
\begin{proof}
Let $\gamma$ be an open geodesic arc of bounded length in $\S$.
If each $\ALPH_y$ transversely intersects
$\gamma$ at least $M$ times, then
$\iota(c,\gamma)\geq M\cdot c(\{x\}\times (y_1,y_2)_\infty)$.
In order to prove the statement, it is thus enough to show that there exists
an arc $\gamma$ which is transversely intersected infinitely many
times by each $\ALPH_y$.

Now fix $y_0\in (y_1,y_2)_\infty$ and let $\ALPH_{y_0}:\RR\rar T^1\S$ be the 
projection of the geodesic that runs from $y_0$ to $x$.
Consider an accumulation point $v\in T^1\S$ for
$\ALPH_0(t)$ as $t\rightarrow+\infty$ and let
$\gamma$ be a small geodesic arc transverse to $v$.
For every $y\in (y_1,y_2)_\infty$ the geodesic $\ALPH_y$ is non-closed and bi-infinite 
and it accumulates at $v$, and so $\ALPH_y$ transversely crosses
$\gamma$ infinitely many times.
\end{proof}

We then have a criterion to determine whether the support of a current $c$
is disjoint from a spike.

\begin{cor}[Currents with support disjoint from a spike]\label{cor:spike}
In the hypotheses of Lemma \ref{lemma:spike},
suppose moreover that the support of $c$ does not transversally cross
$\ALPH_{y_1}$ and $\ALPH_{y_2}$.
Then $\{x\}\times (y_1,y_2)_\infty$ is disjoint from the support of $c$.
\end{cor}

\begin{proof}
Let $X\subset \pa_\infty\tilde{\S}$ be an open neighborhood of $x$ that does not
intersect $[y_1,y_2]_\infty$.
It is easy to see that, for every $x'\in X$ different from $x$ and for every 
$y\in (y_1,y_2)_\infty$, there exists an $i\in\{1,2\}$ such that
the geodesics $\{x',y\}$ and $\{x,y_i\}$ transversely intersect.
By our hypotheses, the geodesic $\{x',y\}$ is not in the support of $c$.
It follows that the open subset $X\times (y_1,y_2)_\infty$ of $\Geod(\tilde{\S})$
does not meet the support of $c$ and the conclusion follows.
\end{proof}

An easy consequence of the above lemma is the existence of geodesics 
not asymptotic to the boundary in the support of any non-zero current.

\begin{cor}\label{cor:non-spiral-geodesic}
Let $0\neq c_0\in\Curr_0(\S)$.
Then there exists a geodesic in $\supp_h(c_0)$, different
from a boundary curve, which is not asymptotic to a boundary curve.
\end{cor}
\begin{proof}
By contradiction, suppose that all geodesics in $\supp_h(c_0)$ are either boundary
curves or asymptotic to them, and so in particular
they are either boundary curves or they are bi-infinite geodesics.
Then $\wti{\supp}(c_0)$ is contained in a countable
union of sets of type $\{x\}\times [y_1,y_2]_\infty$, with $x\in\pa_\infty\tilde{\S}$
and $x\notin [y_1,y_2]_\infty\subset\pa_\infty\tilde{\S}$.
Note that $\ALPH_y$ is not a closed geodesic for all $y\in(y_1,y_2)_\infty$. Hence,
we can apply Lemma \ref{lemma:spike}
and conclude that the whole $\{x\}\times (y_1,y_2)_\infty$ has $c_0$-measure $0$. 
Since $c_0$ is internal in $\S$, it follows that
$\{x\}\times[y_1,y_2]_\infty$ has $c_0$-measure $0$ too.
As a consequence, $c_0=0$ and we have achieved a contradiction.
\end{proof}

\section{An invariant partition of $\Curr$}\label{sec:partitionmain}


In this section we will discuss the key ingredients in the proof of 
Theorem \ref{thm:main}: a partition of the space of currents 
into Borel invariant subsets
and the action of the mapping class group on subsets of $\Curr(\S)$. 

\subsection{Binding currents}\label{sec:binding}

We start by discussing a class of very general currents.

\begin{defi}[Binding currents]\label{def:binding}
A current $c\in\Curr(\S)$ \emph{binds} if every geodesic in $\Geod(\tilde{\S})$
with no endpoint in the closure of $\pa_f\tilde{\S}$
is transversely intersected by a geodesic in the support of $c$. 
Denote by $\Cufill(\S)$ the subspace of binding currents in $\Curr(\S)$.
\end{defi}

We observe that a binding current may well belong to $\Curr_0(\S)$
and that elements of $\Curr_{\pa\S}(\S)$ are never binding.
On the other hand, if $b$ is binding and $c$ is any current,
then $b+c$ is clearly binding.

\begin{example}[Liouville currents on closed surfaces]
If $\S$ is closed and $h$ is a hyperbolic metric on $\S$,
then the associated Liouville current $\Liou_h$ is binding
since it has full support in $\Geod(\tilde{\S})$.
\end{example}

\begin{example}[Binding currents supported on closed geodesics]\label{example:binding}
A binding geodesic current can be obtained
by considering the current associated to the multi-curve
$b=w_1\gamma_1+\dots+w_l\gamma_l$
where all $w_i>0$ and
each $\gamma_i$ is a closed curve in $\S$, such that their geodesic representatives (with respect
to an auxiliary hyperbolic metric) cut $\S$ into a disjoint union of disks
and cylinders homotopic to a boundary curve of $\S$.
Actually, it is possible to have $l=1$.
\end{example}

The following example was proposed by Marc Burger.

\begin{example}[Binding sums of countably many weighted closed curves]
Fix an auxiliary hyperbolic metric on $\S$
and let $(\gamma_k)$ be the set of non-peripheral closed curves in
$\S$, ordered so that $\ell_h(\gamma_k)\leq\ell_h(\gamma_{k+1})$.
Let $b=\sum_k w_k\gamma_k$ where $(w_k)$ is a 
quickly decreasing sequence
of positive numbers, for example $w_k=2^{-k}$.
Then $b$ is certainly binding. In fact, there exists $k$
such that $\S\setminus\gamma_k$ is a disjoint union of disks
and cylinders homotopic to a boundary curve of $\S$.
An analogous binding current can be manufactured by
only adding up simple closed curves.
\end{example}


We begin the analysis of the binding locus by the following simple observation.

\begin{lemma}[Density of binding currents inside $\Curr$]
The subset $\Cufill(\S)$ is non-empty and dense inside $\Curr(\S)$.
\end{lemma}
\begin{proof}
Example \ref{example:binding} shows that
$\Cufill(\S)$ is not empty. 
Concerning the density, just note that,
if $b$ is binding and $c\in\Curr(\S)$, then
$c_k=c+\frac{1}{k}b$ is a sequence of binding currents
that converge to $c$ as $k\rar\infty$.
\end{proof}


Binding currents can be characterized using the intersection pairing:
the following statement was essentially proven in \cite{glorieux:exponents},
up to minor variations. 
%
%
%

\begin{prop}[Positivity of binding currents]\label{prop:positivity}
Let $c$ be a current on $\S$ and fix a hyperbolic metric $h$ on $\S$. The following are equivalent.
\begin{itemize}
\item[(a)]
$c$ binds;
\item[(b)]
$\iota(c,c')>0$ for every $0\neq c'\in\Curr_0(\S)$;
\item[(c)]
for any given compact subset $\cpt\subset\inte{\S}$, the current $c$
satisfies
$\iota(c,c')>0$ for all $0\neq c'\in\Curr_{h,\cpt}(\S)$;
\item[(d)]
$\iota(c,c')>0$ for every $0\neq c'\in\CurrK(\S)$.
\end{itemize}
\end{prop}
\begin{proof}
Suppose first that (a) holds and let $0\neq c'\in\Curr_0(\S)$.
By Corollary \ref{cor:non-spiral-geodesic}, there exists
a geodesic $\gamma'$ in $\supp_h(c')$ which is neither a boundary curve nor asymptotic to a boundary curve. Then $c$ intersects $\gamma'$ transversally and so $\iota(c,c')>0$. Hence, (a) implies (b).

Clearly, (b) implies (c) because $\Curr_{h,\cpt}(\S)\subset\Curr_0(\S)$, and (c) implies (d) because
$\CurrK(\S)\subset\Curr_{h,\cpt}(\S)$ for a suitable $\cpt$.
Thus, we only need to show that (d) implies (a).

Suppose that $c$ is not binding and so there exists a complete geodesic $\gamma'$
which is neither a boundary curve nor asymptotic to a boundary curve, and which is not
transversally intersected by $\supp_h(c)$.
Then Lemma \ref{lemma:non-binding} below guarantees the existence of
a current $0\neq c'\in\CurrK(\S)$ which satisfies $\iota(c,c')=0$.
\end{proof}

The following lemma was essentially proven by Glorieux in \cite{glorieux:exponents}.
A proof taylored to our need is included for completeness.

\begin{lemma}[\cite{glorieux:exponents}]\label{lemma:non-binding}
Let $h$ be a hyperbolic metric on $\S$ and let $\gamma'\subset\S$ be a non-peripheral geodesic
with no end hitting or spiralling about a boundary component of $\S$.
If the support of $c\in\Curr(\S)$ does not transversely intersect $\gamma'$,
then there exists $0\neq c'\in\CurrK(\S)$ supported on the closure of $(\gamma')^{red}$
such that $\iota(c,c')=0$.
\end{lemma}

\begin{proof}
Let $(\gamma')^{red}$ be a reduced geodesic
obtained from $\gamma'$ as in Lemma \ref{lemma:reduced-curve}.
Since the support of $(\gamma')^{red}$ is isotopic to a subset $(\hat{\gamma}')^{red}$ of the support of $\gamma'$,
it follows that $(\gamma')^{red}$ does not transversely intersect $\supp_h(c)$.
In fact, if a curve $\gamma''$ in $\supp_h(c)$ is disjoint from $\gamma'$, then it is also disjoint from $(\hat{\gamma}')^{red}$.
On the other hand, if $\gamma''=\gamma'$, then $\gamma'$ must be simple and so $(\gamma')^{red}=\gamma'$.

%


Since transversality is an open condition, 
every geodesic contained in the closure of $(\gamma')^{red}$ 
is either disjoint from $\supp_h(c)$ or completely contained inside $\supp_h(c)$.
Thus, a geodesic current $c'$ with $\supp_h(c')$ contained in the closure of $(\gamma')^{red}$
satisfies $c'\in\CurrK(\S)$ and $\iota(c,c')=0$. 

In order to construct such non-zero $c'$, 
we produce a measure on $\Xi\subset T^1\S$ supported on the closure of $(\gamma')^{red}$ which is invariant
under the geodesic flow.

More explicitly, consider an arc-length parametrization of $(\gamma')^{red}$, which we denote by
little abuse still by $(\gamma')^{red}:\RR_t\rightarrow \Xi\subset T^1\S$. For every $r>0$, denote by $c'_r$ the probability measure $(\gamma')^{red}_*(\frac{1}{2r}\chi_{[-r,r]}|dt|)$ on $\Xi$, which is
supported on $(\gamma')^{red}([-r,r])$. Then a weak$^\star$-limit $c$ of the measures $c'_r$
as $r\rar\infty$ satisfies the requirements.
%
%
%
%
%
%
%
\end{proof}

The following will be an immediate consequence of
Proposition \ref{prop:positivity} and Proposition \ref{prop:compactness} and it will be proven in the next section.

\begin{cor}[Openness of the binding locus]\label{cor:binding-open}
The locus $\Cufill(\S)$ is open inside $\Curr(\S)$.
\end{cor}

%
%
%

\subsection{Topological properties of $\Curr(\S)$}



The following compactness result is well-known.
In the present form we will directly derive it from
Bonahon's work \cite{bonahon:currents} on closed surfaces.


\begin{prop}[Compactness of sublevels of a binding current]\label{prop:compactness}
The projectivization $\PP\Curr(\S)$ is compact and so are the closed subspaces
$\PP\CurrK(\S)$ and $\PP\Curr_{h,\cpt}(\S)$.
Moreover,
if $b$ is a binding current, then
the restriction of $\ell_b:\Curr(\S)\rar \RR_{\geq 0}$ 
to $\CurrK(\S)$ and to $\Curr_{h,\cpt}(\S)$ is proper.
\end{prop}
\begin{proof}
For $\S$ closed, the first claim was proven by Bonahon in \cite{bonahon:currents}.
Suppose now that $\pa\S\neq\emptyset$.
By embedding $\S$ inside its double $D\S$, which comes endowed with a natural involution
$\sigma$, we can identify $\Curr(\S)$ to the closed subset of $\Curr(D\S)$ consisting
of currents on $D\S$ which are $\sigma$-invariant and which do not intersect $\pa\S$.
Since $\PP\Curr(D\S)$ is compact, it follows that $\PP\Curr(\S)$ is too.

As for the second claim, consider a diverging sequence $(c_k)$ inside $\CurrK(\S)$.
Since $\PP\CurrK(\S)$ is compact, there exists $0\neq c\in\CurrK(\S)$ such that
$[c_k]\rar[c]$, namely there exist $w_k\in\RR_+$ such that
$w_k c_k\rar c$. Since $(c_k)$ is divergent, $w_k\rar 0$.
Moreover, $w_k\ell_b(c_k)=\ell_b(w_k c_k)\rar \ell_b(c)>0$, which implies that
$\ell_b(c_k)\rar\infty$. This shows that the restriction of $\ell_b$ to $\CurrK(\S)$ is proper.

Note that the only properties of $\CurrK(\S)$ we used to prove the second claim are that
$\CurrK(\S)$ is closed inside $\Curr(\S)$ and that a binding current
positively intersects every element of $\CurrK(\S)$.
Thus, an analogous proof works for $\Curr_{h,\cpt}(\S)$.
\end{proof}

The below result
was also proven by Bonahon in \cite{bonahon:currents}
for closed surfaces.

\begin{thm}[Topological properties of $\Curr$]\label{thm:topological} The space $\Curr(\S)$ is
locally compact, $\sigma$-compact and metrizable.
As a consequence, it is also completely metrizable and second countable.
\end{thm}
\begin{proof}
Local compactness and $\sigma$-compactness follow
from Proposition \ref{prop:compactness}.
Moreover, Bonahon \cite{bonahon:currents} showed that
$\Curr(\S)$ is metrizable if $\S$ is a closed surface.

Suppose now that $\S$ has non-empty boundary,
consisting of components $\beta_1,\dots,\beta_n$.
Let $D\S$ be the double of $\S$ so that we can view
$\S$ as naturally embedded inside $D\S$,
and let $\sigma$ be the natural involution of $D\S$ that fixes $\pa\S$.
The space of currents $\Curr(\S)$ can be identified
to the locus of all $c\in\Curr(D\S)$ which are invariant
under $\sigma$ and 
such that $\iota(c,\beta_1+\dots+\beta_n)=0$.
Since $\sigma$ acts as a self-homeomorphism of $\Curr(D\S)$ 
and
the intersection pairing is continuous, $\Curr(\S)$ is a closed subset of $\Curr(D\S)$ and so
the conclusion follows from Bonahon's work.
\end{proof}

We will deal with disjoint unions, products 
and countable-to-one images of Borel subsets of some spaces of currents inside some $\Curr(\S)$.
As mentioned in the introduction, any locally finite measure on such spaces is a Radon measure.\\


%

To conclude this section, we show how the openness
of the binding locus follows from the above results.

\begin{proof}[Proof of Corollary \ref{cor:binding-open}]
Consider a sequence $(c_k)$ in the complement of $\Cufill(\S)$
inside $\Curr(\S)$ that converges to $c\in\Curr(\S)$.
We want to show that $c$ is not binding.

Fix an auxiliary binding current $b$ on $\S$.
By Proposition \ref{prop:positivity},
for every $c_k$ there exists a current $0\neq c'_k\in\CurrK(\S)$ such that
$\iota(c_k,c'_k)=0$. Also,
since $b$ is binding, we have $\iota(b,c'_k)>0$ for all $k$.
Up to rescaling $c'_k$, we can then assume that $\iota(b,c'_k)=1$
for all $k$. Now, $\ell_b^{-1}(1)\cap \CurrK(\S)$
is compact by Proposition \ref{prop:compactness}
and so, up to subsequences, 
$(c'_k)$ converges to some $c'\in\CurrK(\S)$
such that $\iota(b,c')=1$.
In particular, $c'\neq 0$.
By continuity of the intersection pairing, $\iota(c,c')=0$. This shows that $c$ is not binding.
\end{proof}

\subsection{Hull of a current}\label{sec:hull}

%


%

Before defining the hull,
let us first recall the following notion.

\begin{defi}[Simple closed curve components of a current]
A simple closed curve $\gamma\subset\S$ is a connected
component of $c\in\Curr(\S)$ if 
$\iota(\gamma,c)=0$ and there exists $\e>0$
such that $c-\e\gamma$ is a (non-negative) current.
\end{defi}

If no simple closed curve is a connected component of $c$, then
clearly $c$ is {\it{scc-free}} in the sense of Definition
\ref{def:scc-free}. Thus a measured lamination
$\lambda$ is scc-free if and only if it has no closed leaf. On the other hand, an internal binding current is always scc-free.

\begin{remark}\label{rmk:scc-component}
Let $\gamma$ be a simple closed curve which is a connected component
of $c$ and let $w=c(\tilde{\gamma})$ for some lift $\tilde{\gamma}$ of $\gamma$ to $\tilde{\S}$.
Then $c-t\gamma$ is a current (i.e. it is non-negative) if and only if $t\leq w$. In this case, the supports of $c-t\gamma$ and of $\gamma$ are isotopically disjoint. Moreover,
$\gamma$ is not a connected component of $c-t\gamma$ if and only if $t=w$. 
\end{remark}

Fix now an scc-free current $\check{c}\neq 0$ on $\S$.

Let $\R_1,\R_2$ be two closed subsurfaces inside $\S$
and denote by $I_1:\R_1\rightarrow \S$ and $I_2:\R_2\rightarrow \S$
their geodesic realizations with respect to some auxiliary hyperbolic metric $h$ on $\S$.
Suppose that $\supp_h(\check{c})$ is contained inside both
$I_1(\mathring{\R}_1)$ and $I_2(\mathring{\R}_2)$.
Then $\supp_h(\check{c})$ is contained inside their intersection, which is
an open subsurface with piecewise smooth boundary.
We denote by $\R_1\cap\R_2$ the isotopy class of subsurfaces
with smooth boundary homotopic to $I_1(\mathring{\R}_1)\cap I_2(\mathring{\R}_2)$
inside $\S$.
We say that (the isotopy class of) $\R_1$ is smaller than (the isotopy class of) $\R_2$
if $I_1(\mathring{\R}_1)\subseteq I_2(\mathring{\R}_2)$,
and so $\R_1\cap\R_2$ is isotopic to $\R_1$.

Now recall the following from Section \ref{sec:conventions}.

\begin{named}{Definition \ref{def:hull}}[Hull of a current]
The \emph{surface hull} of an scc-free current $\check{c}\in\Curr(\S)$ 
is the isotopy class $\hull(\check{c})$ of the smallest closed subsurface of $\S$ that contains the support of $\check{c}$.
\end{named}

Note that $\hull(\check{c})$ is not necessarily connected.
We denote by $\hull_h(\check{c})$ a surface homeomorphic to $\hull(\check{c})$
endowed with an $h$-geodesic realization $\hull_h(\check{c})\rar \S$
and we remind the reader that the interior of $\hull_h(\check{c})$ is embedded inside $\S$, whereas the realization map can identify couples of boundary
components of $\hull_h(\check{c})$.\\
%
%
%

A general current $c$ can have connected components which are weighted simple closed curves.
A first step toward a standard decomposition of $c$ is the following.

\begin{lemma}[$\Gamma$-summand of a current]\label{lemma:Gamma-summand}
Every $c\in\Curr(\S)$ can be uniquely written as $c=\check{c}+\Gamma$, where
\begin{itemize}
\item[(a)]
$\Gamma$ is a weighted simple multi-curve, for which $\supp_h(\Gamma)$ can be
isotoped to be disjoint from $\supp_h(\check{c})$ (for some hyperbolic metric $h$)
\item[(b)]
$\check{c}$ is scc-free
\item[(c)]
$\check{c}$ is internal in $\hull(\check{c})$.
\end{itemize}
\end{lemma}
\begin{proof}
Consider the set $\{\gamma_i\}$
of all simple closed curves $\gamma_i$ in $\S$
that are connected components of $c$.
Since the $\gamma_i$ must be disjoint,
there exist finitely many of such.
For each $i$, let $\tilde{\gamma}_i\subset\tilde{S}$ be a lift of $\gamma_i$
and let $w_i:=c(\gamma_i)>0$.
Define the weighted multi-curve $\Gamma$ as 
$\Gamma:=\sum_i w_i\gamma_i$ and the non-negative current $\check{c}$ as
$\check{c}:=c-\Gamma$.
Clearly, $\iota(\Gamma,\check{c})=0$ and so (a) holds.
Property (b) is a consequence of Remark \ref{rmk:scc-component}
and (c) follows from (b).

The uniqueness of $\Gamma$ and $\check{c}$ follows from the above construction.
%
%
%
\end{proof}

\begin{defi}\label{defi:full-hull}
A current $c=\check{c}+\Gamma\in\Curr(\S)$ has \emph{full hull} if $\hull(\check{c})=\S$
and $\Gamma$ is supported on $\pa\S$.
The subset of currents on $\S$ with full hull is denoted by $\Cufh(\S)$
and the subset of measured laminations on $\S$ with full hull is denoted by
$\MLfh(\S)$.
\end{defi}


\begin{remark}
A measured lamination $\lambda$ on $\S$
has full hull if and only if it transversely intersects every non-peripheral
simple closed curve. Such laminations are sometimes called ``filling''.
In the literature the term ``filling current'' is sometimes used to denote
what we call a binding current. These two notion of filling are really
different: for example, a measured lamination $\lambda$ cannot be binding
since $\iota(\lambda,\lambda)=0$. For this reason, we choose not use the
word ``filling'' at all.
\end{remark}


Consider now the case of a simple multi-curve
and let $C=\bigcup_{j=1}^l\gamma_j$ be a union of $l$ pairwise disjoint simple closed curves in $\S$.
By analogy with Defintion \ref{defi:full-hull}, we say that a multi-curve $\Gamma$ is supported
on $C$ if $\Gamma=\sum_{j=1}^l w_j\gamma_j$ with all $w_j\geq 0$, and that
it has support equal to $C$ if all $w_j>0$. We will denote by $\Cufh_C(\S)$ the subsets of all
multi-curves with support equal to $C$.\\

\subsection{Complement of the support of a current in its hull}


Fix a hyperbolic metric $h$ on $\S$.
If $\lambda\in\ML(\S)$ is a measured lamination without isolated closed geodesics
in its support, then 
the complement $\mathring{\hull}_h(\lambda)\setminus\supp_h(\lambda)$ consists of a finite union of
\begin{itemize}
\item
geodesic {\it{polygons}} with ideal vertices (and so ends isometric to {{spikes}})
\item
{\it{crowns}}, i.e.
open annuli such that one end is a boundary component
of $\hull_h(\lambda)$ and the other end has finitely many infinite geodesics;
such infinite geodesics come with a cyclic ordering and any two adjacent
ones bound a spike.
\end{itemize}
Clearly, every boundary circle of $\hull_h(\lambda)$ necessarily bounds a crown contained in $\hull_h(\lambda)$.


For a current which is not necessarily a lamination,
polygons must be replaced by topological disks with locally convex (not necessarily smooth) boundary
and possibly spikes, and crowns must be allowed to have locally convex non-peripheral end (and possibly spikes).

\begin{lemma}[Complement of the support of current in its hull]\label{lemma:complement-current}
Let $0\neq \check{c}$ be an scc-free current on $\S$. Then
$\mathring{\hull}_h(\check{c})\setminus \supp_h(\check{c})$
consists of convex disks and locally convex crowns, possibly with spikes.
\end{lemma}
\begin{proof}
Up to looking at the preimage of $\supp_h(\check{c})$ through the geometric realization
map $\hull_h(\check{c})\rar\S$, we can reduce to the case of a $\check{c}$ of full hull.

Assuming then that $\check{c}$ has full hull, let $\S_{\check{c}}$ be the open subsurface $\S\setminus\supp_h(\check{c})$ and let $\bar{\S}_{\check{c}}$ be its metric completion. 

We claim that $\bar{\S}_{\check{c}}$ has locally convex boundary.
In fact, consider a point $x\in\pa\bar{\S}_{\check{c}}$ and let 
$\bar{D}_{\check{c}}(x)$ be a small closed disk of radius $r$ in $\bar{\S}_{\check{c}}$ centered at $x$;
such $\bar{D}_{\check{c}}(x)$ is the metric completion of
a connected component $D_{\check{c}}(x)$ of $D(x)\cap\S_{\check{c}}$,
where $D(x)$ is the closed disk of radius $r$ centered at $x$ in $\S$.
Realize $D(x)$ as a disk inside $\DD^2$
and fix a point $y\in D_{\check{c}}(x)$. Every portion of a geodesic in the support of $\check{c}$ that meets
$D(x)$ can be realized as a portion of a geodesic $\gamma$ in $\DD^2$, and we denote by
$H_\gamma$ the closed half-plane in $\DD^2$ bounded by such $\gamma$ and that contains $y$.
It follows that $\bar{D}_{\check{c}}(x)$ is isometric to the intersection of $D(x)$ with all such $H_\gamma$,
and so it is convex.

Now, let $\bar{\S}'_{\check{c}}$ be a component of $\bar{\S}_{\check{c}}$.
A possible homotopically nontrivial simple closed curve $\gamma$ inside $\bar{\S}'_{\check{c}}$
must be homotopic to some boundary circle of $\S$, because $\check{c}$ has full hull. This
shows that $\bar{\S}'_{\check{c}}$ must be either a topological disk or a topological cylinder homotopic to a boundary component of $\S$.
It is immediate to see that the only possible ends of $\bar{\S}'_{\check{c}}$ are spikes.
\end{proof}

Since the $h$-area of $\S$ is fixed, Gauss-Bonnet theorem ensures
that
\begin{itemize}
\item[(a)]
every convex disk or locally convex crown 
in $\mathring{\hull}_h(\check{c})\setminus \supp_h(\check{c})$
can only have finitely many spikes;
\item[(b)]
there are finitely many components of
$\mathring{\hull}_h(\check{c})\setminus \supp_h(\check{c})$
that are disks with at least $3$ spikes
or crowns with at least $1$ spike.
\end{itemize}

Since both ends
of a bi-infinite geodesic entirely
contained $\mathring{\hull}_h(\check{c})\setminus \supp_h(\check{c})$ must enter a spike,
each such bi-infinite geodesic must be completely
contained in a component above mentioned in (b) and it must be isolated.
We thus have the following consequence.

\begin{cor}[Isolation of geodesics in the complement of a current in its hull]\label{cor:complement-current}
Let $0\neq \check{c}$ be an scc-free current on $\S$. 
Then geodesics completely contained in 
$\mathring{\hull}_h(\check{c})$ whose image in $\S$ does not meet $\supp_h(\check{c})$ are bi-infinite and isolated.
\end{cor}

\subsection{A partition of the space of currents of full hull}\label{sec:partition}

We can now prove that currents of full hull are either laminations of full hull
or binding currents. In particular, Theorem \ref{MLfill-intro} stated in the introduction
is a consequence of Proposition \ref{MLfill} and Lemma \ref{lemma:fh-Borel}.
Such result is a key building block for the construction of a $\Map(\S)$-invariant partition of $\Curr(\S)$.

\begin{prop}[Partition of $\Cufh$]\label{MLfill}
A current of full hull on the connected surface $\S$ is either a measured lamination
or a binding current. In other words,
\[
\Cufh(\S)=\MLfh(S)\dot{\bigcup}\Cufill(\S)
\]
in the set-theoretical sense. Moreover, both $\MLfh(\S)$ and $\Cufill(\S)$ are
$\Map(S)$-invariant.
\end{prop}

In order to prove Proposition \ref{MLfill} we will need the following technical result.

\begin{lemma}[Laminations not intersecting currents of full hull]\label{lemma:sublamination}
Assume $\S$ is connected and let $0\neq \lambda'\in\ML(\S)$ and $c\in\Cufh(\S)$ such that
$\iota(\lambda',c)=0$.
Then $\lambda'$ has full hull.
\end{lemma}
\begin{proof}
Note that no component of $\lambda'$ is a non-peripheral simple closed curve,
because $\iota(\lambda',c)=0$ and $c$ has full hull. Thus,
it is enough to prove that no geodesic in $\supp(c)$ transversely crosses $\pa\hull(\lambda')$,
from which it follows that $\pa\hull(\lambda')=\pa\S$.

By contradiction, suppose that a geodesic $\ALPH\in\supp(c)$
crosses $\pa\hull(\lambda')$ and enters a crown in $\hull(\lambda')\setminus\supp(\lambda')$,
thus ending in a spike bounded by the geodesics $\ALPH_{1},\ALPH_2$.
Let $\tilde{\ALPH},\tilde{\ALPH}_{1},\tilde{\ALPH}_{2}$ be lifts
of $\ALPH,\ALPH_1,\ALPH_2$ on $\tilde{\S}$ 
with endpoints $\{x,y\}$, $\{x,y_1\}$ and $\{x,y_2\}$.
Up to reversing the roles of $\ALPH_1,\ALPH_2$, we can assume
that $y\in (y_1,y_2)_\infty$.
Note that $\tilde{\ALPH}_1,\tilde{\ALPH}_2$ are not transversally intersected
by $\supp(c)$, because $\ALPH_1,\ALPH_2$ belong to $\supp(\lambda')$
and $\iota(\lambda',c)=0$.


If no geodesic in $\{x\}\times (y_1,y_2)_\infty$ projects to a closed curve in $\S$, then
$\{x\}\times (y_1,y_2)_\infty$ is disjoint from $\supp(c)$ by Corollary \ref{cor:spike}, and we achieve a contradiction.
Suppose then that there exists $y_0\in (y_1,y_2)_\infty$ such that the geodesic $\tilde{\ALPH}_{0}$
with endpoints $x,y_0$ projects to a closed curve $\ALPH_{0}$ (which must necessarily be simple).
Since $\ALPH_1,\ALPH_2$ are both asymptotic to $\ALPH_0$,
the support of $c$
cannot transversely intersect $\ALPH_0$. But this contradicts the fact that $c$ has full hull.
%
%
\end{proof}

The above result can be amplified as follows.

\begin{cor}[Currents not intersecting currents of full hull]\label{cor:subcurrent}
Assume $\S$ is connected and let $c'\in\Curr(\S)$
and $c\in\Cufh(\S)$ such that $\iota(c',c)=0$.
Then $c'=\lambda'$ is a measured lamination.
Moreover, if $\supp(\lambda')$ is not contained in $\pa\S$,
then $\lambda'$ has full hull.
\end{cor}
\begin{proof}
The second claim is exactly Lemma \ref{lemma:sublamination}.
Thus, it is enough to show that $c'$ is a measured lamination.

If $\beta_j$ is the $j$-th boundary circle of $\S$, we can write $c'=c'_0+\sum_j w_j \beta_j$, with $c'_0\in\Curr_0(\S)$.
If $c'_0=0$, then $c'$ is a simple multi-curve. Thus, 
we now consider the case $c'_0\neq 0$.

By Corollary \ref{cor:complement-current}, non-peripheral
geodesics in $\S$ disjoint from $\supp_h(c)$
are bi-infinite and isolated.
Hence, they cannot belong to the support of $c'_0$.
It follows that $\supp(c'_0)\subseteq\supp(c)$ and so $\iota(c'_0,c'_0)=0$.
Hence $c'_0$ is a measured lamination and so $c'$ is.
\end{proof}

Now we can complete the proof of the main statement in this subsection.

\begin{proof}[Proof of Proposition \ref{MLfill}]
The last assertion is immediate, so we concentrate on the partition of $\Cufh(\S)$.

Let $c\in\Cufh(\S)$ and let $h$ be a hyperbolic metric on $\S$.
Suppose that $c$ is not a binding current
so that,
%
by Proposition \ref{prop:positivity},
there exists $0\neq c'\in\Curr_0(\S)$ that satifies $\iota(c,c')=0$.
By Corollary \ref{cor:subcurrent}, the current $c'$ is a measured lamination of full hull.
%
%
%
Thus, by reversing the roles of $c$ and $c'$
in Corollary \ref{cor:subcurrent},
we get that $c$ is a measured lamination too.
\end{proof}


\subsection{A partition of the space of currents}\label{sec:partition2}

Let $\R\subset\S$ be a closed subsurface.
If $\R$ is connected, we denote by $\Cufh_{\R}(\S)$ the image of $\Cufh_0(\R)$
via the map $\Curr(\R)\rar\Curr(\S)$ quite analogously to Section \ref{sec:push-forward}.
If $\R=\coprod_i\R_i$ is disconnected, we let $\Cufh_{\R}(\S):=\bigoplus_i \Cufh_{\R_i}(\S)$.
We will also use the symbols
$\Cufill_{\R}(\S)$ and $\MLfh_{\R}(\S)$ with analogous meanings.

\begin{defi}
A {\it{disjoint triple in $\S$}} is an isotopy class of $(\R,C,\A)$, where $\R,\A$ are disjoint subsurfaces of $\S$
and $C$ is an unweighted simple multi-curve disjoint from $\R\cup\A$.
A {\it{type}} is an equivalence class of
triples $(\R,C,\A)$ under the action of $\Map(\S)$. The type of $(\R,C,\A)$ will be denoted by $[\R,C,\A]$.
\end{defi}

In order to construct a decomposition of the space of currents whose parts
are indexed by disjoint triples we first determine the nature of the building blocks of such decomposition.



\begin{lemma} \label{lemma:fh-Borel}
For every subsurface $\R\subset\S$,
the loci $\Cufh_{\R}(\S)$ and $\MLfh_{\R}(\S)$ are Borel subsets of $\Curr(\S)$.
\end{lemma}
\begin{proof} 
By Corollary \ref{cor:R-subset} the
set of currents supported in the interior of $\R$ 
is a Borel subset of $\Curr(\S)$.
By Lemma \ref{lemma:Gamma-summand}, a current in $\R$ can be written
as $c=\Gamma+\check{c}$ in such a way that $\Gamma$ is a multi-curve and
$\check{c}$ is scc-free.
If $c$ does not have full hull in $\R$, then there exists a
proper subsurface $\R'\subset\R$ that contains the support of $\Gamma$
and the hull of $\check{c}$. It follows that
\[
\Cufh_{\R}(\S)=
\Curr_{\R}(\S)\setminus
\bigcup_{\R'\subsetneq\R}\Curr_{\R'}(\S).
\] 
Since $\Curr_{\R'}(\S)$ is Borel,
we deduce that $\Cufh_{\R}(\S)$ is a Borel subset of $\Curr_{\R}(\S)$,
and so of $\Curr(\S)$.

Finally, $\MLfh_{\R}(\S)=\ML(\S)\cap \Cufh_{\R}(\S)$ and so $\MLfh_{\R}(\S)$ is Borel too,
because $\ML(\S)$ is closed.
\end{proof}


The above discussion gives rise to the desired decomposition
of $\Curr(\S)$ as follows.

\begin{cor}[Partition of $\Curr$]\label{partition2}
The space of geodesic currents on $\S$ can be partitioned
into Borel subsets as follows
\[
\Curr(\S)=\dot{\bigcup_{(\R,C,\A)}}\Cutriple_{(\R,C,\A)}(\S)
\quad \text{with}
\quad \Cutriple_{(\R,C,\A)}(\S):=\MLfh_{\R}(\S)\oplus\Cufh_{C}(\S)\oplus \Cufill_{\A}(\S) 
\]
where $(\R,C,\A)$ ranges over all disjoint triples in $\S$.
\end{cor}

As a first consequence of the above corollary, we obtain the existence of
the {\it{standard decomposition of a current}}
in the sense of Definition \ref{def:standard}. Recall 
from Definition \ref{def:scc-free} that we say a current 
is {\it{a-laminational}} if it is scc-free and
no connected component of its support is a lamination.
 
\begin{named}{Proposition \ref{mainprop:standard}}[Standard decomposition of a geodesic current]
Every geodesic current $c\in\Curr(\S)$
admits the following unique {\it{standard decomposition}}
as a sum
\[
c=\lambda+\Gamma+\a
\]
of three currents with isotopically disjoint supports:
an scc-free measured lamination $\lambda$,
a simple multi-curve $\Gamma$ and
an a-laminational current $\a$.
%
\end{named}


We can rearrange the subsets appearing in Corollary \ref{partition2}
in order to obtain a mapping class group invariant partition
by considering
\[
\Cutriple_{[\R,C,\A]}(\S):=\bigcup_{\varphi}
\Cutriple_{(\varphi(\R),\varphi(C),\varphi(\A))}(\S)
\]
of $\Curr(\S)$, where the unions are taken over all $\varphi$
ranging over $\Map(\S)/\stab(\R)\cap\stab(C)\cap\stab(\A)$.
Clearly, each $\Cutriple_{[\R,C,\A]}(\S)$ depends only on the type $[\R,C,\A]$
and it is $\Map(\S)$-invariant.
We have thus shown the following.

\begin{named}{Corollary \ref{cor:invariant-partition}}[$\Map(\S)$-invariant partition of $\Curr$]
The space $\Curr(\S)$ can be decomposed into a union
\[
\Curr(\S)=\dot{\bigcup_{[\R,C,\A]}}\Cutriple_{[\R,C,\A]}(\S)
\]
over all types $[\R,C,\A]$ in $\S$ of the $\Map(\S)$-invariant, disjoint Borel subsets
$\Cutriple_{[\R,C,\A]}(\S)$. 
\end{named}

\section{The action of the mapping class group on $\Curr$}\label{sec:action}

The aim of this section is to study the action of  the mapping class group
$\Map(\S)$ on the space $\Curr(\S)$ of geodesic currents on $\S$.
In particular, we will determine which currents have finite stabilizers,
an invariant locus in $\Curr(\S)$ on which the action is properly discontinuous,
and we will use a result of Lindenstrauss-Mirzakhani \cite{LM08}
to determine all orbit closures.

\subsection{Action on the locus of binding currents}

We recall that, for $\S$
closed, $\Map(\S)$ acts properly discontinuously on Teichm\"uller space,
that is, on the space of Liouville currents associated to hyperbolic metrics on $\S$
and that such Liouville currents bind $\S$.

In this section we will show the following statement.

\begin{prop}[Proper discontinuous action on $\Cufill$]\label{properly}
The mapping class group $\Map(\S )$ acts properly discontinuous on $\Cufill(\S )$
and with closed orbits (as subsets of $\Curr(\S)$).
\end{prop}

We begin by recalling the following two well-known
Lemma \ref{lemma:bound} and Lemma \ref{lemma:div-orbit-binding-multicurve}.

\begin{lemma}[Binding currents bound each other]\label{lemma:bound}
Let $h$ be a hyperbolic metric on $\S$ and $\cpt$ be a compact subset of the interior of $\S$.
Fix $\K\subset\Cufill(\S)$ compact.
There exists a constant $r>0$ (that depends on $h$, $\cpt$ and $\K$) such that 
\[
\frac{1}{r}<\frac{\iota(b_1,c)}{\iota(b_2,c)}<r
\qquad\text{and}\qquad
\frac{1}{r}<\frac{\iota(b_1,c)}{\ell_h(c)}<r
\]
for all $0\neq c\in\Curr_{h,\cpt}(\S )$ and $b_1,b_2\in \K$.
\end{lemma}
\begin{proof}
Since $b_i$ is binding, $\iota(b_i,c)>0$ for $i=1,2$ and all $c\in \Curr_{h,\cpt}(\S)\setminus\{0\}$.
Hence, the function $f:\K\times \K\times\left(\Curr_{h,\cpt}(\S)\setminus\{0\}\right)\rightarrow\RR$ defined as
$$f(b_1,b_2,c):= \frac{\iota(b_1,c)}{\iota(b_2,c)}$$
is continuous and positive, since $\iota$ is continuous.
Moreover, $f$ is homogenous in the third entry in the sense that $f(b_1,b_2,t\cdot c)=t\cdot f(b_1,b_2,c)$ for all $t>0$. Hence $f$ descends to a continuous function
$\ol{f}:\K\times \K\times\PP\Curr_{h,\cpt}(\S)\rar\RR_+$. By Proposition \ref{prop:compactness} the space
$\K\times \K\times\PP\Curr_{h,\cpt}(\S)$ is compact and hence $f$ is bounded from above and below by positive numbers. 
The same proof works for the inequalities on the right, just replacing $\iota(b_2,\cdot)$ by $\ell_h$.
\end{proof}

\begin{remark}[Binding currents are comparable on $\Teich^{\geq s}$]\label{rmk:binding-comparable}
Inequalities analogous to the left ones in Lemma \ref{lemma:bound}
hold if we replace $\iota(c,\cdot)$ by some length function $\ell_h$.
Namely, for every $s>0$ and binding currents $b_1,b_2$
in $\Curr(\S)$ there exists a constant $r>0$
that depends on $s,b_1,b_2$ such that
\[
\frac{1}{r}<\frac{\ell_h(b_1)}{\ell_h(b_2)}<r
\]
for every hyperbolic metric $h$ on $\S$ with $\sys(h)\geq s$.
\end{remark}
\begin{proof}[Proof of Remark \ref{rmk:binding-comparable}]
Consider $\S$ as embedded inside its double $D\S$, which is naturally equipped with an orientation-reversing involution $\sigma$, so that hyperbolic metrics on $\S$ double to
$\sigma$-invariant hyperbolic metrics on $D\S$.
Thus, the space $\Teich^{\geq s}(\S)$ of hyperbolic metrics
on $\S$ with systole at least $s$ can be seen as a subset
of $\Teich(D\S)$.
The closure of $\Teich^{\geq s}(\S)$ inside Thurston compactification of $\Teich(D\S)$ is obtained by adding certain 
projective classes $[\lambda]$ of $\sigma$-invariant measured laminations in $D\S$.
We claim that, for every point $[\lambda]$ in such closure
of $\Teich^{\geq s}(\S)$, we have $\iota(\lambda,b)>0$
for every current $b\in\Cufill(\S)$.
The claim implies that the function $h\mapsto\frac{\ell_h(b_1)}{\ell_h(b_2)}$ continuously and positively extends to the closure of $\Teich^{\geq s}(\S)$. Since such closure is compact, the wished conclusion follows.

In order to prove the claim, 
let $(h_k)\subset\Teich^{\geq s}(\S)$ 
be a sequence that converges to $[\lambda]$ and assume, by contradiction, that $\iota(\beta,\lambda)=0$ and so
$\lambda$ is supported on $\pa\S$.
Fix a pair of pants in $\S$ that is adjacent to a boundary component $\pa_j\S$ in the support of $\lambda$, and which is obtained by doubling a hexagon with edges $\beta,\check{\eta}_2,\eta_1,\check{\beta},\eta_2,\check{\eta}_1$
(where $\beta$ doubles to $\pa_j\S$). Call $\gamma$ the shortest arc inside such hexagon that connects $\beta$ with $\check{\beta}$, 
and which thus splits $\beta$ (resp.~$\check{\beta}$)
into the union $\beta'\cup \beta''$ (resp.~$\check{\beta}'\cup\check{\beta}''$)
of two sub-intervals. 
Up to relabeling and extracting a subsequence of $(h_k)$, we can assume that $\ell_{h_k}(\beta'')\geq \ell_{h_k}(\beta)/2$ and that $\beta'',\check{\eta}_2,\eta_1,\check{\beta}'',\gamma$ form a pentagon $P$ with five right angles.

Note that $\gamma$ doubles to an arc in $\S$ with endpoints in $\pa_j\S$,
which thus doubles to a simple closed curve in $D\S$ that meets $\pa_j\S$ twice. On the other hand, $\eta_1$ doubles to a simple closed curve in $\S$ whose geometric intersection with $\pa\S$ is zero. It follows that 
$\ell_{h_k}(\gamma)\rar\infty$ and that $\frac{\ell_{h_k}(\eta_1)}{\ell_{h_k}(\gamma)}\rar 0$.
On the other hand, $\cosh(\ell_{h_k}(\eta_1))=\sinh(\ell_{h_k}(\gamma))\sinh(\ell_{h_k}(\beta''))\geq \sinh(\ell_{h_k}(\gamma))\sinh(s/2)$ by elementary trigonometry of the pentagon $P$.
Since $\ell_{h_k}(\gamma)\rar\infty$, it follows that
$\limsup_k \frac{\ell_{h_k}(\eta_1)}{\ell_{h_k}(\gamma)}\geq 1$ and we have reached a contradiction.
\end{proof}

\begin{lemma}[Divergence of the orbit of a binding multi-curve]\label{lemma:div-orbit-binding-multicurve}
Let $b'$ be a finite binding multi-curve in $\S$.
Then $\{\varphi\in\Map(\S)\,|\,\iota(b',\varphi(b'))\leq L\}$ is a finite
set for all $L>0$.
\end{lemma}

The above lemma can be proved in a purely topological way. However, we just include
a proof that exploits the properties of the hyperbolic length function of a current in a closed surface.

\begin{proof}[Proof of Lemma \ref{lemma:div-orbit-binding-multicurve}]
Fix a hyperbolic metric $h$ on $\S$ and let $D\S$ be the closed hyperbolic surface
obtained by doubling $\S$.
Any current $c\in\Curr_0(\S)$ can be viewed as a current on $D\S$, invariant under the natural orientation-reversing involution,
and the length of $c$ can be defined as half of the length of such doubled current on $D\S$.
It follows that $\ell_h(c)>0$ for all $0\neq c\in\Curr_0(\S)$.

Clearly, it is enough to prove the statement for a finite
binding multi-curve $b'\in\Curr_0(\S)$.
Recall that $\varphi(b')$ is supported inside some compact subset $\cpt\subset\S$ for all $\varphi$
by Lemma \ref{lemma:curve-selfint-support}(b). By Lemma \ref{lemma:bound} there exists $r>0$ such that
\[
\frac{\ell_h(c)}{\iota(b',c)}<r
\]
for all $c\in\Curr_{h,\cpt}(\S)$.
Thus, taking $c=\varphi(b')$, we obtain $\iota(b',\varphi(b'))>\frac{1}{r}\ell_h(\varphi(b'))$
for all $\varphi$. The result now follows by noting that $(S,h)$ contains finitely many simple closed geodesics
of length at most $rL$.
\end{proof}

We can now prove the main proposition of this section.


\begin{proof}[Proof of Proposition \ref{properly}]
We have to show that, given $\K\subset\Cufill(\S)$ compact and
given $\{\varphi_j\}$ a sequence of distinct elements in $\Map(\S)$,
the union $\bigcup_j \varphi_j(\K)$ is closed and
$\K\cap\varphi_j(\K)=\emptyset$ for $j$ large enough.

Fix $b'$ a finite binding multi-curve in $\S$ and define
\[
m_j=\min \iota(\varphi_j(\K),b')>0.
\]
Here, $\min \iota(\varphi_j(\K),b')=\min_{b\in \K}\iota(\varphi_j(b),b')$.  
Since $\K$ is compact, 
the function $\iota(\cdot,b')$ is bounded on the union $\K\cup\varphi_1(\K)\cup\dots\cup\varphi_k(\K)$ for all $k\geq 1$.
Hence, it is enough to show that $m_j\rar\infty$ as $j\rar\infty$.

%
%

Note that, equivalently, 
$m_j=\min\iota(\K,\varphi^{-1}_j(b'))$. Fix a hyperbolic metric $h$ and note that
there exists a compact subset $\cpt$ in the interior of $\S$ that contains
the geodesic representatives of $\varphi_j^{-1}(b')$ for all $j$ by Lemma \ref{lemma:curve-selfint-support}(b).
Hence, $\varphi_j^{-1}(b')\in\Curr_{h,\cpt}(\S)$ for all $j$.

By applying Lemma \ref{lemma:bound} to the compact subset $\K\cup\{b'\}$,
there exists $r>0$ such that
\[
\frac{\iota(b', \varphi^{-1}_j(b'))}{\iota(b,\varphi^{-1}_j(b'))}<r
\]
for all $j$ and all $b\in\K$. By taking the minimum over $b\in\K$ we obtain
$m_j>\frac{1}{r}\iota(\varphi_j(b'),b')$.
We conclude that $m_j\rar\infty$ because $\iota(\varphi_j(b'),b')\rar\infty$ as $j\rar\infty$
by Lemma \ref{lemma:div-orbit-binding-multicurve}.
\end{proof}

%
%
%
%
%
%
%

In fact,
$\Cufill(\S)$ is the maximal subset of $\Curr(\S)$
on which $\Map(\S)$ acts properly discontinuously
and with closed orbits (as subsets of $\Curr(\S)$).

\begin{prop}\label{prop:non-discrete}
Let $c\in\Curr(\S)$ a current which is not binding.
Then either there are infinitely many $\varphi\in\Map(\S)$ such that $\varphi(c)=c$,
or the orbit $\Map(\S)\cdot c$ is not a closed discrete subset of $\Curr(\S)$.
\end{prop}
\begin{proof}
It follows from Lemma \ref{lemma:stabilizer}(a) and
Theorem \ref{thm:orbit-closure} proven in the next section.
\end{proof}

\subsection{Mapping class group orbits of currents}\label{sec:orbits}

In this subsection we analyze
orbits of currents under the action of $\Map(\S)$, or a
quotient of it, and in particular we determine whether they are closed and
whether they have finite stabilizers.

We begin with a simple observation.

\begin{lemma}[Stabilizer of a current]\label{lemma:stabilizer}
Let $c=\Gamma+\check{c}$ be the sum of 
a simple multi-curve $\Gamma$ and an scc-free current $\check{c}$
with isotopically disjoint supports.
The stabilizer of $c$ for the action of $\Map(\S)$ on $\Curr(\S)$
satisfies the following properties.
\begin{itemize}
\item[(a)]
$\stab(c)$ is finite if and only if $c$ has full hull in $\S$.
\item[(b)]
$\stab(c)$ contains $\Map(\S,\hull(\check{c})\cup C)$ as a finite-index subgroup, where $C$ is the support of $\Gamma$.
\end{itemize}
\end{lemma}
\begin{proof}
Note first that $\Map(\S,\hull(\check{c})\cup C)$ is always contained inside $\stab(c)$
and that $\stab(c)$ is contained inside $\stab(\hull(\check{c})\cup C)$.
Moreover, it is enough to consider $c\in\Curr_0(\S)$.

Let us first prove (a).
%
If $\hull(\check{c})\subsetneq \S$, the group $\Map(\S,\hull(\check{c}))$ is infinite
and so is $\stab(c)$.
Suppose now that $\hull(\check{c})=\S$. Then either $c=\lambda$ is a lamination
or $c=\a$ is binding by Proposition \ref{MLfill}.
If $c=\a$ is binding in $\S$, then its stabilizer is finite by 
Proposition \ref{properly}.
Suppose then that $c=\lambda$ and realize its support
by a geodesic lamination with respect to some hyperbolic metric with geodesic boundary on $\S$.
The stabilizer $\stab(\lambda)$ acts by permuting the components of $\S\setminus\lambda$
and its edges. Since $\lambda$ has full hull, the complement $\S\setminus\lambda$ consists
of finitely many ideal polygons and crowns homotopic to boundary circles of $\S$: hence, the above action
of $\stab(\lambda)$ is faithful. It follows that $\stab(\lambda)$ is finite.

In order to prove (b) we must show that $\stab(c)/\Map(\S,\hull(\check{c})\cup C)$ is finite.
The finite-index subgroup of elements in $\stab(c)/\Map(\S,\hull(\check{c})\cup C)$
that send each component of $C$ to itself
and each component and each boundary circle of $\hull(\check{c})$ to itself
identifies to $\stab_{\Map(\hull(\check{c}))}(\check{c})$.
Since $\check{c}$ has full hull in $\hull(\check{c})$, the group
$\stab_{\Map(\hull(\check{c}))}(\check{c})$ is finite by part (a)
and so $\stab(c)/\Map(\S,\hull(\check{c})\cup C)$ is finite too.
%
%
\end{proof}

As a consequence, we obtain the analogous of Lemma \ref{lemma:div-orbit-binding-multicurve}
for currents of type $\Gamma+\a$.

\begin{lemma}[Orbits of currents $\Gamma+\a$]\label{lemma:Gamma+a-bound}
Let $c=\Gamma+\a\in\Curr(\S)$ the sum of a simple multi-curve $\Gamma$ with support $C$ and an a-laminational current $\a$
with support $A$ such that $A\cap C=\emptyset$. Given a finite binding multi-curve $b'$, the set
\[
\{\varphi\in\Map(\S)/\Map(\S,A\cup C)\ |\ \iota(\varphi(c),b')\leq L\}
\]
is finite for all $L>0$.
In particular, the orbit of $c=\Gamma+\a$ is closed.
\end{lemma}
\begin{proof}
By Lemma \ref{lemma:stabilizer}(b), the quotient $\stab(c)/\Map(\S,A\cup C)$ is finite
and so it is enough to analyze $\{\varphi\in\Map(\S)/\stab(c)\ |\ \iota(\varphi(c),b')\leq L\}$.
Similarly, if $\a_0$ is a finite multi-curve with hull $A$, $\stab(\S,A\cup C)$ has finite index inside $\stab(\Gamma+\a_0)$.

Fix a hyperbolic metric $h$ and $\cpt$ a compact subset of the interior of $\S$ that contains the geodesic representative of $b'$.
As in the proof of Lemma \ref{lemma:bound},
the map $\PP\Curr_{h,\cpt}(\S)\rar\RR_+$ defined as
\[
[b']\mapsto \frac{\iota(\Gamma+\a,\varphi^{-1}(b'))}{\iota(\Gamma+\a_0,\varphi^{-1}(b'))}
\]
takes values in a closed bounded interval of $\RR_+$. 
Hence, it is enough to prove the statement
for $\a$ a finite (non-simple) multi-curve.

For $\a$ a finite multi-curve, the current $c$ can be written as $c=\sum_{i=1}^l w_l\gamma_l$
and $\bigcap_l\stab(\gamma_l)$ has finite index inside $\stab(c)$, so that we only need to prove
the finiteness of
\[
\left\{\varphi\in\Map(\S)/\bigcap_l\stab(\gamma_l)\ |\ \sum_l w_l\cdot\iota(\varphi(\gamma_l),b')\leq L\right\}.
\]
Recall that, for every $\ell>0$, the set of closed curves $\gamma$ in $\S$
with $\ell_h(\gamma)\leq \ell$ is finite. As in the proof of Lemma \ref{lemma:div-orbit-binding-multicurve},
this implies that $\{\gamma\,|\,\iota(\gamma,b')\leq \ell\}$ is finite too and so it concludes the argument.
\end{proof}

We now discuss the closure of the orbits.
The case of a measured lamination was analyzed by Lindenstrauss-Mirzakhani
\cite[Theorem 8.9]{LM08}, and we recall their result here. 

\begin{thm}[Orbit closure of a measured lamination]\label{thm:orbit-lamination}
Let $\lambda+\Gamma\in\ML(\S)$, where $\Gamma$ is a simple multi-curve
with support $C$ and $\lambda$ is a measured lamination with no closed leaves.
Then $\ol{\Map(\S)\cdot(\lambda+\Gamma)}=\Map(\S)\cdot\left(\ML_{\R}(\S)+\Gamma\right)$,
where $\R$ is the union of the components of $\S\setminus C$ that intersect the support of $\lambda$.
\end{thm}

By virtue of the above Theorem \ref{thm:orbit-lamination}
we can complete our analysis of the closure of $\Map(\S)$-orbits of currents.

\begin{named}{Theorem \ref{thm:orbit-closure}}[Orbit closure of a geodesic current]
Let $c\in\Curr(\S)$ be a non-zero geodesic current 
with standard decomposition $c=\lambda+\Gamma+\a$
into a measured lamination $\lambda$ without closed leaves, a simple multi-curve $\Gamma$ with support $C$ and an a-laminational current $\a$ with hull $A$. 
Then
\[
\ol{\Map(\S)\cdot c}=
\Map(\S)
\cdot \left(\ML_{\R}(\S)+\Gamma+\a\right)
\]
where $\R$ is the union of the components of $\S\setminus(C\cup\A)$ that intersect the support of $\lambda$. Moreover, $\stab(\ML_{\R}(\S)+\Gamma+\a)$ contains $\stab(\R,C\cup\A)$ as a finite-index subgroup.
\end{named}
\begin{proof}
Consider the second claim and
note that $\stab(\ML_{\R}(\S)+\Gamma+\a,\pa\R)$ has finite index inside 
$\stab(\ML_{\R}(\S)+\Gamma+\a)$
and $\stab(\R,C\cup\A\cup\pa\R)$ has finite index inside $\stab(\R,C\cup\A)$,
and that 
$\stab(\ML_{\R}(\S)+\Gamma+\a,\pa\R)$ contains $\stab(\R,C\cup\A\cup\pa\R)$.
The conclusion follows, since the restriction to $\S\setminus\R$ identifies
$\stab(\ML_{\R}(\S)+\Gamma+\a,\pa\R)/\stab(\R,C\cup\A\cup\pa\R)$
with $\stab_{\Map(\S\setminus\R)}(\Gamma+\a)/\Map(\S\setminus\R,C\cup\A)$,
which is finite by Lemma \ref{lemma:stabilizer}(b).

As for the first claim,
recall that $\Map(\R)\cdot\lambda$ is dense inside $\ML_0(\R)$ by Theorem \ref{thm:orbit-lamination}.
As a consequence, $\stab(\R,C\cup\A)\cdot\lambda$ is dense inside $\ML_{\R}(\S)$ and so $\stab(\R,C\cup\A)\cdot (\lambda+\Gamma+\a)$
is dense inside $\ML_{\R}(\S)+\Gamma+\a$.
Thus, it is enough to show that $\Map(\S)\cdot \left(\ML_{\R}(\S)+\Gamma+\a\right)$ is a closed
subset of $\Curr(\S)$.

Let then $c_k=\varphi_k\cdot(\lambda_k+\Gamma+\a)$ be sequence
in $\Map(\S)\cdot \left(\ML_{\R}(\S)+\Gamma+\a\right)$ that converges to $\ol{c}\in\Curr(\S)$.
We have to show that 
$\ol{c}\in \varphi(\ML_{\R}(\S)+\Gamma+\a)$
for some $\varphi\in\Map(\S)$.

If $\a=0$, the result follows from Theorem \ref{thm:orbit-lamination}; so we assume $\a\neq 0$.
By Lemma \ref{lemma:Gamma+a-bound}, the convergence of $c_k$ implies that
the subset $\{[\varphi_k]\}\subset\Map(\S)/\Map(\S,A\cup C)$ is finite and so,
up to subsequence, we can assume that it is constant.
This implies that there exists $\varphi\in\Map(\S)$ such that
$\varphi^{-1}\varphi_k\in\Map(\S,A\cup C)$ for all $k$.
%
%
%
Again up to subsequence, we can assume that
the permutation $\sigma$ of the components of
$\S\setminus(A\cup C)$ induced by
$\varphi^{-1}\varphi_k$ is independent of $k$.

Suppose that such permutation $\sigma$ is the identity.
As a consequence, $\varphi^{-1}\varphi_k(\R)=\R$ for all $k$
and so $\varphi^{-1}\varphi_k(\lambda_k)\in\ML_{\R}(\S)$ is converging to some $\ol{\lambda}\in\ML_{\R}(\S)$.
Finally, we conclude that $c_k\rar \varphi(\ol{\lambda}+\Gamma+\a)$.

If $\sigma$ is not the identity, then
$A$ is necessarily empty, $\S$ has no boundary, $\S\setminus C$ consists of two components $\S',\S''$ and each circle in $C$ belongs to $\ol{\S}'\cap\ol{\S}''$. Thus $\sigma$ must exchange $\S'$ and $\S''$, which thus have the same genus. Hence, there exists $\psi\in\Map(\S,C)$ that flips $\S'$ and $\S''$. Up to replacing
$\varphi$ by $\varphi\psi$, we are reduced to the previous case
in which $\sigma$ is the identity, and so we are done.
\end{proof}

\section{Construction of invariant measures}\label{sec:family}



In this section we construct a family of locally finite, ergodic, $\Map(\S)$-invariant measures on the space of geodesic currents and recall the analogous construction by Lindenstrauss-Mirzakhani on the space of measured laminations. 

\subsection{Thurston measure and ergodicity}

We start by giving a brief description of a natural $\Map(S)$-invariant measure on the space of measured laminations, the Thurston measure, and refer the reader to \cite{thurstonnotes} and \cite{PennerHarer} for more details. Recall that the space of measured laminations $\ML_0(\S)$ 
supported in the interior of $\S$
has the structure of a piecewise linear manifold of dimension $N(\S)=-3\chi(\S)-n$ (where, as usual, $n$ is the number of boundary components of $\S$). 
It is also equipped with a $\Map(\S)$-invariant symplectic structure, giving rise to a $\Map(\S)$-invariant measure in the Lebesgue class; this is the {\it{symplectic Thurston measure}} $\Thu^{\mathrm{sympl}}$. Such measure has infinite
total mass, but it is locally finite and it satisfies the following scaling relation
$$\Thu^{\mathrm{sympl}}(L\cdot U)=L^{N(\S)}\cdot \Thu^{\mathrm{sympl}}(U)$$
for all Borel sets $U\subset\ML_0(\S)$ and all $L>0$.

A bit more concretely, the symplectic Thurston measure can be viewed the following way. Fix a maximal bi-recurrent train track $\tau$ on $\S$. The solution set $E(\tau)$ to the switch equations of $\tau$ is a $N(\S)$-dimensional rational cone in a Euclidean space and defines an open set in $\ML_0(\S)$. The restriction of the symplectic Thurston measure on
this open set can be identified to
the natural volume form on $E(\tau)$. 
In fact, the integer points in $E(\tau)$ are in one-to-one correspondence with the set of simple multi-curves with integral weights on $\S$ and we can obtain a multiple of the symplectic Thurston measure as the weak$^\star$ limit 
$$\Thu:=\lim_{L\to\infty}\frac{1}{L^{N(\S)}}\sum_{\gamma}\delta_{\frac{1}{L}\gamma}$$
where the sum is taken over all measured laminations $\gamma\in\ML_0(\S)$ 
corresponding to simple multi-curves with integral weights. 
In fact, the ratio between $\Thu$ and $\Thu^{\mathrm{sympl}}$ 
is a constant factor that only depends on the topology of $\S$
(see \cite{normalization}, for instance).
In what follows, we refer to $\Thu$ as the {\it{Thurston measure}}.

For us one of the most important features of the Thurston measure is the following result due to Masur \cite{Masur85}. 

\begin{thm}[Ergodicity of the Thurston measure]\label{Masur}
The Thurston measure $\Thu$ is ergodic on $\ML_0(\S)$ with respect to the action of $\Map(\S)$. 
\end{thm}

Finally, viewing $\ML_0(\S)$ as a (closed) subset of the space of currents $\Curr(\S)$, we can view $\Thu$ as a measure on $\Curr(\S)$ as well, assigning measure zero to any Borel set $U\subset\Curr_0(S)$ for which $U\cap\ML_0(S)=\emptyset$. Hence $\Thu$ is an example of a locally finite $\Map(\S)$-invariant ergodic measure on $\Curr(\S)$. In Sections \ref{sec:ML} and \ref{sec:other} below we will see further examples of such measures. 


\subsection{Classification of measures on $\ML$}\label{sec:ML}

We briefly discuss the complete classification of locally finite $\Map(\S)$-invariant ergodic measures on $\ML(\S)$ in the terminology used in \cite{LM08}. 

First, recall the definition of a complete pair from Section
\ref{sec:main-results}.

\begin{named}{Definition \ref{def:pair}}[Pairs and complete pairs]
Let $\R\subset \S$ be a subsurface and let $c\in\Curr(\S)$
be a current that standardly decomposes as a sum $c=\Gamma+\a$ of a simple multi-curve $\Gamma$ and an a-laminational $\a$.
The couple $(\R,c)$ is a {\it{pair}} if $\supp(c)$ and $\R$ are isotopically disjoint; such pair $(\R,c)$ is a {\it{complete pair}} if each boundary curve of $\R$ is  homotopic either to a boundary curve of $\S$, or
to a curve in the support of $\Gamma$, or
to a boundary curve of $\hull(\a)$.
\end{named}

In the special case when $c$ is a measured lamination, this definition agrees with the notion of a complete pair introduced in \cite{LM08}: if $c\in\ML(S)$ and $(\R,c)$ is a complete pair, then $c=\Gamma$ for some simple multi-curve $\Gamma$ with support $C$ 
such that $C\cap \R=\emptyset$ and every boundary curve of $\R$ is either a boundary curve of $\S$ or a curve in $C$. 

\begin{remark}
We underline that, in a standard pair $(R,\Gamma+\alpha)$,
the current $\Gamma+\a$ has no scc-free laminational part $\lambda$. 
The above Definition \ref{def:pair} is in fact tailored in such a way
that the couples $(R,\Gamma+\alpha)$
that appear in the orbit classification (Theorem \ref{thm:orbit-closure}) are indeed the complete pairs.
\end{remark}

Consider the map $\ML_{0}(\R)\to\ML(\S)$ defined by $\lambda\mapsto\lambda+\Gamma$. If $\R\neq\emptyset$, define $m^{(\R,\Gamma)}$ to be the push-forward of the Thurston measure through this map, which is then supported
on $\ML_{\R}(\S)+\Gamma$. In the case when $\R=\emptyset$, we define $m^{(\emptyset, \Gamma)}$ to be the Dirac measure $\delta_{\Gamma}$ in $\Curr(\S)$ supported on $\Gamma$. 
Now, define 
$$m^{[\R,\Gamma]} := \sum_{\varphi} m^{(\varphi(\R), \varphi(\Gamma))}$$
where the sum is taken over all $\varphi\in\Map(\S)/\stab(m^{(\R,\Gamma)})$. We note that, when $\Gamma=0$, we have $\R=\S$ and $m^{(\S, 0)}=m^{[\S,0]}=\Thu$.
On the other hand, 
$m^{[\emptyset, \Gamma]}$ is
the counting measure supported on the 
orbit of $\Gamma$.

Lindenstrauss-Mirzakhani \cite{LM08} and Hamenst{\"a}dt \cite{hamenstadt:measures} showed that, for any complete pair $(\R,\Gamma)$, the measures $m^{[\R,\Gamma]}$ are locally-finite, $\Map(\S)$-invariant and ergodic on $\ML(\S)$.
Moreover, the following classification result is proven in \cite{LM08}
(the result in \cite{hamenstadt:measures} is slightly weaker as the author does not show that the pair $(\R,\Gamma)$ must be complete in order for $m^{[\R,\Gamma]}$ to be locally finite).

\begin{thm}[Classification of ergodic invariant measures on $\ML$]\label{thm:LM}
Let $m$ be a locally finite $\Map(\S)$-invariant ergodic measure on $\ML(\S)$. Then $m$ is a multiple of $m^{[\R,\Gamma]}$ for a complete pair $(\R,\Gamma)$. 
\end{thm}

We will later use the following consequence of Theorem \ref{thm:LM},
that also follows from Proposition 8.5 in \cite{LM08}.

\begin{cor}\label{cor:R-in-R'}
For $\R\subseteq\hat{\R}$ consider
the measure on $\Curr(\hat{\R})$ defined as
\[
\sum_\psi m^{\psi(\R),\emptyset}
\]
where $\psi$ ranges over $\Map(\hat{\R},\pa\hat{\R})/\stab(\R)$.
If such measure is locally finite, then $\R=\hat{\R}$.
\end{cor}

Below we will see that any complete pair $(\R, c)$, where $c$ has a standard decomposition
of type $c=\Gamma+\alpha$, gives rise to a locally finite $\Map(\S)$-invariant ergodic measure on $\Curr(\S)$.


\subsection{Subsurface measures on $\Curr$}\label{sec:other}

Since every measure on $\ML(\S)$ can be viewed as a measure on $\Curr(\S)$, the measures  $m^{[\R,\Gamma]}$ defined above are locally finite $\Map(\S)$-invariant ergodic measures also on $\Curr(\S)$. However, one can easily construct other similar measures on $\Curr(\S)$. 

As a first example, consider a binding current $b\in\Cufill(\S)$ and consider the counting measure centered at the $\Map(\S)$-orbit of $b$, i.e. 
$\sum_{\varphi}\delta_{\varphi(b)}$
as $\varphi$ ranges over $\Map(\S)/\stab(b)$.
This defines a $\Map(\S)$-invariant measure by construction and it is also clear that it is ergodic. By Lemma \ref{lemma:div-orbit-binding-multicurve} it is also locally finite. 

As a second example, consider 
a (not necessarily binding) current $c\in\Curr(\S)$
with standard decomposition $c=\Gamma+\alpha$
so that $(\emptyset,c)$ is a complete pair.
%
We define $m^{(\emptyset,c)}:=\delta_c$ and
$
m^{[\emptyset,c]}:=\sum_\varphi m^{(\emptyset,\varphi(c))}
$
as $\varphi$ ranges over all elements of $\Map(\S)/\stab(c)$.
Clearly, $m^{[\emptyset,c]}$ agrees with
the counting measure on the 
orbit of $c$,
which is closed and discrete by Theorem \ref{thm:orbit-closure}.
%
%
It follows that $m^{[\emptyset,c]}$ is a locally finite measure on $\Curr(\S)$ for any current $c$ for which $(\emptyset, c)$ is a complete pair. We record this observation below. 

\begin{lemma}[Locally finite ergodic invariant counting measures]\label{lemma:counting}
Let $\Gamma$ be a simple multi-curve
isotopically disjoint from the a-laminational current $\alpha$.
Then the counting measure $m^{[\emptyset,c]}$ centered at the $\Map(\S)$-orbit of $c=\Gamma+\alpha$
is a locally finite, mapping class group invariant, ergodic measure on $\Curr(\S)$.
\end{lemma}

In the remainder of this section
we consider the general case of a subsurface $\R$ of $\S$
and a current $c\in\Curr(\S)$ and assume that $(\R,c)$ is a pair.
By definition, $c$ admits a standard decomposition $c=\Gamma+\a$,
where $\Gamma$ a simple multi-curve, $\a$ is a-laminational
and the loci $\R$, $C=\supp(\Gamma)$ and $\A=\supp(\a)$
are disjoint up to isotopy.

As in the introduction, and following the construction by Lindenstrauss-Mirzakhani, we define measures $m^{(\R,c)}$ in the following way. If $\R=\emptyset$, we let $m^{(\emptyset,c)}=\delta_{c}$ as above.
If $\R\neq\emptyset$, we let $m^{(\R,c)}$ denote the push-forward of the Thurston measure through the map $\ML_0(\R)\to\Curr(\S)$ defined by $\lambda\mapsto\lambda+c$, which is then supported
on $\ML_{\R}(\S)+c$.
Finally, we define the {\it{subsurface measure $m^{[\R,c]}$ of type $[\R,c]$}} as
$$m^{[\R,c]} := \sum_{\varphi} m^{(\varphi(\R), \varphi(c))}$$
as $\varphi$ ranges over $\Map(\S)/\stab(m^{(\R,c)})$. 

We observe that $\stab(m^{(\R,c)})=\stab(\R)\cap \stab(\Gamma)\cap\stab(\a)$
and that $\stab(\Gamma)\supset \Map(\S,C)$ and $\stab(\a)\supset\Map(\S,\A)$
are finite-index subgroups. Thus,
$\stab(m^{(\R,c)})$ contains $\stab(\R,C\cup\A)$ as a finite-index subgroup.
%
By construction, $m^{[\R,c]}$ is $\Map(\S)$-invariant.

%

Recall from the introduction that the pair $(\R,c)$ is {\it{complete}}
if each boundary curve of $\R$ is either a boundary curve of $\S$ or of $\A$, or a component of $C$. The following lemma highlights the importance of the completeness
property for a pair.

\begin{lemma}[Local finiteness of translates of $\ML_\R$]\label{lemma:local-finiteness}
Let $(\R,c)$ be a pair in $\S$.
Then the following are equivalent:
\begin{itemize}
\item[(a)]
the pair $(\R,c)$ is complete;
\item[(b)]
the quotient $\stab(c)/\stab(\R)\cap\stab(c)$ is finite;
\item[(c)]
the collection of subsets $\ML_{\varphi(\R)}(\S)+\varphi(c)$ of $\Curr(\S)$,
as $\varphi$ ranges over $\Map(\S)/\stab(\R)\cap\stab(c)$,
is locally finite.
\end{itemize}
\end{lemma}
\begin{proof}
Since $\stab(c)/\Map(\S,C\cup \A)$ is finite by Lemma \ref{lemma:stabilizer}(b), it is easy to see that (a) is equivalent to (b).


If $\stab(c)/\stab(\R)\cap\stab(c)$ is infinite, then
the collection of $\ML_{\varphi(\R)}(\S)+\varphi(c)$
as $\varphi$ ranges over $\stab(c)/\stab(\R)\cap\stab(c)$
is not locally finite at $c\in\Curr(\S)$. This shows that (c) implies (b).

Finally, suppose that (b) holds and let
$\varphi_{i}(\lambda_i)+\varphi_{i}(c)\rar \lambda_\infty+c_\infty$,
where $\varphi_i\in\Map(\S)$, $\lambda_i\in\ML_\R(\S)$
and $\lambda_\infty+c_\infty=\lambda_\infty+\Gamma_\infty+\a_\infty$ is the standard decomposition of the limit current.
It is enough to show that, up to extracting a subsequence, all
$\varphi_i$ belong to $\stab(c)\cap\stab(\R)$. 

Fix $b'$ a binding multi-curve on $\S$.
Since $\iota(\varphi_{i}(\lambda_i+c),b')\rar\iota(\lambda_\infty+c_\infty,b')$,
the quantity $\iota(\varphi_{i}(c),b')$ is uniformly bounded.
By Lemma \ref{lemma:Gamma+a-bound}, up to subsequences, we can assume that
$[\varphi_i]\in\Map(\S)/\Map(\S,C\cup \A)$ is constant.
Up to applying $\varphi_{1}^{-1}$ to all involved currents,
we can assume that all $\varphi_i\in\Map(\S,C\cup \A)\subseteq\stab(c)$.
By (b) we can then extract a subsequence such that
all $\varphi_i$ satisfy $\varphi_i(\R)=\R$.
This shows that (b) implies (c).
\end{proof}

We will now show that the subsurface measures $m^{[\R,c]}$ are locally finite and ergodic, provided the pair $(\R,c)$ is complete.

%

\begin{prop}[Local finiteness of subsurface measures]\label{prop:locally-finite}
Let $(\R, c)$ be a pair. Then $m^{[\R,c]}$ is a $\Map(\S)$-invariant ergodic measure on $\Curr(\S)$. Moreover, $m^{[\R,c]}$ is locally finite if and only if $(\R, c)$ is a complete pair. 
\end{prop}
\begin{proof}
%
Mapping class group invariance of $m^{[\R,c]}$
follows from the above discussion.
As for the ergodicity, note that the support of $m^{(\R,c)}$ is
$\ML_{\R}(\S)+c$ and that the support of $m^{[\R,c]}$
is the union of all translates $\varphi\cdot(\ML_{\R}(\S)+c)$.
Thus, $m^{[\S',c]}$ is ergodic if and only if $\stab(\ML_{\R}(\S)+c)$
acts ergodically on $\ML_{\R}(\S)+c$ with respect to the measure $m^{(\R,c)}$.


Recall that $\Map(\R)$ acts ergodically on $\ML_0(\R)$ with respect to the Thurston measure (Theorem \ref{Masur})
and so $\Map(\S,\S\setminus\R)$ acts ergodically on
$\ML_\R(\S)+c$.
Since $\stab(\ML_\R(\S)+c)$ contains $\Map(\S,\S\setminus\R)$,
it follows that $\stab(\ML_\R(\S)+c)$ acts ergodically
on $\ML_\R(\S)+c$ too. We conclude that $m^{[\R,c]}$ is 
ergodic for the action of $\Map(\S)$.

It remains to show that $m^{[\R,c]}$ is locally finite if and only if $(\R,c)$ is a complete pair. 

Suppose first that $(\R,c)$ is complete. 
By Lemma \ref{lemma:local-finiteness}(c),
the union of all translates
$\varphi\cdot(\ML_{\R}(\S)+c)$ 
as $\varphi$ ranges over $\Map(\S)/\stab(\R)\cap\stab(c)$
is locally finite.
Since the Thurston measure on $\ML_0(\R)$ is locally finite, so is $m^{[\R,c]}$.

Suppose conversely that $m^{[\R,c]}$ is locally finite.
If $\R=\emptyset$, the conclusion follows from Lemma \ref{lemma:Gamma+a-bound}. Assume then
$\R\neq\emptyset$ and
let $\hat{\R}$ be the union of the components of $\S\setminus(C\cup\A)$ that
intersect $\R$, so that $(\hat{\R},c)$ is a complete pair.
The pull-back of $m^{[\R,c]}$ via the map $\ML_0(\hat{\R})\rar \ML_{\hat{\R}}(\S)+c$
is locally finite and it is greater or equal than
\[
\sum_\psi m^{\psi(\R),\emptyset}
\]
where $\psi$ ranges over $\Map(\hat{\R},\pa\hat{\R})/\stab(\R)$.
By Corollary \ref{cor:R-in-R'}, it follows that $\R=\hat{\R}$ and so $(\R,c)$ is a complete pair.
\end{proof}

In the next section we will see that any locally finite, $\Map(S)$-invariant, ergodic measure on $\Curr(S)$ must be a positive multiple of $m^{[\R,c]}$ for some complete pair $(\R, c)$. 

\section{Classification of invariant measures}\label{sec:proof}




\subsection{Measures on $\Cufh$}

Suppose $m$ is a locally finite, $\Map(\S)$-invariant, ergodic measure on $\Curr(\S)$. 
If $m(\{0\})>0$, then $m$ is a positive multiple of $\delta_{\{0\}}$, the Dirac measure centered at $0$.
From now on, we therefore assume that $m(\{0\})=0$ and so $m$ is the push-forward of a measure
on $\Curr(\S)\setminus\{0\}$.

Note that since $\Cufh(S)$ is $\Map(S)$-invariant, it follows by ergodicity that if $m(\Cufh(S))>0$ then
$\Cufh(\S)$ has in fact full $m$-measure and so we can interpret $m$ as
(the push-forward of) a measure on $\Cufh(\S)$. The classification of such measures
is provided by the following proposition, which partially relies on \cite{LM08}.

\begin{prop}[Ergodic measures supported on $\Cufh$]\label{fullhull}
Suppose $\S$ is connected and $m$ is a locally finite, $\Map(\S)$-invariant, ergodic 
measure on $\Curr(\S)$ such that $m(\Cufh(\S))>0$. Then exactly one of the following holds:
\begin{itemize}
\item[(i)]
$\MLfh(\S)$ has full $m$-measure and
$m$ is a positive multiple of a translate of the Thurston measure given by $m^{[\S,\Gamma]}$ where $\Gamma$ is any simple multi-curve with
support contained inside $\pa\S$, or
\item[(ii)]
$\Cufill(\S)$ has full $m$-measure and
$m$ is a positive multiple of the Dirac measure $m^{[\emptyset,b]}$ for some binding current $b\in\Cufill(\S)$.
\end{itemize}  
\end{prop}

\begin{proof}
Recall that, by Theorem \ref{MLfill-intro}, the space
$\Cufh(\S)$ is the union of $\MLfh(S)$ and $\Cufill(\S)$
and these sets are disjoint, Borel, and $\Map(S)$-invariant.

If $\MLfh(\S)$ has full measure, then it follows from \cite[Theorem 7.1]{LM08} that $m$ is a multiple of the Thurston measure on 
a translate of $\ML_0(\S)$ and so we are in case (i).
Otherwise, Lemma \ref{lemma:hausdorff} below
shows that we are in case (ii), since $G=\Map(\S)$ acts in a properly discontinuous way
on $X=\Cufill(\S)$ by Proposition \ref{properly}.
\end{proof}

The following lemma is well-known; we include it for completeness.

\begin{lemma}\label{lemma:hausdorff}
Let $X$ be a locally compact Hausdorff topological space 
and let $G$ be a discrete group that acts properly discontinuously
on $X$ via self-homeomorphisms.
Then a locally finite $G$-invariant ergodic measure $m$ on $X$
is a positive multiple of the counting measure on a $G$-orbit.
\end{lemma}
\begin{proof}
Let $x\in X$ be a point in the support of 
the locally finite, ergodic, $G$-invariant measure
$m$ on $X$.
It is enough to show that, if $x'\notin G\cdot x$, then
$x'$ does not belong to the support of $m$.

Note that, since $G$ acts properly discontinuously
and $X$ is Hausdorff and locally compact, the quotient $X/G$ is Hausdorff.
Moreover, $[x]\neq [x']$ as points of $X/G$. 
Thus, there exist disjoint open neighbourhoods $U,U'\subset X/G$ of $[x]$ and $[x']$, respectively,
and we denote by $\tilde{U},\tilde{U'}$ their preimages in $X$, which are disjoint, open
and $G$-invariant. Since $x$ belongs to the support of $m$, we must have
$m(\tilde{U})>0$ and so $m(\tilde{U}')=0$ by ergodicity.
As a consequence, the support of $m$ is contained inside $X\setminus\tilde{U}'$ and
so, in particular, it does not contain $x'$.
\end{proof}

\subsection{Classifying ergodic invariant measures on $\Curr$}

Before proving our main result, we recall the following useful lemma
by Lindenstrauss-Mirzakhani \cite[Lemma 8.4]{LM08}.

\begin{lemma}[Ergodic action on a product]\label{lemma:product-measures}
Let $X'$ and $X''$ be locally compact, second countable, metric spaces and let $G'$ and $G''$ be discrete, countable groups,
acting continuously on $X'$ and $X''$ respectively. Then any locally finite $(G'\times G'')$-invariant ergodic measure $m$ on $X'\times X''$ is of the form $m=m'\otimes m''$, where $m'$ (resp. $m''$) is a locally finite $G'$-invariant ergodic measure on $X'$ (resp. $G''$-invariant ergodic measure on $X''$). 
\end{lemma}

\begin{remark}
Let $\R,\A\subset\S$ be subsurfaces.
By the work of Thurston, $\ML_0(\R)$ is homeomorphic
to a finite-dimensional Euclidean space.
As a consequence, $\ML_0(\R)$ is metrizable
and it has a countable exhaustion by compact subsets.
Bonahon showed (Theorem \ref{thm:topological}) that
the same properties hold for $\Curr_0(\A)$.
Thus,  both $\ML_0(\R)$ and $\Curr_0(\A)$ are metrizable, locally compact and second countable.
Note now that the locus $\ML_0(\R)^*$ of measured laminations whose support intersects all
connected components of $\R$ is open inside $\ML_0(\R)$; in particular,
if $\R$ is connected, then $\ML_0(\R)^*=\ML_0(\R)\setminus\{0\}$.
Moreover, $\Cufill_0(\A)$ is open inside $\Curr_0(\A)$ by Corollary \ref{cor:binding-open}.
It follows that $\ML_0(\R)^*$ and $\Cufill_0(\A)$ are locally compact, metrizable
and second countable as well.
\end{remark}

We can now restate Theorem \ref{thm:main} and
finally complete the classification of all $\Map(\S)$-invariant
ergodic measures on $\Curr(\S)$.

\begin{named}{Theorem \ref{thm:main}}[Classification of locally finite invariant measures on $\Curr$]
The measure $m^{[\R,c]}$ on $\Curr(\S)$ is ergodic, $\Map(\S)$-invariant and locally finite
for every complete pair $(\R,c)$. Moreover, if $m$ is a locally finite, $\Map(\S )$-invariant, ergodic measure on $\Curr(\S )$,
then $m$ is a positive multiple of $m^{[\R,c]}$ for some complete pair $(\R,c)$.
\end{named}

\begin{proof}
By Proposition \ref{prop:locally-finite} every
$m^{[\R,c]}$ associated to a complete pair $(\R,c)$ is
$\Map(\S)$-invariant, ergodic and locally finite:
this is exactly the first claim.

In order to prove the second claim, consider
a $\Map(\S)$-invariant, locally finite, ergodic measure
$m\neq 0$ on $\Curr(\S)$. We want to show that
$m$ is a positive multiple of $m^{[\R,c]}$, for some
pair $(\R,c)$. By Proposition \ref{prop:locally-finite}
it will automatically follow that $(\R,c)$ is complete.

We recall that, by Corollary \ref{cor:invariant-partition}, the space of currents
can be decomposed into a union of the $\Map(\S)$-invariant, disjoint Borel subsets
$\Curr_{[\R,C,\A]}(\S)$.

By ergodicity, there exists a unique triple $(\R,C,\A)$ such that
$\Curr_{[\R,C,\A]}(\S)$ has full $m$-measure.
Thus, 
it is enough to analyze the restriction 
$m_G$ of $m$ to a single component
$\Curr_{(\R,C,\A)}(\S)$ of $\Curr_{[\R,C,\A]}(\S)$, which is
ergodic with respect to the stabilizer
$G=\stab(\R)\cap\stab(C)\cap\stab(\A)$ of $\Curr_{(\R,C,\A)}(\S)$.
Indeed, the conclusion will follow by $\Map(\S)$-invariance. 

Let $H\subset G$ be the finite-index subgroup of elements
that send every component of $\R,\A,C$ and of $\pa \R,\pa\A$ to itself.
Then $m_G$ can be written as $m_G=\frac{1}{|G/H|}\sum_{g\in G/H}g_* m_H$,
where $m_H$ is an $H$-invariant ergodic measure on $\Curr_{(\R,C,\A)}(\S)$.
Thus, it is enough to show that $m_H$ is a multiple of
the restriction of $m^{[\R,c]}$ to $\Curr_{(\R,C,\A)}(\S)$
for some $c=\Gamma+\a$, with $\supp(\Gamma)=C$ and 
$\a$ an a-laminational current that binds $\A$.

Recall that $\Curr_{(\R,C,\A)}(\S)$ is the product of the three 
factors $\MLfh_{\R}(\S)$, $\Cufh_C(\S)$ and $\Cufill_{\A}(\S)$
and that the push-forward maps identify
$\MLfh_0(\R)$ to $\MLfh_{\R}(\S)$ and $\Cufill_0(\A)$ to $\Cufill_{\A}(\S)$.
Since $H$ acts trivially on $C$, we can write $m_H=\check{m}\otimes\delta_\Gamma$,
where $\Gamma$ is a simple multi-curve with support $C$ and $\check{m}$
can be viewed as a locally finite $\left(\Map(\R)\times\Map(\A)\right)$-invariant
ergodic measure on $\MLfh_0(\R)\times\Cufill_0(\A)$ of full support, and so in particular on $\ML_0(\R)^*\times\Cufill_0(\A)$.

%
%

Applying Lemma \ref{lemma:product-measures} to
$X'=\ML_0(\R)^*$ and $X''=\Cufill_0(\A)$ 
with $G'=\Map(\R)$ and $G''=\Map(\A)$,
we obtain the decomposition $\check{m}=\check{m}_R\otimes\check{m}_A$,
where $\check{m}_R$ is a $\Map(\R)$-invariant ergodic
measure of full support on $\ML_0(\R)^*$ and $\check{m}_A$ is
a $\Map(\A)$-invariant ergodic measure of full support on $\Cufill_0(\A)$.

Since $\ML_0(\R)^*=\bigoplus_i \ML_0(\R_i)^*$ is acted on
by $\Map(\R)=\prod_i\Map(\R_i)$ and
$\Cufill_0(\A)=\bigoplus_j \Cufill_0(\A_j)$
is acted on by $\Map(\A)=\prod_j\Map(\A_j)$,
we can iteratively apply Lemma \ref{lemma:product-measures}
and we obtain that $\check{m}_{\R}=\bigotimes_i \check{m}_{\R_i}$ and
$\check{m}_{\A}=\bigotimes_j \check{m}_{\A_j}$,
where $\check{m}_{\R_i}$ is a locally finite $\Map(\R_i)$-invariant ergodic measure of full support
on $\ML_0(\R_i)^*$ and $\check{m}_{\A_j}$ is a locally finite $\Map(\A_j)$-invariant ergodic
measure of full support on $\Cufill_0(\A_j)$.

It is also easy to see that each $\check{m}_{\R_i}$ is indeed the push-forward of a measure on $\MLfh_0(\R_i)$ of full support.
By Lemma \ref{fullhull},
it follows that $\check{m}_{R_i}$ is a multiple of the Thurston measure on $\ML_0(\R_i)^*$
and that $\check{m}_{\A_j}$ is a multiple of the counting measure on the $\Map(\A_j)$-orbit
of some $\a_j\in\Cufill_0(\A_j)$.
We have then obtained that $m_H$ is a multiple of $m^{(\R,c)}$ with
$c=\Gamma+\a$ and $\a=\sum_j \a_j\in\Cufill_0(\A)$,
and so the proof is complete.
%
\end{proof}

\section{Homogeneous invariant measures}\label{sec:homogeneous}

Consider the natural action of $\RR_+$ on $\Curr(\S)$
by multiplication.

\begin{defi}[Homogeneous measures]
A measure $m$ on $\Curr(\S)$ is {\it{$d$-homogeneous}} for some $d\in\RR$
if $m(t\cdot U)=t^d\cdot m(U)$ for all Borel subsets $U$ of $\Curr(\S)$
and all $t\in\RR_+$.
\end{defi}

A $\Map(\S)$-invariant, $d$-homogeneous measure $m$ on $\Curr(\S)$
is ergodic (as a $d$-homogeneous measure) if
it is ergodic for the action
of $\Map(\S)\times \RR_+$ on $\Curr(\S)$
defined as
$\left((\varphi,t)\cdot m\right)(U):=t^d\cdot m(\varphi(U))$.

Lindenstrauss-Mirzakhani \cite[Proposition 8.2]{LM08}
showed that, if a measure $m$
is locally finite, $\Map(\S)$-invariant, $d$-homogeneous 
and supported on $\ML(\S)$, then $d\geq N(\S)$.

The aim of this section is to give an almost complete
classification of $d$-homogeneous, $\Map(\S)$-invariant,
ergodic measures on $\Curr(\S)$.

\subsection{Construction of the homogeneous measures}

Let $(\R,c)$ be a complete pair in $\S$, that is
$\R\subset\S$ is a (possibly empty) subsurface and $c=\Gamma+\alpha$
is the sum of a multi-curve $\Gamma$ and an a-laminational current
$\alpha$
such that $\R$, $C=\supp(\Gamma)$ and $\A=\hull(\alpha)$ are disjoint.
It will be useful to decompose $c$ as $c=c_{\pa\R}+c'$, where
$c_{\pa\R}$ is supported on $\pa\R$ and $\supp(c')\cap\R=\emptyset$.

If $\R\neq\emptyset$ has genus $g(\R)$ and $n(\R)$ boundary components, we define by
$N(\R):=6g(\R)-6+2n(\R)$
and by $N'(\R):=6g(\R)-6+3n(\R)$. 
If $\R=\emptyset$, we let $N(\emptyset)=N'(\emptyset)=0$.

For every $d\in\RR$ 
consider the measure
\[
m^{(\R,c)}_d:=
\begin{cases}
m^{(\S,0)} & \text{if $(\R,c)=(\S,0)$ and $d=N(\S)$}\\
\int_0^{+\infty}t^{d-N(\R)-1}m^{(\R,tc)}dt & \text{if $c\neq 0$}
\end{cases}
\]
on $\Curr(\S)$.
Notice that we do not define $m^{(\S,0)}_d$ with
$d\neq N(\S)$.

In order to study the local finiteness of the measures $m^{(\R,c)}_d$,
we  fix an auxiliary hyperbolic metric $h$ on $\S$
and we let $\ell_h:\Curr(\S)\rar\RR$ be the proper continuous length function
attached to $h$ as in Remark \ref{rmk:hyperbolic-length-function}.
Moreover, we denote by $\Thu^\R$ the Thurston measure
on $\Curr(\S)$ which is supported on $\ML_{\R}(\S)$,
and by $B_h$ the $h$-unit ball of currents 
$B_h:=\{c\in\Curr(\S)\,|\,\ell_h(c)\leq 1\}$. 


Analogously to what is done in Section \ref{sec:family},
let
\[
m^{[\R,c]}_d:=\sum_{\varphi}m_d^{(\varphi(\R),\varphi(c))}
\]
as $\varphi$ ranges over $\Map(\S)/\stab(\R,c)$, where we
recall that $\stab(\R,c)$ is a finite-index subgroup of $\stab(c)$, because $(\R,c)$ is a complete pair.

The main result of this section is the following classification theorem.

\begin{named}{Theorem \ref{thm:locally-finite-homogeneous}}[Locally finite invariant homogeneous measures]
Every locally finite $\Map(\S)$-invariant $d$-homogeneous ergodic measure
on $\Curr(\S)$ is a positive multiple of one of the following:
\begin{itemize}
\item[(i)]
the Thurston measure $m^{[\S,0]}_{N(\S)}=\Thu^{\S}$
\item[(ii)]
the measure $m^{[\R,c]}_d$ with $c\neq 0$ and $d>N(\S)$ large enough.
\end{itemize}
In part (ii) every $d>N(\S)+N(\R)$ works.
\end{named}

An immediate consequence of the above result, we have the following useful observation.

\begin{cor}[Invariant measures of homogeneity $N(\S)$]\label{cor:homog-N}
A locally finite $\Map(\S)$-invariant $N(\S)$-homogeneous measure
on $\Curr(\S)$ is a multiple of the Thurston measure $\Thu^{\S}$.
\end{cor}

We will prove Theorem \ref{thm:locally-finite-homogeneous}
through a series of lemmas in the next subsection.

\subsection{Local finiteness of the measures $m_d^{[\R,c]}$}

Before studying the invariant measures $m_d^{[\R,c]}$, we
analyze the (non-invariant) homogeneous measures $m_d^{(\R,c)}$.

\begin{lemma}[Local finiteness of $m_d^{(\R,c)}$]\label{lemma:hom-measures}
The measure $m^{(\R,c)}_d$ on $\Curr(\S)$ is $d$-homogeneous.
Moreover,
$m^{(\R,c)}_d$ with $c\neq 0$ 
is locally finite if and only if $d>N(\R)$,
in which case $m^{(\R,c)}_d(\ell_h\leq L)
=\Thu^\R(B_h)\,\frac{N(\R)! }{d(d-1)\cdots(d-N(\R))}\,
\frac{L^d}{\ell_h(c)^{d-N(\R)}}$.
\end{lemma}
\begin{proof}
The $d$-homogeneity of $m^{(\R,c)}_d$ is clear by construction, since
$\Thu^\R$ is $N(\R)$-homogeneous.

As for the second claim, we fix a hyperbolic metric $h$ on $\S$
and an $L>0$, and we want to determine for which $d$
the quantity $m^{(\R,c)}_d(\ell_h\leq L)$ is finite.
Clearly,
\[
m^{(\R,tc)}(\ell_h\leq L) =\Thu^{\R}(\ell_h\leq L-t\ell_h(c))
=\begin{cases}
\Thu^R(B_h)\cdot (L-t\ell_h(c))^{N(\R)} & \text{if $L>t\ell_h(c)$}\\
0 & \text{if $L\leq t\ell_h(c)$.}
\end{cases}
\]
Thus,
\[
m^{(\R,c)}_d(\ell_h\leq L)=\Thu^\R(B_h) \int_0^{L/\ell_h(c)}t^{d-N(\R)-1}(L-t\ell_h(c))^{N(\R)}dt
\]
is finite if and only if $d>N(\R)$. In this case,
\[
m^{(\R,c)}_d(\ell_h\leq L)
=\Thu^\R(B_h) \,\frac{N(\R)! }{d(d-1)\cdots(d-N(\R))}\,
\frac{L^d}{\ell_h(c)^{d-N(\R)}}.
\]
is obtained just integrating by parts.
\end{proof}

In the following proposition we analyze local finiteness
of homogeneous measures of non-Thurston type.

\begin{prop}[Finiteness of invariant homogeneous measures of non-Thurston type]\label{prop:finiteness}
For a complete pair $(\R,c)\neq (\S,0)$  the following holds:
\begin{itemize}
\item[(i)]
The measure $m^{[\R,c]}_d$ is $d$-homogeneous,
$\Map(\S)$-invariant and ergodic.
\item[(ii)]
For every $d> N(\S)+N(\R)$, the measure
$m^{[\R,c]}_d$ is locally finite.
\item[(iii)]
For every $d\leq N(\S)$, the measure $m^{[\R,c]}_d$ is not locally finite.
\end{itemize}
\end{prop}

Now we show how Theorem \ref{thm:locally-finite-homogeneous} follows from the above result.

\begin{proof}[Proof of Theorem \ref{thm:locally-finite-homogeneous}]
The measure $\Thu$ and the measures $m_d^{[\R,c]}$ on $\Curr(\S)$
are homogeneous, ergodic and $\Map(\S)$-invariant. We wish to show that these are the only ones.

Let $m_d$ be a locally finite, $d$-homogeneous, $\Map(\S)$-invariant measure on $\Curr(\S)$.
Using the ergodic decomposition of locally finite (not necessarily homogeneous) $\Map(\S)$-invariant measures, we can write
\[
m_d=r\cdot\Thu+\sum_{\R\subsetneq\S}\int_{\mathcal{C}_\R}m^{[\R,c]}\,\mu_\R(c) 
\]
where $\mathcal{C}_\R$ is the space of currents $c$ of type $c=\Gamma+\alpha$
such that $(\R,c)$ is a complete pair, $\mu_\R$ is a measure on $\mathcal{C}_\R$ and $r\in\RR_{\geq 0}$.
Clearly, we must have $r=0$ unless $d=N(\S)$.

Fix an auxiliary hyperbolic metric $h$ on $\S$ and let $\mathcal{C}^1_\R$ be the subset of $\mathcal{C}_\R$
consisting of currents of $h$-length $1$.
The map $\RR_+\times \mathcal{C}^1_\R\rar \mathcal{C}_\R$ given by $(t,c)\mapsto tc$ is clearly a homeomorphism.
For every $\R\subsetneq\S$, define the measure $\mu^1_\R$ on $\mathcal{C}^1_\R$ as $\mu^1_\R(U):=d\cdot \mu_\R(\hat{U})$ for all Borel subsets $U\subseteq\mathcal{C}^1_\R$, where $\hat{U}:=(0,1)\cdot U$.
Since $m_d$ is $d$-homogeneous, it can be rewritten as
\begin{align*}
m_d &=r\cdot\Thu+\sum_{\R\subsetneq\S}\int_{\mathcal{C}^1_\R} \left(\int_0^{+\infty} m^{[\R,tc]}t^{d-N(\R)-1}dt\right)\mu^1_\R(c)=\\
&=r\cdot\Thu+\sum_{\R\subsetneq\S}\int_{\mathcal{C}^1_\R}m_d^{[\R,c]}\,\mu^1_\R(c).
\end{align*}
Hence, locally finite, $d$-homogeneous, ergodic, $\Map(\S)$-invariant measures are multiples either of $\Thu$
(if $d=N(\S)$), or of measures of type $m_d^{[\R,c]}$.
The result now follows from the analysis of the local finiteness of the measures $m_d^{[\R,c]}$ in Proposition \ref{prop:finiteness}.
\end{proof}

Now we turn to Proposition \ref{prop:finiteness},
which relies on the following estimate, which will be proven
in Appendix \ref{app}.

\begin{lemma}[Volume of the unit ball in $\ML_0(\R)$]\label{lemma:estimate-ball}
Let $(\R,h_\R)$ be a hyperbolic surface with geodesic boundary with 
systole $\sys(\R)\geq s>0$.
Then
\[
\frac{\hat{k}}{\ell_{h_\R}(\pa\R)^{N(\R)}}
<\Thu^\R(B_{h_\R})<\hat{k}'.
\]
for suitable constants $\hat{k},\hat{k}'>0$ that depend only on $s$ and
the topology of $\R$.
\end{lemma}


\begin{proof}[Proof of Proposition \ref{prop:finiteness}]
Part (i) is immediate by construction and by the ergodicity of
$m^{[\R,c]}$ proven in Theorem \ref{thm:main}.

As for (ii) and (iii), fix an auxiliary hyperbolic metric $h$ on $\S$ and an $L>0$.
We want to determine for which $d$ the quantity 
$m^{[\R,c]}_d(\ell_h\leq L)$ is finite.

Let $s:=\sys(\S,h)>0$ and let $\varphi\in\Map(\S)$.
Observe preliminarily that $\sys(\varphi(\R))\geq s$ and that
the ratio $\ell_h(\varphi(\pa \R))/\ell_h(\varphi(c_{\pa\R}))$ can be bounded from above
and from below by positive constants that are independent of $\varphi$
(they can be chosen to depend only on $c_{\pa\R}$, on the boundary length of $\S$
and on $s$).
Moreover, applying Lemma \ref{lemma:bound-orbit-c} with $L=q^u$, we obtain a constant $v>0$
that depend only on $S$, $h$ and $c$ such that
\[
{\textstyle\frac{1}{v}} \cdot q^{uN(\S)}<\#\left\{\varphi\in\Map(\S)/\stab(\R,c)\,|\,\ell_h(\varphi(c))\in[q^u,q^{u+1})\right\}<v \cdot q^{uN(\S)}
\]
for all $u$.

By Lemma \ref{lemma:hom-measures},
we need to study the finiteness of the following
\[
m^{[\R,c]}_d(\ell_h\leq L) =
\frac{N(\R)! }{d(d-1)\cdots(d-N(\R))}
L^d
\sum_\varphi 
\frac{\Thu^{\varphi(\R)}(B_h)}{\ell_h(\varphi(c))^{d-N(\R)}}
\]
By Lemma \ref{lemma:estimate-ball},
the quantity $m^{[\R,c]}_d(\ell_h\leq L)/L^d$
can be bounded from below as
\[
m^{[\R,c]}_d(\ell_h\leq L)/L^d
\geq 
\sum_{\varphi}
\frac{k_1}{\ell_h(\varphi(c_{\pa\R}))^{N(\R)}\ell_h(\varphi(c))^{d-N(\R)}}
\geq
 \sum_{\varphi}
 \frac{k_2}{\ell_h(\varphi(c))^d}
\]
for suitable constants $k_1,k_2$ independent of $\varphi$.
Hence,
\[
m^{[\R,c]}_d(\ell_h\leq L)/L^d\geq
\sum_{l\in\ZZ_+}
\sum_{\varphi\in\Phi_u}
 \frac{k_2}{\ell_h(\varphi(c))^d}
 \geq
 k'_2
\sum_{u\in\ZZ_+} 
q^{u(N(\S)-d)}
\]
for a suitable constant $k'_2$, where
$\Phi_u=\{\varphi\in\Map(\S)\,|\,\ell_h(\varphi(c))\in [q^u,q^{u+1})\}$.
It follows that, if $m^{[\R,c]}_d(\ell_h\leq L)$ is finite, then
$d>N(\S)$ and so (iii) is proven.

Analogously,
using again Lemma \ref{lemma:estimate-ball},
the quantity $m^{[\R,c]}_d(\ell_h\leq L)/L^d$ can be bounded
from above as
\[
m^{[\R,c]}_d(\ell_h\leq L)/L^d
\leq \sum_{\varphi}
\frac{k_3 }{\ell_h(\varphi(c))^{d-N(\R)}}
\]
for a suitable constant $k_3$ independent of $\varphi$.
Thus,
\[
m^{[\R,c]}_d(\ell_h\leq L)/L^d
\leq
\sum_{u\in\ZZ_+}
\sum_{\varphi\in\Phi_u}
 \frac{k_3}{\ell_h(\varphi(c))^{d-N(\R)}}\leq
k'_3
\sum_{u\in\ZZ_+} 
q^{u(N(\S)+N(\R)-d)}
\]
for a suitable $k'_3$. Since the last series is convergent for $d>N(\S)+N(\R)$,
we obtain (ii).
\end{proof}

\subsection{Counting curves} We conclude by exploring an application to Theorem \ref{thm:locally-finite-homogeneous}. For simplicity of exposition, from here on we assume that $\S$ is a closed surface of genus $g\geq 2$. As was mentioned in the introduction, one of the motivations to studying invariant measures on $\Curr(\S)$ is the use of certain such measures as a tool for counting curves on surfaces. More precisely, let us fix $\gamma$ to be a closed curve on $S$ and consider, for $L>0$, the family of {\it{curve-counting measures}} on $\Curr(S)$ defined by
\begin{equation*}
m^{[\gamma/L]} := \frac{1}{L^{N(\S)}}m^{[\emptyset,\gamma/L]}=\frac{1}{L^{N(\S)}}
\sum_{\gamma'\in\Map(S)\cdot\gamma}\delta_{\frac{1}{L}\gamma'}.
\end{equation*}
Let $f:\Curr(\S)\to\mathbb{R}_{\geq0}$ be any continuous, $1$-homogeneous function (for instance, a hyperbolic length function) and $B_f:=\{c\in\Curr(\S)\,|\,f(c)\leq1\}$. 
We include the following easy observation without proof.

\begin{remark}
The ball $B_f$ is closed. Moreover, every $d$-homogeneous measure $m_d$ (with $d\neq 0$) satisfies $m_d(\pa B_f)=0$.
\end{remark}

The reason for considering the ball $B_f$ is that
$$\frac{1}{L^{N(\S)}}\#\{\gamma'\in\Map(\S)\cdot\gamma\,\vert\,f(\gamma')\leq L\} = m^{[\gamma/L]} (B_f).$$
Hence, counting curves in a $\Map(\S)$-orbit of bounded length reduces to understanding the asymptotics of the curve-counting
measures $m^{[\gamma/L]}$ as $L\to\infty$. 

\begin{remark}[How to detect precompactness of curve-counting measures]\label{rmk:compactness}
Note that $B_h$ is compact and that probability measures on a compact space (such as $B_h$) are weak$^\star$-compact.
Hence, 
if $\limsup_{L\rar\infty}\frac{1}{L^{N(\S)}}\#\{\gamma'\in\Map(\S)\cdot\gamma\,\vert\,\ell_h(\gamma')\leq L\}<\infty$,
then the set $\{m^{[\gamma/L]}\}$ is precompact in the space of locally finite measures on $\Curr(\S)$ with the weak$^\star$-topology.
Moreover, if $\liminf_{L\rar\infty}\frac{1}{L^{N(\S)}}\#\{\gamma'\in\Map(\S)\cdot\gamma\,\vert\,\ell_h(\gamma')\leq L\}>0$,
then $(m^{[\gamma/L]})$ does not accumulate at the zero measure.
\end{remark}

Theorem \ref{thm:locally-finite-homogeneous} results in a short proof of the following result from \cite{ES}, which is very different from the original one.

\begin{thm}[Curve-counting measures accumulate at positive multiples of $\Thu$]\label{thm:application}
Let  $\gamma$ be a closed curve on $\S$ and let $(L_i)$ be any sequence of real numbers such that $L_i\to\infty$. 
Then, up to passing to a subsequence, the curve-counting measures $m^{[\gamma/L_i]}$ associated to $\gamma$ converge in the weak$^\star$-topology
to a positive multiple of the Thurston measure $\Thu$.
\end{thm}

Before giving the new proof, we point out the implications of the theorem. First recall a celebrated result by Mirzakhani \cite{Mir08} \cite{mirzakhani-orbit}
implies that, for any hyperbolic metric $h$ on $\S$ 
and any curve $\gamma$, the limit
$$\lim_{L\to\infty}\frac{1}{L^{N(\S)}}\#\{\gamma'\in\Map(\S)\cdot\gamma\,\vert\,\ell_h(\gamma')\leq L\}$$
exists and is positive. This result together with Theorem \ref{thm:application} imply the following.

\begin{cor}[Convergence of curve-counting measures]\label{cor:count}
Let $\gamma$ be a closed curve on $S$. Then there exists $u=u(g,\gamma)>0$ such that 
$$m^{[\gamma/L]}\to u\cdot \Thu$$
as $L\to\infty$. In particular, 
$$\lim_{L\to\infty}\frac{1}{L^{N(\S)}}\#\{\gamma'\in\Map(\S)\cdot\gamma\,\vert\,f(\gamma')\leq L\} = u\cdot\Thu(B_f)$$
for any continuous, $1$-homogeneous function $f:\Curr(S)\to\mathbb{R}_{\geq0}$. 
\end{cor}

Note that if $f$ denotes hyperbolic length, the last assertion of the corollary is only repeating Mirzakhani's result. However, there are many other such functions, including any length coming from a metric on $S$ which has an associated Liouville current such as any negatively curved metric \cite{otal} or Euclidean cone metric \cite{duchin-leininger-rafi:flat} \cite{bankovic-leininger}. In fact, in \cite{EPS} it was shown that $f$ can be replaced with any (possibly singular) Riemannian metric. 

Here we include a proof of the corollary, and we refer to \cite{ES} \cite{EPS} for more details. 

\begin{proof}[Proof of Corollary \ref{cor:count}]
Let $\gamma$ be a closed curve and let $(L_i)$ be a sequence such that $L_i\to\infty$. By Theorem \ref{thm:application}, up to passing to a subsequence, there exists $u>0$ such that $m^{[\gamma/L_i]}\to u\cdot\Thu$. 
Such a $u$ does a priori depend on $g$, $\gamma$ and the subsequence.
Note that, if we can show that $u$ is independent of the subsequence, then $m^{[\gamma/L]}$ converges as $L\to\infty$, and has limit $u\cdot \Thu$. To that end, let $f$ be as above and note that on the one hand 
$$\lim_{i\to\infty}m^{[\gamma/L_i]}(B_f)=u\cdot \Thu(B_f)$$
because $\Thu(\partial B_f)=0$. On the other hand,  
\begin{equation}\label{eq:mirzakhani}
\lim_{i\to\infty}m^{[\gamma/L_i]}(B_f) = \lim_{i\to\infty} \frac{1}{L_i^{N(\S)}}\#\{\gamma'\in\Map(\S)\cdot\gamma\,\vert\,f(\gamma')\leq L_i\}.
\end{equation}
If we consider the particular case when $f=\ell_h$, the length function with respect to a hyperbolic metric $h$ on $\S$, it follows by Mirzakhani's result above that the limit on the right hand side of \eqref{eq:mirzakhani} exists and it is independent of the increasing sequence $(L_i)$. 
Hence, the constant $u$ is independent of the sequence $(L_i)$ too and we can conclude that 
\begin{equation}\label{eq:limit}
\lim_{L\to\infty}m^{[\gamma/L]}\to u\cdot \Thu
\end{equation}
and in particular   
$$ \lim_{L\to\infty} \frac{1}{L^{N(\S)}}\#\{\gamma'\in\Map(\S)\cdot\gamma\,\vert\,f(\gamma')\leq L\} = u\cdot\Thu(B_f)$$
as desired. 
\end{proof}

We point out that Theorem \ref{thm:application} was used by Rafi-Souto \cite{rafi-souto} to prove that Corollary \ref{cor:count} also holds in the case when $\gamma$ is a current. 

\begin{proof}[Proof of Theorem \ref{thm:application}]
We first prove that the family $(m^{[\gamma/L]})$ is precompact, and note that this follows by similar logic as is used in \cite{ES} \cite{rafi-souto}. Fix a hyperbolic metric $h$ on $\S$ and recall that $\{c\in\Curr(S)\,\vert\,\ell_h(c)\leq r\}$ is compact for all $r>0$. 
By Remark \ref{rmk:compactness}, it is enough to show that
%
$$
0<\liminf_{L\to\infty} m^{[\gamma/L]}( B_h )\qquad\text{and}\qquad
\limsup_{L\to\infty} m^{[\gamma/L]}( B_h ) < \infty.$$
To that end, note that 
$$m^{[\gamma/L]}(B_h) = \frac{1}{L^{N(\S)}}\#\{\gamma'\in\Map(\S)\cdot\gamma\,\vert\,\ell_h(\gamma')\leq L\}$$ 
and the quantity on the right is bounded from above and from below by positive constants independent of $L$ by Lemma \ref{lemma:bound-orbit-c}. 

Now, let $m\neq 0$ be any accumulation point of $(m^{[\gamma/L]})$. The above shows that $m$ is locally finite.
Since each $m^{[\gamma/L]}$ is $\Map(\S)$-invariant, so is $m$.
Moreover, it is easy to check that $m$ is $N(\S)$-homogeneous.
As a consequence, $m$ is a positive multiple of $\Thu$ by Corollary \ref{cor:homog-N}.
%
\end{proof}

Not all counting problems will lead to $\Map(\S)$-invariant $d$-homogeneous measures
with $d=N(\S)$, and so proportional to the Thurston measure. Here is an example that came up
after a discussion with Fran\c{c}ois Labourie.

\begin{example}[A counting problem with higher homogeneity]\label{ex:higher}
Let $\pi'$ be a characteristic subgroup of finite index inside $\pi$
(for instance, we can take $\pi'$ to be the kernel of
the homomorphism $\pi\rar H_1(\S;\ZZ/k)$ for any integer $k\geq 2$),
which is thus invariant under the action of $\Map(\S)$.
Such a $\pi'$ corresponds to a finite regular cover $p:\S'\rar\S$.
For every simple closed curve $\gamma'$ in $\S'$, we denote by $p_*\gamma'$ 
the corresponding integral multi-curve in $\S$ (where the concatenation of a curve
$\gamma\subset\S$ with itself $w$ times is identified to the multi-curve $w\gamma$),
and we say that {\it{the multi-curve $p_*\gamma'$ comes from $\S'$}}.

Fix a hyperbolic metric $h$ on $\S$ and let $h'$ be its pull-back on $\S'$.
For every $L>0$, consider the sets
\[
\mathcal{C}:=\{\Gamma\in\Curr(\S)\,|\,\text{$\Gamma$ integral multi-curve that comes from $\S'$}\}
\ \text{and}
\ \mathcal{C}_L:=\{\Gamma\in\mathcal{C}\,|\,\ell_h(\Gamma)\leq L\}
\]
and
\[
\mathcal{C}'_L:=\{\gamma'\ \text{simple closed curve in $\S'$}\,|\,\ell_{h'}(\gamma')\leq L\}.
\]
and define the locally finite measures $m^L:=\frac{1}{L^{N(\S')}}\sum_{\Gamma\in\mathcal{C}}\delta_{\Gamma/L}$
on $\Curr(\S)$. Note that $m^L$ is $\Map(\S)$-invariant because $\pi'$ is a characteristic subgroup of $\pi$.
By \cite{rivin}, we know that $|\mathcal{C}'_L|/L^{N(\S')}$ is bounded above and below by positive constants.
Since the map $p_*:\mathcal{C}'_L\rar\mathcal{C}_L$ is surjective, with fiber of cardinality at most
$[\pi:\pi']$, the quantity $|\mathcal{C}_L|/L^{N(\S')}$ is bounded above and below by positive constants too.
By Remark \ref{rmk:compactness},
there exists a sequence $(L_i)$ with $L_i\rar\infty$ such that
$m^{L_i}\rar m$ and such measure $m\neq 0$ is locally finite and $\Map(\S)$-invariant.
Moreover, it is immediate to see that $m$ is $N(\S')$-homogeneous, with $N(\S')>N(\S)$.
As a consequence, $m$ cannot be a multiple of the Thurston measure.
\end{example}

\appendix

\renewcommand{\thesection}{\Roman{section}}

\section{Estimates}\label{app}

In the present section we collect some estimates that will be employed in the proof
of Proposition \ref{prop:finiteness}, namely a bound on the Thurston volume 
of the $h_\R$-unit ball $B_{h_\R}$ inside $\ML_0(\R)$ (Lemma \ref{lemma:estimate-ball}) and
an asymptotic bound on the number of currents in the $\Map(\S)$-orbit
of $c=\Gamma+\alpha$ of $h$-length at most $L$ (Lemma \ref{lemma:bound-orbit-c}).\\

We begin by introducing some notation.

Consider the subset $\ML^{\ZZ}(\R)\cong\ML_0^{\ZZ}(\R)\times\left(\bigoplus_j \ZZ_{\geq 0}\cdot\pa_j\R\right)$
of integral simple multi-curves in $\R$
inside $\ML(\R)\cong \ML_0(\R)\times\left(\bigoplus_j \RR_{\geq 0}\cdot\pa_j\R\right)$, and define the measures
\[
\ol{m}^{\R/L}:=\frac{1}{L^{N'(\R)}}\sum_{\gamma\in\ML^{\ZZ}(\R)}\delta_{\frac{1}{L}\gamma}\quad
\text{and}\qquad
\bThu^\R:=\Thu^\R\otimes\lambda^{\pa\R}
\]
on $\ML(\R)$, where $\lambda^{\pa\R}$ is the Lebesgue measure
on $\bigoplus_j \RR_{\geq 0}\cdot\pa_j\R$.
Since $\bigoplus_j \ZZ_{\geq 0}\cdot\pa_j\R$ is a lattice of unit co-volume in $\bigoplus_j \RR_{\geq 0}\cdot\pa_j\R$,
we have $\ol{m}^{\R/L}\rar \bThu^\R$ in the weak$^\star$-topology.

We denote by $b_h^c(L)$ the number of currents in the $\Map(\S)$-orbit of $c$
of $h$-length at most $L$
and by $b'_{h_\R}(L)$ the number of points in $\ML^{\ZZ}(\R)$
of $h_\R$-length at most $L$, so that 
\[
b'_{h_\R}(L)=L^{N'(\R)}\ol{m}^{\R/L}(B_{h_\R}).
\]
Similarly, we denote by $v_{h_\R}(L)$ the volume of the subset of laminations in $\ML(\R)$
of $h_\R$-length at most $L$, so that
\[
v_{h_\R}(L)=L^{N'(\R)}\bThu^\R(B_{h_\R}).
\]

%
\subsection{The proof of Lemma \ref{lemma:estimate-ball}}

We proceed in three steps:
we relate first $\Thu^\R(B_{h_\R})$ to $v_{h_\R}(1)$, then
$v_{h_\R}(1)$ to $b'_{h_\R}(L)/L^{N'(\R)}$, and finally
we estimate $b'_{h_\R}(L)/L^{N'(\R)}$ in terms of the $h_\R$-lengths
of the boundary components of $\R$ and of the $h_\R$-systole of $\R$.

In the following lemma we relate $v_{h_\R}(L)$ and $\Thu^\R(B_{h_\R})$.

\begin{lemma}[Volume of balls in $\ML(\R)$ and in $\ML_0(\R)$]\label{lemma:volume}
For every hyperbolic surface $(\R,h_\R)$ with geodesic boundary,
the following holds
\[
v_{h_{\R}}(L)= \Thu^\R(B_{h_\R}) \frac{L^{N'(\R)}}{N'(\R)\cdot(n-1)!\,\prod_j \ell_{h_\R}(\pa_j \R)}.
\] 
where $\ML_0(\R)$ is endowed with the Thurston measure and $\RR_{\geq 0}\pa_j\R$ with the Lebesgue measure.
\end{lemma}
\begin{proof}
Since every lamination in $\R$
can be written as $\lambda+\sum_j x_j\pa_j\R$ for a certain $\lambda$ supported in the interior of $\R$ and $x_i\in\RR$,
the volume $v_{h_{\R}}(L)$ can be computed as
\[
v_{h_{\R}}(L) =\int_{\hat{D}} \Thu^\R(B_h) \cdot \big(L-\sum_j x_j\ell_{h_\R}(\pa_j\R)\big)^{N(\R)}\,dx_1\cdots dx_n
\]
where the integration is performed over the domain $\hat{D}=\{(x_1,\dots,x_n)\in\RR^n_{\geq 0}\,|\,L-\sum_j x_j \ell_{h_\R}(\pa_j\R)\geq 0\}$.
By the change of variables $t_j:=x_j \ell_{h_\R}(\pa_j\R)/L$, we obtain
\[
v_{h_{\R}}(L)= \Thu^\R(B_h) \frac{L^{N'(\R)}}{\prod_j \ell_{h_\R}(\pa_j\R)} \int_D \big(1-\sum_j t_j\big)^{N(\R)}\,dt_1\cdots dt_n
\]
where the integration is performed over the domain $D=\{(t_1,\dots,t_j)\in\RR^n_{\geq 0}\,|\,\sum_j t_j\leq 1\}$.
One can easily check that $\int_D (1-\sum_j t_j)^{N(\R)}dt_1\cdots dt_n=\frac{1}{(N(\R)+n)(n-1)!}$.
\end{proof}

We will also need that the asymptotic $b'_h(L)\sim v_h(L)$ as $L\rar+\infty$
is uniform over all hyperbolic metrics $h$ whose systole is bounded from below.

\begin{lemma}[Volume of balls in $\ML(\R)$ and simple integral multi-curves]\label{lemma:simple-volume}
Let $\R$ be a hyperbolic surface with geodesic boundary and let $s>0$.
Given $\e>0$ there exists $L_0>0$ (that may depend on $\e$ and $s$) such that
\[
\left|\frac{b'_{h_\R}(L)}{v_{h_\R}(L)}-1\right|<\e
\]
for all $L\geq L_0$ and for all metrics $h_\R$ on $\R$ with $\sys(h_\R)\geq s$.
\end{lemma}
\begin{proof}
Note that all involved quantities are invariant under action of $\Map(\R)$.
Thus it is enough to consider $h$ in a fundamental domain $\mathcal{F}(\R)$
for the action of $\Map(\R)$ on the Teichm\"uller space of hyperbolic metrics on $\R$.
We denote by $\mathcal{F}^s(\R)$ the subset of metrics $h_\R$ on $\R$
with $\sys(h_\R)\geq s$, which is well-known to be compact.

%

Since the length function associated to $h_\R$ depends continuously on $h_\R$
and $\mathcal{F}^s(\R)$ is compact, the union
$B^s$ of the balls $B_{h_\R}\subset\ML(\R)$ as $h_\R$ ranges in $\mathcal{F}^s(\R)$
is a compact subset. We also denote by $v^s>0$ the minimum of $v_{h_\R}(1)=L^{-N'(\R)}v_{h_\R}(L)$ as $h_\R$ ranges in $\mathcal{F}^s(\R)$.

Since $\ol{m}^{\R/L}\rar \bThu^\R$
in the weak$^\star$-topology of $\ML(\R)$, there exists $L_0$ such that
$|\ol{m}^{\R/L}-\bThu^\R|(B^s)<v^s\e$ for all $L>L_0$.
Thus,
\begin{align*}
L^{-N'(\R)}|b'_{h_\R}(L)-v_{h_\R}(L)| &=|\ol{m}^{\R/L}(B_{h_\R})-\bThu^\R(B_{h_\R})|\leq |\ol{m}^{\R/L}-\bThu^\R|(B_{h_\R})
\leq \\
&\leq |\ol{m}^{\R/L}-\bThu^\R|(B^s)<v^s\e
\end{align*}
for all $L>L_0$, and we conclude that
\[
\left|\frac{b'_{h_\R}(L)}{v_{h_\R}(L)}-1\right|=\frac{|b'_{h_\R}(L)-v_{h_\R}(L)|}{L^{N'(\R)}} \frac{L^{N'(\R)}}{v_{h_\R}(L)}
\leq \frac{v^s\e}{v_{h_\R}(1)}\leq \e
\]
for all $L>L_0$.
\end{proof}

The last estimate needed
to prove Lemma \ref{lemma:estimate-ball} 
concerns $b'_{h_\R}(L)$ and it is a minor variation
of Proposition 3.6 in \cite{Mir08}.

\begin{lemma}[Simple integral multi-curves and boundary lengths]\label{lemma:estimate-b'}
Let $(\R,h_\R)$ be a hyperbolic surface with $\sys(h_\R)\geq s>0$.
Then
\[
\frac{k_4}{\ell_{h_\R}(\pa\R)^{N'(\R)}}\leq
\frac{b'_{h_{\R}}(L)}{L^{N'(\R)}}\leq \frac{k_5}{\prod_j\ell_{h_\R}(\pa_j\R)}
\]
where $k_4,k_5$ are constants that depends only on the topology of $\R$
and on $s$.
\end{lemma}
\begin{proof}
Let $\mathcal{P}=\{\eta_1,\dots,\eta_{N'(\R)}\}$ be a maximal
set of simple, pairwise disjoint geodesic arcs in $\R$ that meet $\pa\R$ orthogonally.
Since $\sys({\R})\geq s>0$, such arcs $\eta_i$
can always be chosen to be shorter than a constant that only
depends on $s$.
Up to relabeling, we can assume that $\ell_{h_\R}(\eta_1)\leq \ell_{h_\R}(\eta_2)\leq\dots\leq\ell_{h_\R}(\eta_{N'(\R)})$.

The trivalent ribbon graph embedded in $\R$
dual to $\mathcal{P}$ has set $V$ of vertices 
corresponding to components of
$\R\setminus\bigcup_i\eta_i$
and set $E=\{1,\dots,N'(\R)\}$ of edges corresponding to 
the $\eta_i$'s.

For every integral simple multi-curve $\gamma\in\ML^{\ZZ}(\R)$, let
$DT_j(\gamma):=\iota(\gamma,\eta_j)$ 
for all $j\in E$
and $DT(\gamma)=(DT_1(\gamma),\dots,DT_{N'(\R)}(\gamma))\in\NN^E$.
For every $v\in V$,
let $E_v$ be the subset of indices $\{i_1,i_2,i_3\}\subset E$ such that
$v$ is bounded by the arcs $\eta_{i_1},\eta_{i_2},\eta_{i_3}$
and denote by $DT_v(\gamma)$ the sum $DT_{i_1}(\gamma)+DT_{i_2}(\gamma)+DT_{i_3}(\gamma)$.

This easier version of Dehn-Thurston coordinates establishes a bijection
\[
DT:
\ML^{\ZZ}(\R)
\lra \Zcal=
\left\{m\in\NN^{E}
\ \Big|\ 
\begin{array}{c}
\text{$m_v\geq 2m_i$ is even}\\
\text{for every $v\in V$ and $i\in E_v$}
\end{array}
\right\}
\]
and we define $\ell_{\mathcal{P}}(m):=\sum_i m_i \Coll(\ell_{h_\R}(\eta_i))$
for all $m\in \NN^E$, where $\Coll$
is the decreasing function $\Coll(\ell):=\arcsinh(\sinh(\ell/2)^{-1})$.

We notice that 
\begin{itemize}
\item
$\Zcal_M\subset \Zcal\subset\NN^E$ for every $M>0$,
where
$\Zcal_M:=(2\mathbb{N}\cap [2M,3M])^E$
\item
by Proposition 3.5 in \cite{Mir08}, the following estimate holds
\[
\frac{1}{k}\leq \frac{\ell_{h_\R}(\gamma)}{\ell_{\mathcal{P}}(DT(\gamma))}\leq k
\]
for a constant $k>1$ that depends only on $\R$ and on $s$.
\end{itemize}
Since $\iota(\pa\R,\eta_i)=2$
and $\Coll(\ell_{h_\R}(\eta_1))\geq \Coll(\ell_{h_\R}(\eta_i))$ for all $i$, we in particular obtain
\[
\frac{2}{k}\Coll(\ell_{h_\R}(\eta_1))
\leq \ell_{h_\R}(\pa\R)\leq
2kN'(\R)\Coll(\ell_{h_\R}(\eta_1)).
\]
By the above considerations
\[
\Zcal_M(L/k)\subseteq DT(\ML^\ZZ(\R)(L))\subseteq \Zcal(kL)
\]
where $\Zcal(kL)=\{m\in\Zcal\,|\,\sum_i m_i\Coll(\ell_{h_\R}(\eta_i))\leq kL\}$
and $\Zcal_M(L/k)=\{m\in\Zcal_M\,|\,\sum_i m_i\Coll(\ell_{h_\R}(\eta_i))\leq L/k\}$.

On one hand, we observe that $\Zcal(kL)\subset \NN^E\cap \prod_i [0,\frac{kL}{\Coll(\ell_{h_\R}(\eta_i))}]$.
Since there exists a constant $k'=k'(s)>0$
such that $k'<\Coll(\ell_{h_\R}(\eta_i))<\ell_{h_\R}(\pa_j\R)$ whenever the arc $\eta_i$ meets $\pa_j\R$,
we obtain
\[
|\Zcal(kL)|\leq \frac{k'_5 L^{N'(\R)}}{\prod_i \Coll(\ell_{h_\R}(\eta_i))}
\leq k_5\frac{L^{N'(\R)}}{\prod_j\ell_{h_\R}(\pa_j\R)}
\]
where $k'_5$ and $k_5$ depend on $\R$ and $s$ only
and the last product is taken over all the boundary components $\pa_j \R$ of $\R$.

On the other hand, we can take $M=\lfloor \frac{L}{3kN'(\R)\Coll(\ell_{h_\R}(\eta_1))}\rfloor$ in such a way that
$\Zcal_M(L/k)\supset (2\NN\cap [2M,3M])^E$. It follows that
\[
|\Zcal_M(L/k)|\geq (M/2-1)^{N'(\R)}\geq k_4\frac{L^{N'(\R)}}{\ell_{h_\R}(\pa\R)^{N'(\R)}}.
\]
We conclude that there are constants $k_4,k_5>0$ such that
\[
\frac{k_4}{\ell_{h_\R}(\pa\R)^{N'(\R)}}\leq
\frac{b'_{h_{\R}}(L)}{L^{N'(\R)}}\leq \frac{k_5}{\prod_j\ell_{h_\R}(\pa_j\R)}
\]
with $k_4,k_5$ as desired.
\end{proof}

We have now all the ingredients to estimate the volume of the unit ball $B_{h_{\R}}$.

\begin{named}{Lemma \ref{lemma:estimate-ball}}[Volume of the unit ball in $\ML_0(\R)$]
Let $(\R,h_\R)$ be a hyperbolic surface with geodesic boundary with 
systole $\sys(\R)\geq s>0$.
Then
\[
\frac{\hat{k}}{\ell_{h_\R}(\pa\R)^{N(\R)}}
<\Thu^\R(B_{h_\R})<\hat{k}'.
\]
for suitable constants $\hat{k},\hat{k}'>0$ that depend only on $s$ and
the topology of $\R$.
\end{named}
\begin{proof} 
From Lemma \ref{lemma:volume} we have
\[
\Thu^\R(B_{h_\R})=N'(\R)\cdot(n-1)!\prod_j \ell_{h_\R}(\pa_j \R)\cdot \frac{v_{h_\R}(L)}{L^{N'(\R)}}
\]
for any $L>0$. As $s>0$ is fixed, 
we can find a suitably large $L$
such that Lemma \ref{lemma:simple-volume} gives
$\frac{1}{2}b'_{h_\R}(L)<v_{h_\R}(L)<2b'_{h_\R}(L)$, and so
\[
\frac{1}{2}N'(\R)\cdot(n-1)!\prod_j \ell_{h_\R}(\pa_j \R)\cdot \frac{b'_{h_\R}(L)}{L^{N'(\R)}}\leq
\Thu^\R(B_{h_\R})\leq 2N'(\R)\cdot(n-1)!\prod_j \ell_{h_\R}(\pa_j \R)\cdot \frac{b'_{h_\R}(L)}{L^{N'(\R)}}.
\]
Together with Lemma \ref{lemma:estimate-b'} we obtain the wished estimate.
\end{proof}

\subsection{The proof of Lemma \ref{lemma:bound-orbit-c}}

We recall the following lemma by Ivanov.

\begin{lemma}[Dehn twists and intersection numbers {\cite[Lemma 4.2]{ivanov:subgroups}}]\label{lemma:ivanov}
Let $\{\eta_j\}_{j=1}^p$ be disjoint simple closed curves on the surface $\S$
and let $\{t_j\}$ be integers. Then
\[
-\iota(\gamma,\beta)+\sum_{j=1}^p \left(|t_j|-2\right)\iota(\gamma,\eta_j)\iota(\eta_j,\beta)
\leq \iota(\psi(\gamma),\beta)
\]
where $\psi=T^{t_1}_{\eta_i}\circ\dots\circ T^{t_p}_{\eta_p}$ is the composition of $t_j$ right Dehn twists about each $\eta_j$ and
$\beta,\gamma$ are measured laminations.
\end{lemma}

Actually, we will employ Ivanov's above inequality in the following form.

\begin{cor}[Dehn twists and hyperbolic lengths]\label{cor:ivanov}
Let $(\S,h)$ be a hyperbolic surface with $\sys(\S)\geq s>0$ and
let $\{\eta_j\}$ be a pair of pants decomposition.
Let $\gamma=\gamma_1+\dots+\gamma_p$ be a finite sum of simple closed curves
(but we do not exclude that $\iota(\gamma,\gamma)>0$)
such that 
\begin{itemize}
\item[(a)]
each $\gamma_j$ intersects $\eta_j$ in one or two points for $j=1,\dots,p$,
\item[(b)]
$\iota(\gamma_j,\eta_i)=0$ for $i\neq j$.
\end{itemize}
Moreover, let $c=\gamma+\eta_1+\dots+\eta_q$ with $q\geq p$.\\
Then there exists $\tilde{k}>0$ that depends only on $h$ such that the following holds:
\begin{itemize}
\item
for every $\psi=T^{t_1}_{\eta_i}\circ\dots\circ T^{t_p}_{\eta_p}$ composition of $t_j$ Dehn twists about $\eta_j$ for each $j=1,\dots,p$ such that $\ell_h(c),\ell_h(\psi(c))\leq L$, we have
$\sum_{j=1}^p |t_j|\ell_h(\eta_j)\leq \tilde{k}L$ and,
in particular, $\dis |t_j|\leq \tilde{k}\frac{L}{\ell_h(\eta_j)}$ for all $j=1,\dots,p$.
\end{itemize}
\end{cor}
\begin{proof}
%
Note that the formula in Lemma \ref{lemma:ivanov} is linear both in $\beta$ and in $\gamma$.
We choose $\beta$ to be a binding current in $\S$, which is a finite sum of simple closed curves,
and we choose $\gamma$ as in the statement. 
By Lemma \ref{lemma:ivanov}, we have
\[
-\iota(\beta,\gamma)+\sum_{j=1}^p\left(|t_j|-2\right)\iota(\gamma,\eta_j)\iota(\beta,\eta_j)\leq \iota(\beta,\psi(\gamma))
\]
and so
\[
\sum_{j=1}^p |t_j|\iota(\beta,\eta_j)\leq \iota(\beta,\gamma)+\iota(\beta,\psi(\gamma))+4\sum_{j=1}^p\iota(\beta,\eta_j).
\]
Recall that there exists a compact subset $\cpt$ of the interior of $\S$ that 
depends only on $h$
and that contains all simple closed $h$-geodesics, and so in particular it contains $\psi(\gamma)$ for all $\psi\in\Map(\S)$.
By Lemma \ref{lemma:bound}, there is a constant $r>0$ that depends on $h$
such that $\frac{1}{r}<\frac{\iota(\beta,\cdot)}{\ell_h}<r$.
It follows that
\[
\frac{1}{r}\sum_{j=1}^p |t_j|\ell_h(\eta_j) \leq r \ell_h(\gamma)+r\ell_h(\psi(\gamma))+4r\sum_{j=1}^p\ell_h(\eta_j)
\leq 5rL
\]
because $\gamma+\sum_{j=1}^p\eta_j\leq c$ and $\psi(\gamma)\leq \psi(c)$ as currents.
It follows that $\sum_{j=1}^p |t_j|\ell_h(\eta_j)\leq 5r^2 L$ and so we can take $\tilde{k}=5r^2$.
\end{proof}

Finally, we will also need the following observation.

\begin{remark}\label{rmk:transverse}
Let $(\S,h)$ be a hyperbolic surface with $\sys(h)\geq s>0$.
Let $\{\eta_j\}$ be a pair of pants decomposition of $\S$
and let $\gamma_l$ be a simple closed curve on $\S$ that intersects $\eta_l$ in one or two points
and such that $\iota(\gamma_l,\eta_j)=0$ for $j\neq l$.
There exists a constant $a=a(s)\geq 1$ that depends on $s$ only and an integer $t$ such that 
$\ell_h(T^t_{\eta_l}\gamma_l)\leq a\cdot \max_j\{\ell_h(\eta_j)\}$.
\end{remark}
\begin{proof}[Idea of the proof of Remark \ref{rmk:transverse}]
Consider the complement $\S_l$ of $\bigcup_{j\neq l}\eta_j$ inside $\S$, which can 
be a torus with one boundary component, or 
a sphere with four boundary components. By elementary trigonometry it is easy to see
that there exists an integer $t$ such that
$\ell_h(T^t_{\eta_l}\gamma_l)\leq a\cdot \ell_l$, where $\ell_l$ is the maximum
$h$-length of a boundary component of $\S_l$.
\end{proof}

The key idea needed to prove Lemma \ref{lemma:bound-orbit-c}
is to relate the number of currents in the orbit of $c$ to 
the number of suitable integral simple multi-curves.
The following is essentially borrowed from Lemma 5.6 in \cite{mirzakhani-orbit}.

\begin{lemma}[Comparing orbits of a current with orbits of simple multi-curves]\label{lemma:c-nonsimple}
Let $(\S,h)$ be a hyperbolic surface (possibly with boundary)
and let $c$ be a current of type $c=\Gamma+\alpha$.
There exists a constant $k=k(\S,s,c)>1$
and an integral simple multi-curve $\Gamma'$
such that
\[
b_h^{\Gamma'}(L/k)\leq b_h^c(L)\leq b'_h(kL)
\]
for all $L$.
\end{lemma}
\begin{proof}
Let $\mathcal{P}=\{\eta_1,\dots,\eta_{N(\S)+n}\}$ be
a pair of pants decomposition of $\S$ such that
\begin{itemize}
\item
$\eta_1\cup\dots\cup\eta_p$ sits inside $\A:=\hull(\alpha)$ (note that $p=-\frac{1}{2}\left(3\chi(\A)+n(\A)\right)$,
where $n(\A)$ is the number of boundary components of $\A$),
\item
$\eta_{p+1}\cup\dots\cup\eta_{p'}$ is equal to $\supp(\Gamma)\cap\pa\A$, and
\item
$\eta_{p'+1}\cup\dots\cup\eta_q$ is equal to $\supp(\Gamma)\setminus\pa\A$.
\end{itemize}

{\it{Reduction to the case of $c$ equal to a multi-curve.}}
Note that there exists simple closed curves $\gamma_1,\dots,\gamma_p$ in the interior of $\A$
 such that
\begin{itemize}
\item
$\gamma=\gamma_1+\dots+\gamma_p$ satisfies hypotheses (a) and (b) in Corollary \ref{cor:ivanov};
\item
$\eta_1\cup\dots\eta_p\cup\gamma$ binds $\A$.
\end{itemize}

Denote by $c_\A$ the sum of $\alpha$ and of the summands of $c$ supported on the components of $\pa\A$.
Since the systole of $h$ is bounded from below by $s$, 
for all $\varphi\in\Map(\S)$ the ratio
\[
\frac{\ell_h(\varphi(c_\A))}{\ell_h(\varphi(\gamma+\sum_{j=1}^{p'}\eta_j))}=
\frac{\ell_{\varphi^*h}(c_\A)}{\ell_{\varphi^*h}(\gamma+\sum_{j=1}^{p'}\eta_j)} 
\]
is bounded below and above by positive constants that depend on
$\A$, $s$ and $\alpha$ only by Remark \ref{rmk:binding-comparable}.
Moreover, $c-c_\A$ is the sum of positive multiples of $\eta_{p'+1},\dots,\eta_q$ and so
$\ell_h(\varphi(c-c_\A))/\ell_h(\sum_{j={p'+1}}^q \varphi(\eta_j))$ is also bounded below and above by positive constants that only depend on $c$.
As a consequence, for all $\varphi\in\Map(\S)$ the ratio
\[
\frac{\ell_h(\varphi(c))}{\ell_h(\varphi(\gamma+\sum_{j=1}^q\eta_j))}
\]
is bounded above and below by constants that depend on $\S$, $s$ and $c$ only.

Thus, we can assume that $c=\gamma+\sum_{j=1}^q\eta_j$
and we choose $\Gamma':=\sum_{j=1}^q\eta_j$.

Let  $G:=\bigcap_{j=1}^q \stab(\eta_j)$ and $H:=\Map(\S,\A)\cap G$ and note
that $H$ has finite index in $\stab(c)$, because $\gamma+\sum_{j=1}^p\eta_j$ fills $\A$.
In order to give upper and lower bounds for the cardinality $b_h^c(L)$ of the set
\[
\{\varphi\in\Map(\S)/\stab(c)\,|\,\ell_h(\varphi(c))\leq L\}
\]
we follow Mirzakhani's idea and we construct a map 
\[
\xi:\{\varphi\in\Map(\S)/H\,|\,\ell_h(\varphi(c))\leq L\}\lra \ML^{\ZZ}(\S)
\]
with fibers of bounded cardinality.
The wished conclusion follows as we show that
there exist constants $k_1,k_2>1$ such that 
\begin{itemize}
\item[(i)]
the image of $\xi$
is contained
inside the subset of $\ML^{\ZZ}(\S)$
consisting of multi-curves with $\ell_h\leq k_1 L$, and
\item[(ii)]
the image of $\xi$ 
contains all points in $\Map(\S)\cdot\Gamma'$
with $\ell_h\leq L/k_2$.
\end{itemize}


{\it{Definition of the map $\xi$.}}
Note that $G/H$ is isomorphic to the free
Abelian group $\Psi$ of diffeomorphisms of type $\psi=T^{t_1}_{\eta_1}\circ\dots\circ T^{t_p}_{\eta_p}$
for suitable $t_1,\dots,t_p\in\ZZ$.
In every $[\varphi]\in\Map(\S)/G$ we choose a representative $\varphi^0\in\Map(\S)/H$ as follows.
Pick any representative $\varphi\in\Map(\S)/H$ in the class $[\varphi]$.
The map $\Psi\ni \psi\mapsto \ell_{\varphi^*h}(\psi(c))\in\RR_+$ is proper, because $\gamma+\eta_1+\dots+\eta_p$ fills $\A$,
and so it achieves a minimum at some element $\psi_\varphi$; we define $\varphi^0:=\varphi\circ\psi_\varphi$.
Now, every $\varphi'\in\Map(\S)/H$ in the class $[\varphi]$ can be uniquely written as 
$\varphi'=\varphi^0\circ T^{t_1(\varphi')}_{\eta_1}\circ\dots\circ T^{t_p(\varphi')}_{\eta_p}$ and we let $w_j(\varphi'):=|t_j(\varphi')|$.
We then define $\xi$ as
\[
\varphi\longmapsto\xi_\varphi:=\sum_{j=1}^p w_j(\varphi)\cdot\varphi(\eta_j)+
\sum_{j=1}^q \varphi(\eta_j).
\]
Note that $\xi_\varphi$  determines an element $\varphi\in\Map(\S)/H$ up to finitely many choices (namely, up to
classes of diffeomorphisms that possibly exchange the signs of the $t_j(\varphi)$ and permute the $\eta_j$'s).
Hence, each fiber of $\xi$ has cardinality bounded above by a constant that depends only on the topology of $\S$.

{\it{Upper bound (i).}}
We want to show that $\ell_h(\xi_\varphi)\leq k_1 L$ with $k_1:=\tilde{k}p+1$.

Since $\ell_h(\varphi^0(c))\leq \ell_h(\varphi(c))\leq L$, we have
$w_j(\varphi)=|t_j(\varphi)|\leq \tilde{k} \frac{L}{\ell_h(\varphi(\eta_j))}$
by Corollary \ref{cor:ivanov}.
Thus,
$\ell_h(\xi_\varphi)\leq \sum_{j=1}^p (\tilde{k}L)+\sum_{j=1}^q\ell_h(\varphi(\eta_j))$.
Since $\sum_{j=1}^q\ell_h(\varphi(\eta_j))\leq L$, it follows that $\ell_h(\xi_\varphi)\leq (\tilde{k}p+1) L=k_1 L$.


{\it{Lower bound (ii).}}
Let $a=a(s)\geq 1$ be the constant that appears in Remark \ref{rmk:transverse}
and take $k_2=1+ap$.
Fix $\varphi\in\Map(\S)/H$ such that $\ell_h(\varphi^0(\Gamma'))=\ell_h(\varphi(\Gamma'))\leq L/k_2$.
We want to show that $\ell_h(\varphi^0(c))\leq L$ and so $\varphi(\Gamma')=\xi_{\varphi^0}$ belongs to the image of $\xi$.

For every $j=1,\dots,p$ by hypothesis
$\ell_h(\varphi(\eta_j))\leq L/k_2$,
and so Remark \ref{rmk:transverse} implies that there exists $t_j\in\ZZ$
such that $\ell_h(\varphi(T^{t_j}_{\eta_j}\gamma_j))\leq a(L/k_2)$.
Hence, $\ell_h(\varphi\circ\psi(\gamma_j))\leq a(L/k_2)$ with $\psi=T^{t_1}_{\eta_1}\circ\dots\circ T^{t_p}_{\eta_p}$,
and so $\ell_h(\varphi\circ\psi(c))\leq (1+ap)(L/k_2)$ because $c=\Gamma'+\gamma_1+\dots+\gamma_p$.
As a consequence, $\ell_h(\varphi^0(c))\leq \ell_h(\varphi\circ\psi(c))\leq (1+ap)(L/k_2)=L$.
\end{proof}

Using the bounds proven in Lemma \ref{lemma:c-nonsimple}, we can now
obtain the wished estimate for $b_h^c$.

\begin{proof}[Proof of Lemma \ref{lemma:bound-orbit-c}]
It follows from \cite{Mir08}  
that the quantities $b'_h(L)/L^{N(\S)}$ 
and $b_h^{\Gamma'}(L)/L^{N(\S)}$ are bounded from below and above
by constants that depend on $\S$, $h$ and $c$.
Thus, by Lemma \ref{lemma:c-nonsimple}
there exists $\tilde{v}>1$ such that $(1/\tilde{v})L^{N(\S)}<b_h^c(L)<\tilde{v}L^{N(\S)}$.
In particular,
\[
\left(q^{N(\S)}/\tilde{v}-\tilde{v}\right)L^{N(\S)}<b_h^c([L,qL])<\left(\tilde{v} q^{N(\S)}\right)L^{N(\S)}
\]
Let $q>1$ be large enough so that $v:=\tilde{v} q^{N(\S)}$
satisfies $(1/v)<(q^{N(\S)}/\tilde{v})-\tilde{v}$.
It follows that $(1/v)L^{N(\S)}<b_h^c([L,qL])<vL^{N(\S)}$
and, clearly, $(1/v)L^{N(\S)}<b_h^c(L)<vL^{N(\S)}$.
\end{proof}

%
%
%


\bibliographystyle{amsalpha} 
\bibliography{ref-rev}

\end{document}